\documentclass[12pt]{amsart}
\usepackage{a4}
\usepackage{amsfonts}
\usepackage{amssymb}
\usepackage[latin1]{inputenc}

\author{François Courtès} 
\title[Steinberg representations for groups of parahoric types]{Steinberg representations for groups of parahoric types: the special case}
\address{Universit\'e de Poitiers\\D\'epartement de Math\'ematiques\\UMR 6086 du CNRS\\T\'el\'eport 2\\Boulevard Marie et Pierre Curie\\86962 Futuroscope Chasseneuil Cedex}
\email{courtes@math.univ-poitiers.fr}
\begin{document}
\def \mth{\mathbb}
\begin{abstract}
In this paper, we define and study a kind of Steinberg representation for linear algebraic groups of a particular kind, called groups of parahoric type, defined overa finite field; in particular, when $G$ is the group of $F$-points of a connected reductive quasisplit group defined over $F$ which splits over an unramified extension of $F$, the quotients of parahoric subgroups of $G$ by their congruence subgroups are groups of parahoric type. In particular, under certain conditions on the residual characteristic $p$ of $F$, we determine the irreducible factors of the Steinberg representation of a group ${\mth{G}}$ of parahoric type associated to a pseudo-Borel subgroup of ${\mth{G}}$ in the case when ${\mth{G}}$ is special, that is a quotient of a maximal special parahoric subgroup of $G$.
\end{abstract}
\maketitle
\vskip 1cm
\font\teuf=eufm10
\font\seuf=eufm7
\font\sseuf=eufm6
\newfam\euffam
\textfont\euffam=\teuf
\scriptfont\euffam=\seuf
\scriptscriptfont\euffam=\sseuf
\def \got{\fam\euffam}
\newtheorem{theo}{Theorem}[section]
\newtheorem{prop}[theo]{Proposition}
\newtheorem{lemme}[theo]{Lemma}
\newtheorem{cor}[theo]{Corollary}
\newtheorem{cond}{Condition}

\section{Introduction}

This paper is devoted to the extension of the notion of Steinberg representation to groups defined over a finite ring, namely the quotient of the ring of integers of a $p$-adic field by a non-maximal ideal. The Steinberg representation was first introduced by Steinberg in the context of finite groups of Lie type (see \cite{st1} and \cite{st2} for classical groups, \cite{st3} and \cite{st4} for twisted groups); he alspo proves its irreducibility. Later on, Matsumoto (\cite{mat}) and Shalika (\cite{shal}) studied a "special" irreduciblie admissible representation for reductive $p$-adic groups, and it was soon pointed out (see for example \cite[5.10]{bos}) that it was in fact the $p$-adic equivalent of the Steinberg representation; it thus took the same name. 

We'll proceed the other way round here, starting from the $p$-adic representation. Let $F$ be a local complete field with discrete valuation; let $\mathcal{O}$ be its ring of integers, ${\mathfrak{p}}$ the maximal ideal of $\mathcal{O}$, $k=\mathcal{O}/{\mathfrak{p}}$ its residual field; we'll assume $k$ is finite, and let $p$ be its characteristic. Let $F_{nr}$ be the maximal unramified extension of $F$.

Let $\underline{G}$ be a connected reductive algebraic group defined over $F$; we'll assume $\underline{G}$ is quasisplit and splits over an unramified extension of $F$.
Let $G=\underline{G}(F)$ (resp. $G_{nr}=\underline{G}(F_{nr})$) be the group of $F$-points (resp. $F_{nr}$-points) of $\underline{G}$.

The Steinberg representation $St$ of the group $G$ can be defined as follows: let's choose a maximal torus $T$ of $G$.and a Borel subgroup $B$ of $G$. For every parabolic subgroup $P$ of $G$ containing $B$, let $1_P$ be the induced representation from $P$ fo $G$ of the trivial character on a Levi component $M$ of $P$. Obviously we have $1_P\supset 1_{P'}$ for every $P\subset P'$, and $St=St_B$ is the quotient of $1_B$ by the sum of the $1_P$, $P\subsetneq B$. Up to equivalence, the representation $St$ doesn't depend on the choice of $B$; moreover, as previously stated, $St$ is irreducible. Replacing $B$ by a parabolic subgroup $P$ of $G$, we can similarly define the generalized Steinberg representation $St_P$ of $G$, which depends only on the conjugacy class of $P$, but which is in general not irreducible.

Let $H$ be any subgroup of $G$, let $H'$ be a normal subgroup of $H$ and let $L$ be the quotient group $H/H'$. Let $\pi$ be a representation of $G$; the group $H$ acts on the set of $H'$-fixed vectors of $\pi$, and this can be viewed as a representation of $L$, which we will denote by $\pi^L$.

Now let $K$ be a parahoric subgroup of $G$, let $h$ be a positive integer and let $K^h$ be its $h$-th congruence subgroup; this is a normal subgroup of $K$. Assume first $h=1$: the quotient ${\mth{G}}=K/K^1$ is the group of $k$-points of a reductive group $\underline{\mth{G}}$ defined over the residual field $k$ of $F$. Let $T$ be a maximal torud of $G$ whose parahoric subgroup is contained in $K$, and let $B$ be a Borel subgroup of $G$ containing $T$; the image ${\mth{B}}$ of $B\cap K$ in $G$ is a Borel subgroup of ${\mth{G}}$. We can define the Steinberg representation $St_{\mth{B}}$ of ${\mth{G}}$ the same way as for $G$; it is irreducible as well; Moreover, it is easy to check that $(1_B)^{\mth{G}}=1_{\mth{B}}$, we have a similar equality for any parabolic subgroup of $G$ containing $B$, from wich we deduce that $(St)^{\mth{G}}=St_{\mth{B}}$. Similarly, if ${\mth{P}}$ is a parahoric subgroup of ${\mth{G}}$, the generalized Steinberg representation $St_{\mth{P}}$ of ${\mth{G}}$ is the reduction to ${\mth{G}}$ of a generalized Steinberg representation of $G$. See for example \cite{sam} for a study of such representations.

Assume now $h>1$. The group ${\mth{G}}$ is still the group of $k$-points of a connected linear algebraic group $\underline{\mth{G}}$ defined over $k$, but this group is not reductive anympre; its unipotent radical is $K^1/K^h$. It anyway has some properties which make it look like a kind of reductive group; we'll call an algebraic group having such properties a group of parahoric type.
Moreover, the image ${\mth{B}}$ (resp ${\mth{P}}$) in ${\mth{G}}$ of the intersection with $K$ of $B$ (resp. $P$) is not a Borel (resp. parabolic) subgroup of $G$ but also looks quite similar to such a subgroup; we'll call it a pseudo-Borel (resp. pseudo-parabolic) subgroup of ${\mth{G}}$.

We can define the Steinberg representation $St_{\mth{B}}$ of ${\mth{G}}$ relatively to ${\mth{B}}$ in a similar way as before. Note that contrary to Borel subgroups, two pseudo-Borels are not necessarily conjugate to each other, hence the representation $St_{\mth{B}}$ may depend on ${\mth{B}}$.

Moreover, this representation is not irreducible. One first reason for that is that the pseudo-parabolic subgroups of ${\mth{G}}$ are not the only subgroups of ${\mth{G}}$ containing a Borel subgroup. That's why we have to consider the full set of such subgroups, which we will call generalized pseudo-parabolic subgroups. We can then define the (small) Steinberg representation $st_{\mth{B}}$ of ${\mth{G}}$ as the quotient of $1_{\mth{B}}$ by all the $1_{\mth{P}}$, where ${\mth{P}}$ runs through the set of generalized pseudo-parabolic subgroups of ${\mth{G}}$ strictly containing ${\mth{B}}$. For every generalized pseudo-parabolic subgroup ${\mth{P}}$ of ${\mth{G}}$, we can also define the generalized small Steinberg representation $st_{\mth{P}}$ of ${\mth{G}}$ in a similar fashion.

Note that in \cite{lees}, P. Lees has defined a Steinberg representation for $GL_n$ over a finite ring in a quite similar way, but his representation is different from ours; he uses a smaller set of generalized parahoric subgroups, and thus obtains a larger Steinberg representation. Our representation $st_{\mth{B}}$ is clearly isomorphic to a subrepresentation of his Steinberg representation.

But the representations $st_{\mth{P}}$ are still not small enough to be irreducible in general; they can even have a quite large number of irreducible factors. In theorem \ref{generic}, we prove that for a generic ${\mth{P}}$, the number of factors of $st_{\mth{P}}$ is at least $(q-1)^{\frac{\#(\Phi)}2-rank(\Phi)}$, where $q$ is the cardinal of $k$ and $\Phi$ is the relative root system of $G$.

However, the number of factors is not always that large. In this article, we focus more on the case when $K$ is a maximal special parahoric subgroup of $G$; then all pseudo-Borel subgroups of ${\mth{G}}$ are conjugate to each other.  Our main result is the following one: let $\Delta$ be the set of negative simple roots of $\Phi$ w.r.t ${\mth{B}}$, and let ${\mth{U}}_\Delta$ be a certain abellian unipotent subgroup of ${\mth{G}}$ associated to $\Delta$ which will be defined in section 4. Then we have:

\begin{theo}
Assume the following conditions on $\Phi$ and $p$ are true:
\begin{itemize}
\item $p\neq 2$;
\item if $\Phi$ has at least one irreducible component of type $A_n$, then $p$ doesn't divide the adjoint index of $\underline{\mth{G}}$;
\item if $\Phi$ has at least one irreducible component of type $E_6$, $E_7$ or $F_4$, then $p\neq 3$;
\item if $\Phi$ has at least one irreducible component of type $E_8$, then $p>5$.
\end{itemize}
Then $st_{\mth{B}}$ is multiplicity-free, and its irreducible components are in $1-1$ correspondance with the orbits of the action of ${\mth{T}}$ on the regular characters of ${\mth{U}}_\Delta/$. In particular, if $\underline{\mth{G}}$ is adjoint modulo its center, then $st_{\mth{B}}$ is irreducible.
\end{theo}

Section 2 is devoted to the study of algebraic groups of parahoric type and their pseudo-Borel, pseudo-parabolic and generalized pseudo-parabolic subgroups in a general setting. In section 3, we prove some results about the double classes of groups of parahoric type modulo their generalized pseudo-parabolic subgroups; we'll need those results in section 4 to deal with the Steinberg representations and in particular to prove the main result. The last section is devoted to the (quite long) proof of the proposition \ref{clbi}, which constitutes the most part of the proof of our main theorem.

\section{Some facts about groups of parahoric type}

\subsection{Generalities}

Let $\underline{\mth{G}}$ be a connected algebraic group defined over any field $k$, and let $h$ be a positive integer. The group $\underline{\mth{G}}$ is said to be {\em of parahoric type} of depth $h$ if it satisfies the following conditions:

\begin{itemize}
\item let $\underline{\mth{T}}$ be a maximal torus of $\underline{\mth{G}}$; the set of roots of $\underline{\mth{G}}$ with respect to $\underline{\mth{T}}$ is a root system $\underline{\Phi}$;
\item let $R_u(\underline{\mth{G}})$ be the unipotent radical of $\underline{\mth{G}}$, and let $\underline{\Psi}$ be the root system of the reductive group $\underline{\mth{G}}/R_u(\underline{\mth{G}})$ with respect to $\underline{\mth{T}}$ (viewed as a maximal torus of the quotient); $\underline{\Psi}$ is a root subsystem of $\underline{\Phi}$, and for every $\alpha\in\underline{\Phi}$, the root subgroup $\underline{\mth{U}_\alpha}$ of $\underline{\mth{G}}$ with respect to $\alpha$ is of dimension $h$ (resp. $h-1$) if $\alpha\in\underline{\Psi}$ (resp. if $\alpha\not\in\underline{\Psi}$). We will call $\underline{\mth{G}}$ {\em special} if $\underline{\Psi}=\underline{\Phi}$;
\item let $\underline{\mth{H}}$ be the centralizer of $\underline{\mth{T}}$; $\underline{\mth{H}}$ is abelian, and for every $\alpha\in\underline{\Phi}$, the intersection of $\underline{\mth{H}}$ with the subgroup of $\underline{\mth{G}}$ generated by $\underline{\mth{U}_\alpha}$ and $\underline{\mth{U}_{-\alpha}}$ is of dimension $h$ (resp. $h-1$) if $\alpha\in\underline{\Psi}$ (resp. $\alpha\not\in\underline{\Psi}$);
\item there exists a concave function $f_0$ from $\underline{\Phi}$ to ${\mth{Z}}$ and, for every $\alpha\in\underline{\Phi}$ and every integer $i\geq f_0(\alpha)$, a subgroup $\underline{\mth{U}_{\alpha,i}}$ of $\underline{\mth{U}}_\alpha$ satisfying the following conditions:
\begin{itemize}
\item $\underline{\mth{U}_{\alpha,f_0(\alpha)}}=\underline{\mth{U}_\alpha}$;
\item for every $i$, $\underline{\mth{U}}_{\alpha,i+1}\subset\underline{\mth{U}}_{\alpha,i}$, and if $\underline{\mth{U}_{\alpha,i}}$ is nontrivial, $dim(\underline{\mth{U}_{\alpha,i+1}})=dim(\underline{\mth{U}_{\alpha,i}})-1$;
\item the commutator relations: for every $\alpha,\beta,i,j$ such that $\alpha+\beta\in\underline{\Phi}$, we have $[Lie(\underline{\mth{U}_{\alpha,i}}),Lie(\underline{\mth{U}_{\beta,j}})]=Lie(\underline{\mth{U}_{\alpha+\beta,i+j}})$,
\item for every $\alpha,i$ such that $i\geq f_0(\alpha)+f_0(-\alpha)$, the dimension of the subalgebra $L_{\alpha_i}=[Lie(\underline{\mth{U}}_\alpha),Lie(\underline{\mth{U}}_{-\alpha,i-f_0(\alpha)})]$ of $Lie(\underline{\mth{H}})$ is $Sup(h-i,0)$, and for every $j\geq f_0(\alpha)$, we have $[L_{\alpha_i}, Lie(\underline{\mth{U}}_{\alpha,j})]=Lie(\underline{\mth{U}}_{\alpha,i+j})$.
\end{itemize}
\end{itemize}

Remember that a concave function $f_0$ on $\underline{\Phi}$ is a function satisfying the following properties:
\begin{itemize}
\item for every $\alpha\in\underline{\Phi}$, $f_0(\alpha)+f_0(-\alpha)\geq 0$;
\item for every $\alpha,\beta\in\underline{\Phi}$ such that $\alpha+\beta\in\underline{\Phi}$, $f_0(\alpha+\beta)\leq f_0(\alpha)+f_0(\beta)$.
\end{itemize}
In this article all concave functions will be ${\mth{Z}}$-valued.

Since all maximal tori of $\underline{\mth{G}}$ are conjugated, these properties don't depend on the choice of $\underline{\mth{T}}$.

The last two properties are equivalent to the following ones in terms of subgroups:
\begin{itemize}
\item the commutator relations: for every $\alpha,\beta,i,j$ such that $\alpha+\beta\in\underline{\Phi}$, the group $[\underline{\mth{U}_{\alpha,i}},\underline{\mth{U}_{\beta,j}}]$ is contained in the product of the $\underline{\mth{U}}_{\lambda\alpha+\mu\beta,\lambda i+\mu j}$, $\lambda,\mu>0$, and its canonical projection on $\underline{\mth{U}_{\alpha+\beta,i+j}}$ is injective,
\item for every $\alpha,i$ such that $i\geq f_0(\alpha)+f_0(-\alpha)$, the canonical projection $\underline{\mth{H}}_{\alpha,j}$ of $[\underline{\mth{U}}_\alpha),\underline{\mth{U}}_{-\alpha,i-f_0(\alpha)}]$ on $\underline{\mth{H}}$ is of dimension $Sup(h-i,0)$, and for every $j\geq f_0(\alpha)$, we have $[\underline{\mth{H}}_{\alpha_i}, \underline{\mth{U}}_{\alpha,j}]=\underline{\mth{U}}_{\alpha,i+j}$.
\end{itemize}

We deduce from the second and third conditions that $R_u(\underline{\mth{G}})$ is generated by $R_u(\underline{\mth{H}})$ and subgroups of dimension $h-1$ of the $\underline{\mth{U}}_\alpha$, $\alpha\in\underline{\Phi}$. In particular, when $h=1$, $\underline{\mth{G}}$ is simply a reductive group.

For every $\alpha\in\underline{\Phi}$ and every $u\in\underline{\mth{U}}_\alpha$, we'll define the valuation $v(u)$ of $u$ as the largest integer $v$ such that $u\in \underline{\mth{U}}_{\alpha,v}$. By convention the valuation of the identity element is infinite.

An immediate consequence of the last condition is that for every $\alpha\in\underline{\Phi}$, we have $f_0(\alpha)+f_0(-\alpha)=0$ (resp. $f_0(\alpha)+f_0(-\alpha)=1$) if $\alpha\in\underline{\Psi}$ (resp. $\alpha\not\in\underline{\Psi}$). Moreover, $f_0$ is entirely determined by its values on the elements of a given set of simple roots of $\underline{\Phi}$, and these values can be chosen arbitrarily. We can for example choose $f_0(\alpha)=0$ for every $\alpha$ belonging to our set of simple roots; if $\underline{\mth{G}}$ is special, we then have $f_0=0$.

Another consequence is that for every $\alpha\in\underline{\Phi}$ and every integer $i$, $\underline{\mth{U}}_{\alpha,i}=\{1\}$ if and only if $i\geq h-f_0(-\alpha)$.

\begin{prop}\label{gfparah}
Assume $k$ is perfect. Let $\underline{G}$ be a connected reductive algebraic group defined over a local field $F$ admitting $k$ as its residual field, and split over the unramified closure $F_{nr}$ of $F$, and let $G_{nr}$ be the group of $F_{nr}$-points of $\underline{G}$. Assume the residual characteristic $p$ of $F$ and the root system $\underline{\Phi}$ of $\underline{G}$ satisfy one of the following conditions:
\begin{itemize}
\item $p>3$;
\item $p=3$ and $\underline{\Phi}$ has no irreducible component of type $G_2$;
\item $p=2$ and all irreducible components of $\underline{\Phi}$ are of type $A_n$ for some $n$.
\end{itemize}
Let $K$ be a parahoric group of $G_{nr}$; assume $K$ is stable by the action of $Aut(F_{nr}/F)$ over $G_{nr}$. For any integer $h>0$, let $K^h$ be the $h$-th congruence subgroup of $K$. Then $\underline{\mth{G}}=K/K^h$ is an algebraic group of parahoric type of depth $h$ defined over the residual field $k$ of $F$.
\end{prop}

The fact that $\underline{\mth{G}}$ satisfies the first three properties of groups of parahoric type is an easy consequence of the definitions. The existence of the concave function $f$ and the subgroups $\underline{\mth{U}}_{\alpha,i}$ follows from the existence of a valued root datum on $G_{nr}$ (see \cite[I; 6.1]{bt}). When $p$ satisfies the required conditions, the commutator relations come from \cite[Theorem 1]{chev}, and the condition on the $L_(\alpha_j]$ is easy to check directly.

\subsection{Pseudo-Borel and pseudo-parabolic subgroups}

Let $\underline{\Phi}^+$ be a set of positive roots in $\underline{\Phi}$, and let $\underline{\mth{B}}$ be the subgroup of $\underline{\mth{G}}$ generated by $\underline{\mth{H}}$ and the $\underline{\mth{U}_\alpha}$, $\alpha\in\underline{\Phi}^+$; $\underline{\mth{B}}$ is called a pseudo-Borel subgroup of $\underline{\mth{G}}$. We'll also write $\underline{\Phi}^-$ for the corresponding set of negative toots, which are the opposites of the elements of $\underline{\Phi}^+$.

A pseudo-Borel subgroup of $\underline{\mth{G}}$ is solvable. It is easy to check this fact directly with the commutator relations; we can also make the simple remark that a pseudo-Borel subgroup of $\underline{\mth{G}}$ is obviously contained in a true Borel subgroup.

Note that in the general case, two pseudo-Borel subgroups of $\underline{\mth{G}}$ are not necessarily conjugated. In fact we have the following result:

\begin{prop}
Let $\underline{\mth{B}}$ and $\underline{\mth{B}'}$ be two pseudo-Borel subgroups of $\underline{\mth{G}}$ containing $\underline{\mth{T}}$, and let $\underline{\Phi}^+$ and $\underline{\Phi}'^+$ be the corresponding sets of positive roots in $\underline{\Phi}$. The subgroups $\underline{\mth{B}}$ and $\underline{\mth{B}}'$ are conjugated in $\underline{\mth{G}}$ if and only if there exists an element $w$ of the Weyl group of $\underline{\Psi}$ such that $w\underline{\Phi}^+=\underline{\Phi'}'^+$.
\end{prop}

Assume $\underline{\mth{B}}$ and $\underline{\mth{B}}'$ are conjugated. Since $\underline{\mth{T}}$ is a maximal torus of both, there exists an element $n$ of its normalizer $N_{\underline{\mth{G}}}(\underline{\mth{T}})$ in $\underline{\mth{G}}$ such that $n\underline{\mth{B}}n'=\underline{\mth{B}}'$, and $n$ must belong to some element $w$ of the Weyl group $W$ of $\underline{\mth{G}}$ relatively to $\underline{\mth{T}}$, which is also the Weyl group of $\underline{\Psi}$; $w$ then satisfies $w\underline{\Phi}^+=\underline{\Phi'}'^+$. Since the converse is obvious, the proposition follows. $\Box$

\begin{cor}
Assume $\underline{\mth{G}}$ is special. Then all its pseudo-Borel subgroups are conjugated.
\end{cor}

Let $\underline{\mth{B}}$ and $\underline{\mth{B}}'$ be two pseudo-Borel subgroups of $\underline{\mth{G}}$. By replacing them by conjugates we may assume both of them contain $\underline{\mth{T}}$; let then $\underline{\Phi}^+$ and $\underline{\Phi}'^+$ be the corresponding sets of positive roots in $\underline{\Phi}$. By \cite[1, cor. 1 to prop. 20 and th. 2]{bou}, they are conjugated by an element of the Weyl group of $\underline{\Phi}$; on the other hand, since $\underline{\mth{G}}$ is special, we have $\underline{\Phi}=\underline{\Psi}$. The above proposition then immedately implies that $\underline{\mth{B}}$ and $\underline{\mth{B}}'$ are conjugated. $\Box$

Now let $\underline{\Phi}'$ be a parabolic subset of $\underline{\Phi}$; that is a closed subset of $\underline{\Phi}$ containing a set of positive roots. The group $\underline{\mth{P}}$ generated by $\underline{\mth{H}}$ and tue $\underline{\mth{U}_\alpha}$, $\alpha\in\underline{\Phi}'$, is a pseudo-parabolic subgroup of $\underline{\mth{G}}$.

Note that our definition of pseudo-parabolic subgroups is compatible with the definition in \cite[15.1]{spr}: let $\underline{\Delta}$ be a set of simple roots of $\underline{\Phi}$ contained in $\underline{\Phi}'$, and let $\underline{\Delta}'$ be the subset of $\underline{\Delta}$ containing the $\alpha$ such that $-\alpha\in\underline{\Phi}'$. Then $\underline{\Delta}'$ doesn't depend on the choice of $\underline{\Delta}$, and with the notations of \cite{spr}, $\underline{\mth{P}}$ is the group $\underline{\mth{P}}(\lambda)$, where $\lambda$ is a cocharacter of $\underline{\mth{G}}$ such that for every $\alpha$ in $\underline{\Delta}'$ (resp. $\underline{\Delta}-\underline{\Delta}'$), with the usual pairing, $<\alpha,\lambda>=0$ (resp.$<\alpha,\lambda>\ >0$).

We have a pseudo-Levi decomposition:
\[\underline{\mth{P}}=\underline{\mth{M}}\underline{\mth{U}}\]
where $\underline{\mth{M}}$ is the subgroup of $\underline{\mth{G}}$ generated by $\underline{\mth{H}}$ and the $\alpha\in\underline{\Phi'}$ whose opposite is also in $\underline{\Phi'}$, and $\underline{\mth{U}}$ is the subgroup of $\underline{\mth{G}}$ generated by the $\alpha\in\underline{\Phi'}$ whose opposite is not in $\underline{\Phi'}$. The group $\underline{\mth{M}}$ is a group of parahoric type, which will be called a pseudo-Levi subgroup of $\underline{\mth{G}}$, or a pseudo-Levi component of $\underline{\mth{P}}$.

We can easily see with the help of the commutator relations that the group $\underline{\mth{U}}$ is normal in $\underline{\mth{P}}$. Since all maximal tori of $\underline{\mth{P}}$ are conjugated in $\underline{\mth{P}}$, all pseudo-Levi decompositions of $\underline{\mth{P}}$ yield the same group $\underline{\mth{U}}$, which thus behaves like a sort of unipotent radical; it will be called the pseudo-unipotent radical of $\underline{\mth{G}}$.

\begin{prop}\label{gfpsb}
Let $G_{nr}$, $K$, $h$ and $\underline{\mth{G}}$ be as in proposition \ref{gfparah}. Let $B_{nr}$ (resp. $P_{nr}$) be a Borel (resp. parabolic) subgroup of $G_{nr}$, and let $\underline{\mth{B}}$ (resp. $\underline{\mth{P}}$) be the image in $\underline{\mth{G}}$ of its intersection with $K$. Then $\underline{\mth{B}}$ (resp. $\underline{\mth{P}}$) is a pseudo-Borel (resp. pseudo-parabolic) subgroup of $\underline{\mth{G}}$.
\end{prop}

Let $T_{nr}$ be a maximal torus of $G_{nr}$ whose parahoric subgroup $K_T$ is contained in $K$, and let $I$ be an Iwahori subgroup of $G_{nr}$ contained in $K$ and also containing $K_T$. Let $W$ be the Weyl group of $G_{nr}$ relatively to $T_{nr}$; we have the following decomposition:

\begin{lemme}
For a suitable set $R$ of representants of the elements of $W$ in the normalisator of $T_{nr}$ in $G_{nr}$, we have $G=\bigcup_{n\in R}nIB_{nr}$ (resp. $G=\bigcup_{n\in R}nIP_{nr}$).
\end{lemme}

It is obviously enough to prove it for $B_{nr}$. Let $U^(_{nr}$ be the unipotent radical of the unique Borel subgroup $B_{nr}^-$ of $G$ opposite to $B_{nr}$ and containing $T_{nr}$; for every $n\in R$, we have the Iwahori decomposition $n^{-1}In=(n^{-1}In\cap B_{nr})(n^{-1}In\cap U_{nr}^-$). Consider now the Iwasawa decomposition $G_{nr}=\bigcup_{n\in R}InB_{nr}$; we can rewrite it the following way:
\[G_{nr}=\bigcup_{n\in R}n(n^{-1}In)B_{nr}=\bigcup_{n\in R}n(n^{-1}In\cap U^-_{nr})B_{nr}.\]
Since we always can choose the set $R$ in such a way that $I\cap U^-_{nr}\subset n^{-1}In\cap U^-_{nr}$ for every $n\in R$, the result follows. $\Box$

We deduce from this lemma that if $T'_{nr}$ is any maximal torus of $B_{nr}$, there exist $l\in I$ and $b\in B_{nr}$ such that $T'_{nr}=b^{-1}l^{-1}T_{nr}lb$, which means that by eventually replacing $T_{nr}$ by $l^{-1}T_{nr}l$, we can assume it is contained in $B_{nr}$ (resp. $P_{nr}$).  The group $\underline{\mth{B}}$ (resp. $\underline{\mth{P}}$) is then generated by the image $\underline{\mth{T}}$ of $K_T$ in $\underline{\mth{G}}$ and the root subgroups corresponding to the roots of $B_{nr}$ (resp. $P_{nr}$) relatively to $T_{nr}$; the fact that it is a pseudo-Borel (resp. pseudo-parabolic) subgroup of $\underline{\mth{G}}$ follows now immediately from the definitions. $\Box$

\subsection{Generalized pseudo-parabolic subgroups}

It is obvious from the definitions that all pseudo-parabolic subgroups of $\underline{\mth{G}}$ are groups containing a pseudo-Borel subgroup. However, the converse is not true: for example, any pseudo-Borel subgroup of $\underline{\mth{G}}$ is solvable, hence contained in a true Borel subgroup, but the Borels of $\underline{\mth{G}}$ are not pseudo-parabolics unless $\underline{\mth{G}}$ himself is solvable. We will thus define a larger class of subgroups of $\underline{\mth{G}}$, which we will call generalized pseudo-parabolic subgroups.

Let $\underline{\mth{T}}$, $\underline{\Phi}$ and $f_0$ be defined as in the previous subsections. Let $f$ be a concave function on $\underline{\Phi}$ such that $f(\alpha)\geq f_0(\alpha)$ for every $\alpha$. Then $f$ defines the subgroup $\underline{\mth{P}}_f$ (resp. $\underline{\mth{U}}_f$) of $\underline{\mth{G}}$, generated by $\underline{\mth{H}}$ (resp. the unipotent radical of $\underline{\mth{H}}$) and the $\underline{\mth{U}_{\alpha,f(\alpha)}}$, $\alpha\in\underline{\Phi}$. We deduce from the concavity of $f$ and the commutator relations that for every $\alpha\in\underline{\Phi}$, $\underline{\mth{U}}_{\alpha}\cap\underline{\mth{P}}_f=\underline{\mth{U}}_{\alpha,f(\alpha)}$

The groups $\underline{\mth{P}}_f$ and $\underline{\mth{U}}_f$ are connected algebraic groups.
The group $\underline{\mth{P}}_f$ is in general not reductive; its unipotent radical is $\underline{\mth{U}}_{f'}$, where $f'$ is the concave function defined by $f'(\alpha)=f(\alpha)+1$ if $\alpha\in\underline{\Psi}$ and $f(\alpha)=f_0(\alpha)$ and $f'(\alpha)=f(\alpha)$ if $\alpha\not\in\underline{\Psi}$ or $f(\alpha)>f_0(\alpha)$. It is solvable if and only if the following equivalent conditions are satisfied:
\begin{itemize}
\item for every $\alpha\in\underline{\Phi}$, $f(\alpha)+f(-\alpha)>0$;
\item for every $\alpha\in\underline{\Psi}$, either $f(\alpha)>f_0(\alpha)$ or $f(-\alpha)>f_0(-\alpha)$.
\end{itemize}
In this case, we have a decomposition $\underline{\mth{P}}_f=\underline{\mth{H}}\prod_{\alpha\in\Phi}\underline{\mth{U}}_{\alpha,f(\alpha)}$. In the general case, this product is an open subset of $\underline{\mth{P}}_f$.

We'll say $\underline{\mth{P}}_f$ is a generalized pseudo-parabolic subgroup of $\underline{\mth{G}}$ if it contains at least one pseudo-Borel subgroup of $\underline{\mth{G}}$ containing $\underline{\mth{T}}$, or in other words if there exists a set of positive roots $\underline{\Phi}^+$ of $\underline{\Phi}$ such that $f(\alpha)=f_0(\alpha)$ for every $\alpha\in\underline{\Phi}^+$.

In particular, the Borel and parabolic subgroups of $\underline{\mth{G}}$ are generalized pseudo-parabolic subgroups of $\underline{\mth{G}}$; a Borel (resp. parabolic) subgroup containing $\underline{\mth{T}}$  is a subgroup of the form $\underline{\mth{P}}_f$ such that $f(\alpha)+f(-\alpha)=1$ (resp. $f(\alpha)+f(-\alpha)\leq 1$) for every $\alpha\in\underline{\Phi}$. We then have the following result:

\begin{prop}
Every solvable generalized pseudo-parabolic subgroup of $\underline{\mth{G}}$ is contained in a unique Borel subgroup of $\underline{\mth{G}}$.
\end{prop}

Let $\underline{\mth{P}}$ be a solvable generalized pseudo-parabolic subgroup of $\underline{\mth{G}}$. By eventually replacing it by a conjugate we can assume it contains $\underline{\mth{T}}$; let then $f$ be the concave function on $\underline{\Phi}$ associated to $\underline{\mth{P}}$. Let $S$ be the subset of the $\alpha\in\underline{\Psi}$ such that $f(\alpha)=f_0(\alpha)$; by the definition of generalized pseudo-parabolic subgroups, $S$ contains a set of positive roots in $\underline{\Psi}$; on the other hand, since $\underline{\mth{P}}$ is solvable, for every $\alpha\in S$, $f(\alpha)+f(-\alpha)>0$, hence $f'(-\alpha)>f_0(-\alpha)$, which shows that $S$ is precisely that set of positive roots. Let $\underline{\mth{B}}$ be the Borel subgroup of $\underline{\mth{G}}$ generated by $R_u(\underline{\mth{G}})$ and the $\underline{\mth{U}}_\alpha$, $\alpha\in S$; obviously $\underline{\mth{P}}$ is contained in $\underline{\mth{B}}$. If now $\underline{\mth{B}}'$ is another Borel subgroup of $\underline{\mth{G}}$ containing $\underline{\mth{P}}_f$, then $\underline{\mth{B}}'$ must contain $\underline{\mth{T}}$ and the $\underline{\mth{U}}_\alpha$, $\alpha\in S$, which is only possible if $\underline{\mth{B}}'=\underline{\mth{B}}$. $\Box$

Let $f,f'$ be two concave functions on $\Phi$: we obviously have $\underline{\mth{P}}_f\supset\underline{\mth{P}}_{f'}$ if and only if $f(\alpha)\leq f'(\alpha)$ for every $\alpha\in\Phi$. We'll write $f\leq f'$ if these propositions are true; this defines a partial order on the set of concave functions on $\Phi$.

Let $\underline{\Phi}^+$ be a set of positive roots in $\Phi$ and let $\underline{\mth{B}}$ be the corresponding pseudo-Borel subgroup of $\underline{\mth{G}}$. It is obvious that if $f(\alpha)=f_0(\alpha)$ for every $\alpha\in\underline{\Phi}^+$, $\underline{\mth{P}}_f$ contains $\underline{\mth{B}}$; conversely, we have the following result:

\begin{prop}
Assume $p\neq 2$ or $\underline{\Phi}$ has no irreducible components of type $A_1$ or $C_n$.
\begin{itemize}
\item Every subgroup of $\underline{\mth{G}}$ containing $\underline{\mth{B}}$ is of the form $\underline{\mth{P}}_f$ for some $f$.
\item For every $f$ such that $\underline{\mth{B}}\subset \underline{\mth{P}}_f$, $\underline{\mth{P}}_f$ is its own normalizer in $\underline{\mth{G}}$.
\end{itemize}
\end{prop}

Let $\underline{\mth{P}}$ be a subgroup of $\underline{\mth{G}}$ containing $\underline{\mth{B}}$. According to the previous remark, in order to prove that $\underline{\mth{P}}$ is of the form $\underline{\mth{P}}_f$ for some $f$, it is enough to prove that if $u=\prod_{\alpha\in\underline{\Phi}^-}u_\alpha$ belongs to $\underline{\mth{P}}$, then $\underline{\mth{P}}$ contains all the subgroups $\underline{\mth{U}}_{\alpha,v(u_\alpha)}$, $\alpha\in\underline{\Phi}^-$.

For every concave function $f$ on $\underline{\Phi}^-$, let $\underline{\mth{U}}^-_f$ be the product of the $\underline{\mth{U}}_{\alpha,f(\alpha)}$, $\alpha\in\underline{\Phi}^-$. We'll prove the claim by descending induction on the largest concave function $f$ on $\underline{\Phi}^-$ such that $u$ is contained in $\underline{\mth{U}}^-_f$, given that when $f$ is the restriction to $\underline{\Phi}^-$ of the concave function $f_{\underline{\mth{B}}}$ associated to $\underline{\mth{B}}$, then $u=0$ and the assertion is trivial. First assume there exists an unique element $\alpha$ of $\underline{\Phi}^-$ such that $f(\alpha)<f_{\underline{\mth{B}}}(\alpha)$. Consider an isomorphism $\phi$ between the quotient $\underline{\mth{U}}_{\alpha,v(u_\alpha)}/\underline{\mth{U}}_{\alpha,v(u_\alpha)+1}$ and $\overline{k}$; lifting this morphism to $\underline{\mth{U}}_{\alpha,v(u_\alpha)}$, we claim that its restriction to $\underline{\mth{P}}\cap\underline{\mth{U}}_\alpha$ is surjective. Assume first there exists a one-parameter group $\xi$ of $\underline{\mth{T}}$ such that $<\alpha,\xi>=1$ for the usual pairing; this is always the case if $\underline{\Phi}$ has no irreducible components of type $A_1$ or $C_n$. Then for every $x\in\overline{k}^*$, we have $\phi(\xi(x)u_\alpha\xi(x^{-1})=x\phi(u_\alpha)$, hence the image of $\underline{\mth{P}}\cap\underline{\mth{U}}_\alpha$ by $\phi$ contains all nonzero elements of $\overline{k}$, hence contains $0$ as well. Assume now such a $\xi$ doesn't exist; we will consider now the coroot $\alpha^\vee$ associated to $\alpha$ in the root datum $(X^*(\underline{\mth{T}},\underline{\mth{\Phi}},X_*(\underline{\mth{T}},\underline{\mth{\Phi}}^\vee)$ of $\underline{\mth{G}}$; it always satisfies $<\alpha,\alpha^\vee>=2$. Writing $x=\phi(u_\alpha)$, for every $y\in\overline{k}^*$, we have $\phi(\alpha^\vee(y)u_\alpha\alpha^\vee(y^{-1})=y^2x$. Hence the image of $\underline{\mth{P}}\cap\underline{\mth{U}}_\alpha$ contains the subgroup of $\overline{k}$ generated by the $y^2x$.
 On the other hand, by the hypotheses we have made, $p\neq 2$, and for every $z\in\mathcal{O}$, we have:
\[z=(\frac{z+1}2)^2-(\frac{z-1}2)^2,\]
which proves the claim.

Moreover, since $\underline{\mth{P}}$ contains $\underline{\mth{H}}$, it contains in particular all the $\underline{\mth{H}}_{\alpha,i}$, fron which we deduce that $\underline{\mth{P}}$ contains elements of $\underline{\mth{U}}_\alpha$ of any valuation greater than or equal to the valuation of $u_\alpha$. Applying the above claim to all these valuations, we finally obtain that $\underline{\mth{P}}$ contains the whole group $\underline{\mth{U}}_{\alpha,v(u_\alpha)}$.

Now let's consider the general case. We have a canonical partial order on $\underline{\Phi}^-$ defined by $\alpha\leq\beta$ if there exist $\alpha_0=\beta,\alpha_1,\dots,\alpha_n=\alpha$ such that $\alpha_i-\alpha_{i-1}\in\underline{\Phi}^-$ for every $i$. Let $\alpha$ be an element of $\underline{\Phi}^-$ satisfying the following conditions:
\begin{itemize}
\item $f(\alpha)-f_0(\alpha)$ is minimal;
\item $\alpha$ is minimal for the above condition.
\end{itemize}
For such an $\alpha$ we obviously have $f(\alpha)=v(u_\alpha)$. We'll prove that $\underline{\mth{U}}^-_f$ contains an element $u'_\alpha$ of $\underline{\mth{U}}_\alpha$ of the same valuation as $u_\alpha$. If this is true, by the same argument as above it contains the whole group $\underline{\mth{U}}_{\alpha,v(u_\alpha)}$, and in particular $u_\alpha$, and we can then apply the induction hypothesis to $u_\alpha^{-1}u$, which is contained in $\underline{\mth{P}}_f'$ with $f'$ such that $f'(\alpha)>f(\alpha)$, to get the desired result.

Let $\alpha'$ be an element of $\underline{\Phi}^-\{\alpha\}$ such that $f(\alpha')-f_0(\alpha)$ is minimal and $\alpha'$ is minimal for that property, and let $f''$ be the function on $\underline{\Phi}^-$ defined by $f''(\alpha)=f(\alpha)+1$, $f''(\alpha')=f(\alpha')+1$ and for every $\beta\neq\alpha, \alpha'$, $f''(\beta)=f(\beta)$. It is easy to check that $f'$ is concave;
moreover, for every $\beta,\gamma\in\underline{\Phi}^-$ such that $\beta+\gamma=\alpha$, $f(\beta)-f_0(\beta)$ and $f(\gamma)-f_0(\gamma)$ are strictly greater than $f(\alpha)-f_0(\alpha)$, hence:
\[f'(\alpha)=f(\alpha)+1=(f(\alpha)-f_0(\alpha))+f_0(\alpha)+1\]
\[\leq (f(\beta)-f_0(\beta))+(f(\gamma)-f_0(\gamma))-1+f_0(\alpha)+1\]
\[\leq f'(\beta)+f(\gamma)+(f_0(\alpha)-f_0(\beta)-f_0(\gamma))\leq f'(\beta)+f(\gamma).\]
We have a similar assertion for $\alpha'$ (which holds even when either $\beta$ or $\gamma$ is equal to $\alpha$); hence $\underline{\mth{U}}_{f''}$ is normal in $\underline{\mth{U}}_f$.

Moreover, the quotient $\underline{\mth{U}}_f/\underline{\mth{U}}_{f''}$ is abelian and isomorphic to $\overline{k}^2$, hence can be viewed as a $2$-dimensional $\overline{k}$-vector space on which $\underline{\mth{T}}$ acts, and its weight subspaces are the images of respectively $\underline{\mth{U}}_{\alpha,f(\alpha)}$ and $\underline{\mth{U}}_{\alpha',f(\alpha')}$, the corresponding weights being of course $\alpha$ and $\alpha'$. Consider now the quotient $(\underline{\mth{P}}\cap\underline{\mth{U}}_f)/(\underline{\mth{P}}\cap\underline{\mth{U}}_{f''})$: it is a $\underline{\mth{T}}$-stable subspace of $\underline{\mth{U}}_f/\underline{\mth{U}}_{f''}$, which contains at least one element whose projection on the weight subspace associated to $\alpha$ is nonzero; it then contains that weight subspace, which amounts to say that $\underline{\mth{P}}$ contains an element of $\underline{\mth{U}}_{\alpha,v(u_\alpha)}\underline{\mth{P}}_{f''}=\underline{\mth{P}}_{f'}$, with $f'$ being a concave function strictly larger than $f$ and such that $f'(\alpha)=f(\alpha)$. By iterating the process, after a finite number of steps we reach the point where $f'(\beta)=f_{\underline{\mth{B}}}(\beta)$ for every $\beta\neq\alpha$, which proves the desired assertion.

Let now $\underline{\mth{L}}$ be the normalizer of $\underline{\mth{P}}$ in $\underline{\mth{G}}$. $\underline{\mth{L}}$ also contains $\underline{\mth{B}}$, hence is of the form $\underline{\mth{P}}_{f'}$ for some $f'\leq f$. Assume there exists an $\alpha$ such that $f'(\alpha)<f(\alpha)$. Let $t$ be an element of $\underline{\mth{T}}$ such that $\alpha(t)\neq 1$; for every $u\in \underline{\mth{U}_{\alpha,f'(\alpha)}}-\underline{\mth{U}_{\alpha,f(\alpha)}}$  we must then have $tut^{-1}u^{-1}\in \underline{\mth{U}_{\alpha,f'(\alpha)}}-\underline{\mth{U}_{\alpha,f(\alpha)}}$, which is impossible. Hence $f'=f$ and the proposition is proved. $\Box$

Now we would like to generalize the notion of pseudo-Levi component to generalized pseubo-parabolic subgroups. We can define the pseudo-unipotent radical of a generalized pseudo-parabolic subgroup of $\underline{\mth{G}}$ the following way: let $\underline{\mth{P}}$ be a generalized pseudo-parabolic subgroup and let $f$ be the concave function defining it; set for every $\alpha\in\underline{\Phi}^+$:
\[f'(\alpha)=h-f(-\alpha)\]
and let $\underline{\mth{U}}$ be the subgroup of $\underline{\mth{G}}$ generated by the $U_{\alpha,f'(\alpha)}$, $\alpha\in\underline{\Phi}^+$. For every $\alpha,\beta\in\underline{\Phi}$ such that $\alpha+\beta\in\underline{\Phi}$, we have:
\[f'(\alpha+\beta)=h-f(-\alpha-\beta)\leq h-f(-\alpha)+f_0(\beta)=f'(\alpha)+f_0(\beta),\]
hence $\underline{\mth{U}}$ is normal in $\underline{\mth{P}}$. The quotient $\underline{\mth{P}}/\underline{\mth{U}}$ will be called the pseudo-Levi quotient of $\underline{\mth{P}}$. Note that in general there is no subgroup of $\underline{\mth{P}}$ playing the role of a pseudo-Levi component of $\underline{\mth{P}}$.

Now we'll determine the smallest generalized pseudo-parabolic subgroups of $\underline{\mth{G}}$ strictly containing a given generalized pseudo-parabolic subgroup $\underline{\mth{P}}$. Let $f$ be the concave function associated to $\underline{\mth{P}}$, and for every $\alpha\in\underline{\Phi}$ such that $f(\alpha)>f_0(\alpha)$ (which implies in particular $\alpha\in\underline{\Phi}^-$), let $f_\alpha$ be the concave function satisfying the following conditions:
\begin{itemize}
\item $f_\alpha\leq f$;
\item $f_\alpha(\alpha)<f(\alpha)$;
\item $f_\alpha$ is the largest concave function on $\Phi$ satisfying the above conditions.
\end{itemize}
The function $f_\alpha$ can be directly defined the following way: let $f'_\alpha$ be the function on $\underline{\Phi}$ such that $f'_\alpha(\alpha)=f(\alpha)-1$ and $f'_\alpha(\beta)=f(\beta)$ for every $\beta\neq\alpha$, and set:
\[f_\alpha(\beta)=Inf_{\beta_1+\dots+\beta_s=\beta}\sum_{i=1}^sf'_\alpha(\beta_i)\]
for every $\beta\in\underline{\Phi}$. In particular, we have $f_\alpha(\alpha)=f(\alpha)-1$. We'll denote by $\underline{\mth{P}}_\alpha$ the generalized parabolic subgroup of $\underline{\mth{G}}$ associated to $f_\alpha$.

Let's now consider another root $\beta$ such that $f_\alpha(\beta)<f(\beta)$; we then have $f_\alpha\leq f'_\beta$, hence by maximality of $f_\beta$, $f_\alpha\leq f_\beta$. Assume now $f_\alpha=f_\beta$; this implies:
\begin{itemize}
\item there exist $\alpha_1,\dots,\alpha_s\in\underline{\Phi}$ such that $\alpha_1=\alpha$, $\alpha_1+\dots+\alpha_s=\beta$ and $f(\beta)\leq\sum_{i=1}^sf(\alpha_i)$;
\item there exist $\beta_1,\dots,\beta_r\in\underline{\Phi}$ such that  $\beta_1=\beta$, $\beta_1+\dots+\beta_r=\beta$  and $f(\alpha)\leq\sum_{i=1}^rf(\beta_i)$.
\end{itemize}
Hence $\sum_{i=2}^s\alpha_i+\sum_{i=2}^r\beta_i$ is zero and $\sum_{i=1^sf(\alpha_i)+\sum_[i=2}^rf(\beta_i)\leq 0$; the concavity of $f$ implies this last sum is zero as well, and we must then have $f(\alpha_i)=f_0(\alpha_i)$ for every $i\geq 2$, and similarly for the $\beta_i$. In other words, if $\underline{\mth{P}}_0$ is the largest pseudo-parabolic subgroup of $\underline{\mth{G}}$ contained in $\underline{\mth{P}}$, all the $\alpha_i$ and the $\beta_i$ are both roots of its Levi component $\underline{\mth{M}}_0$ and elements of $\underline{\Psi}$, hence $\alpha-\beta$ is a sum of such roots. Conversely, it is easy to see that if $\alpha-\beta$ is of that form, $f_\alpha=f_\beta$.

Let $\underline{\Delta}_f$ be the set of elements $\alpha\in\underline{\Phi}^-$ such that:
\begin{itemize}
\item $f_\alpha$ is maximal among the $f_\beta$, $\beta\in\underline{\Phi}^-$;
\item $\alpha$ is maximal among the negative roots contained in its equivalence class modulo the subgroup of $X_*(\underline{\mth{T}})$ generated by the roots of $\underline{\mth{M}}_0$ contained in $\underline{\Psi}$.
\end{itemize}
We deduce from the above discussion that the $\underline{\mth{P}}_\alpha$, $\alpha\in\underline{\Delta}_f$, are exactly the minimal generalized parahoric subgroups of $\underline{\mth{G}}$ strictly containing $\underline{\mth{P}}$.

In particular, let's suppose $\underline{\mth{G}}$ special and consider the case of the pseudo-Borel subgroup $\underline{\mth{B}}$ corresponding to the concave function $f$ such that $f(\alpha)=0$ if $\alpha>0$ and $f(\alpha)=h$ fs $\alpha<0$. It is easy to see that for every $\alpha\in\underline{\Phi}^-$ and every $\beta\in\underline{\Phi}$, $f_\alpha(\beta)=h-1$ if $\alpha\leq\beta<0$ and $f_\alpha(\beta)=f(\beta)$ otherwise; hence if $\alpha$ is a simple root in $\underline{\Phi}^-$, $f_\alpha$ is obviously maximal; on the other hand, if $\alpha$ is not simple, by \cite[1, prop. 19]{bou}, there exists at least one simple root $\beta$ in $\underline{\Phi}^-$ such that $\alpha-\beta$ is a root. Hence $\underline{\Delta}_f$ is exactly the set of simple roots in $\underline{\Phi}^-$.

For every subset $I$ of $\underline{\Delta}_f$, we will also define the generalized parahoric subgroup $\underline{\mth{P}}_I$ as the subgroup of $\underline{\mth{G}}$ generated by the $\underline{\mth{P}}_\alpha$, $\alpha\in I$. Note that those subgoups are not necessarily distinct; it may happen that $\underline{\mth{P}}_I=\underline{\mth{P}}_J$ for two distincts subsets $I$ and $J$ of $\underline{\Delta}_f$.

\subsection{Groups of rational points}

Let now ${\mth{G}}$ be the set of $k$-points of $\underline{\mth{G}}$; we will also assume $\underline{\mth{T}}$ is defined over $k$. Then the group $\Gamma=Aut(\overline{k}/k)$  acts on the root systems $\underline{\Phi}$ and $\underline{\Psi}$.

Let $\underline{\mth{S}}$ be a maximal $k$-split torus of $\underline{\mth{G}}$. Assume $\underline{\mth{T}}$ contains $\underline{\mth{S}}$; we have the following result, which is an extension of the well-known similar result for reductive groups:

\begin{prop}
The following conditions are equivalent:
\begin{itemize}
\item the centralizer of $\underline{\mth{S}}$ is a Cartan subgroup of $\underline{\mth{G}}$;
\item there exists a $\Gamma$-stable set of positive roots in $\underline{\Phi}$;
\item there existe a $\Gamma$-stable pseudo-Borel subgroup of $\underline{\mth{G}}$ containing $\underline{\mth{T}}$.
\end{itemize}
\end{prop}

The last two conditions are obviously equivalent. Assume they are true. Since the centralizer of $\underline{\mth{S}}$ contains $\underline{\mth{T}}$, it is generated by $\underline{\mth{H}}$ and the $\underline{\mth{U}}_\alpha$, where $\alpha$ belongs to the set of elements of $\underline{\Phi}$ which are trivial on $\underline{\mth{S}}$. Let $\alpha$ be such a root; for every $\sigma\in\Gamma$, $\sigma(\alpha)$ satisfies the same condition. Let $\Gamma'$ be a subgroup of $\Gamma$ of finite index acting trivially on $\underline{\Phi}$; such a subgroup exists since $\underline{\Phi}$ is finite, and the $\Gamma$-stable element $\sum_{\sigma\in\Gamma/\Gamma'}\sigma(\alpha)$ of $X^*(\underline{\mth{T}})$ is trivial on $\underline{\mth{S}}$, hence must be zero. Hence every set of positive roots of $\underline{\Psi}$ must contain some of the $\sigma(\alpha)$ but not all of then, which contradicts the second assertion. Hence $\alpha$ cannot exist and the centralizer of $\underline{\mth{S}}$ is reduced to $\underline{\mth{H}}$.

Conversely, assume the first assertion is true. Consider the restrictions of the elements of $\underline{\Phi}$ to $\underline{\mth{S}}$; they are all nonzero. Let $H$ be an hyperplane of the vector space $X^*(\underline{\mth{S}})\otimes{\mth{R}}$ which doesn't contain any one of these restrictions, let $C$ be a half-space delimited by $H$, and let $\underline{\Phi}^+$ be the set of elements of $\underline{\Phi}$ whose restriction belongs to $C$. The set $\underline{\Phi}^+$ is closed, and for each pair of opposite roors it contains one of the two roots but not both; by \cite[1, cor. 1 to prop. 20]{bou}, it is then a set of positive roots of $\underline{\Phi}$. Moreover, it is obviously $\Gamma$-stable, which completes the proof of the proposition. $\Box$

We'll say $\underline{\mth{G}}$ is quasi-split (as a group of parahoric type) if it satisfies the above conditions.

Obviously, if $\underline{\mth{G}}$ is quasi-split, then its reductive quotient is quasi-split too. The converse is not always true; in particular, contrary to what happens with reductive groups, a group of parahoric type defined over a finite field is not necessarily quasi-split. For example, assume $k$ is finite; let $F$ be a local field whose residual field is $k$, let $D$ be a division algebra of center $F$ and of finite degree $d>1$ over $F$, and let $G=GL_n(D)$; $G$ is the group of $F$-points of an inner form of the linear group $GL_{dn}$ which splits over a nonramified extension $F'$ of $F$ of degree $d$. If $K$ is a $Gal(F'/F)$-stable parahoric subgroup of $G$ and if ${\mth{G}}=K/K^h$, $h\geq 2$, then ${\mth{G}}$ is the group of $k$-points of a connected linear algebraic group of parahoric type $\underline{\mth{G}}$ defined over $k$ but not quasi-split.

In the sequel, we will assume $\underline{\mth{G}}$ is quasi-split and $\underline{\mth{B}}$ is $\Gamma$-stable, which implies $\underline{\Phi}^+$ is $\Gamma$-stable as well. Let ${\mth{G}}$ (resp. ${\mth{T}}$, ${\mth{H}}$, ${\mth{B}}$) be the set of $k$-points of $\underline{\mth{G}}$ (resp. $\underline{\mth{T}}$, $\underline{\mth{H}}$, $\underline{\mth{B}}$).

Let $\Phi$ (resp. $\Psi$) be the set of images of the elements of $\underline{\Phi}$ (resp. $\underline{\Psi}$) in $X^*(\underline{\mth{S}})$; we deduce from \cite[5.8]{bot} that $\Phi$ and $\Psi$ are both root systems.

The concave function $f_0$ is entirely determined by its values on the set $\underline{\Delta}^+$ of simple roots of $\underline{\Phi}^+$, which is $\Gamma$-stable. Since these values can be chosen arbitrarily, we can always assume that $f_0$ is constant on every $\Gamma$-orbit in $\underline{\Delta}^+$; we can easily check that it implies that $f_0$ is constant on every $\Gamma$-orbit in $\underline{\Phi}$. (It is in particular obviously the case for the concave function $f_0=0$ when $\underline{\mth{G}}$ is special.) Then $f_0$ reduces to a function on $\Phi$, which is clearly concave and which by a slight abuse of notation we will denote by $f_0$ too.

For every $\alpha\in\Phi$, let ${\mth{U}}_\alpha$ be the group of $k$-points of the subgroup of $\underline{\mth{G}}$ generated by the $\underline{\mth{U}}_\beta$, where $\beta$ runs through the elements of $\underline{\Phi}$ whose image in $\Phi$ is $\alpha$. We can similarly define, for every integer $i\geq f_0(\alpha)$, the subgroup ${\mth{U}}_{\alpha,i}$ as the group of $k$-points of the subgroup of $\underline{\mth{G}}$ generated by the $\underline{\mth{U}}_{\beta,i}$, with $\beta$ as above.

Note that it may happen that $2\alpha$ is also an element of $\Phi$. In this case, we see from the commutator relations that for every $i\geq f_0(\alpha)$, the group ${\mth{U}}_{2\alpha,2i}$ is contained in ${\mth{U}}_{\alpha,i}$. However, the group ${\mth{U}}_{\alpha}$ contains the full group ${\mth{U}}_{2\alpha}$ only if $f_0(2\alpha)=2f_0(\alpha)$.

Let $\underline{\mth{P}}$ be a $\Gamma$-stable generalized parahoric subgroup of $\underline{\mth{G}}$, and let ${\mth{P}}$ be the group of $k$-points of $\underline{\mth{P}}$. Let $f$ be the concave funcion associated to $\underline{\mth{P}}$; it is easy to see that $f$ is $\Gamma$-stable, and we will also denote by $f$ the corresponding function on $\Phi$. We have:

\begin{prop}
The group ${\mth{P}}$ is the subgroup of ${\mth{G}}$ generated by ${\mth{H}}$ and the ${\mth{U}}_{\alpha,f(\alpha)}$, $\alpha\in\Phi$.
\end{prop}

Obviously, ${\mth{H}}$ and the ${\mth{U}}_{\alpha,f(\alpha)}$ are contained in ${\mth{P}}$; we'll prove the other inclusion. Assume first $\underline{\mth{P}}$ is solvable, and let $p$ be an element of ${\mth{P}}$. We then deduce from the commutator relations the following decomposition:
\[p=c\prod_{\alpha\in\underline{\Phi}}u_\alpha=c\prod_{\beta\in\Phi}(\prod_{\alpha\mapsto\beta}u_\alpha);\]
where $c$ is an element of $\underline{\mth{H}}$ and the $u_\alpha$ are elements of $U_\alpha$. Moreover, since $p\in{\mth{P}}$, we also have for every $\sigma\in\Gamma$:
\[p=\sigma(p)=\sigma(c)\prod_{\beta\in\Phi}(\prod_{\alpha\mapsto\beta}\sigma(u_\alpha)).\]
For every $\sigma,\alpha$, $\sigma(u_\alpha)$ belongs to $\underline{\mth{U}}_{\sigma(\alpha)}$. Let $\beta$ be a nonmultipliable element of $\Phi$; for every $\alpha,\alpha'$ whose image in $\Phi$ is $\beta$, $\alpha+\alpha'$ is not a root (if it was, its image in $\Phi$ would be $2\beta$, which contradicts our hypothesis), hence $u_\alpha$ and $u_\alpha'$ commute. We can then reorder the terms in the second product in such a way that the order of the elements of $\underline{\mth{U}}_\alpha$, with $\alpha$ having $\beta$ as an image in $\Phi$, is the same as in the first one, and by unicity of the decomposition, we obtain $\sigma(u_\alpha)=u_{\sigma(\alpha)}$ for every such $\alpha$ and every $\sigma$. Hence $\prod_{\alpha\mapsto\beta}u_\alpha$ is $\Gamma$-fixed and belongs to ${\mth{U}}_{\beta,f(\beta)}$.

Now let's look at the multipliable roots in $\Phi$. When they exist, they generate a subsystem of type $BC_1^r$. We are thus reduced to prove the assertion when $\Phi$ is of type $BC_1$ and $p$ is of the form $c\prod_{\alpha\mapsto\pm\beta}u_\alpha$, where $\beta$ is a multipliable root in $\Phi$, in which case it is an immediate consequence of the existence and unicity of the Iwahori decomposiion in the compact subgroup of $K$ whose image in $\underline{\mth{G}}$ is $\underline{\mth{P}}$.

Now look at the general case. The previous discussion shows that ${\mth{H}}$ and the ${\mth{U}}_\alpha$ generate a Zariski-open subgroup of ${\mth{P}}$, which is enough to prove the proposition when $k$ is infinite. Assume then $k$ finite, and let $\mathbf{F}$ be the Frobenius map on $\underline{\mth{G}}$ such that ${\mth{G}}$ is the set of $\mathbf{F}$-fixed points of $\underline{\mth{G}}$. Let $\underline{\mth{B}}$ be a pseudo-Borel subgroup of $\underline{\mth{G}}$ contained in $\underline{\mth{P}}$, let $\underline{\mth{B}}_0$ be the unique Borel subgroup of $\underline{\mth{G}}$ containing $\underline{\mth{B}}$, and let $p$ be an element of ${\mth{P}}$. Since $\underline{\mth{B}}_0\cap\underline{\mth{P}}$ is then a Borel subgroup of $\underline{\mth{P}}$, by Bruhat decomposition we can write $p$ as $p=bnb'$, where $b,b'$ are elements of $\underline{\mth{B}}_0\cap\underline{\mth{P}}$ and $n$ is an element of the normalizer of $\underline{\mth{T}}$ in $\underline{\mth{P}}$, which can always be chosen $\mathbf{F}$-fixed and thus, as a well-known result, is a product of elements of ${\mth{H}}$ and the ${\mth{U}}_{\alpha,f(\alpha)}$, since it is an element of a $w$ which belongs to the relative Weyl group of ${\mth{P}}$ relatively to ${\mth{T}}$. Hence we have $bnb'=\mathbf{F}(b)n\mathbf{F}(b')$, which can be rewritten as:
\[b^{-1}\mathbf{F}(b)=nb'\mathbf{F}(b')^{-1}n^{-1}.\]
Since both sides of the above equality are elements of the connected group $\underline{\mth{B}}_0\cap n\underline{\mth{B}}_0n^{-1}\cap\underline{\mth{P}}$, by Lang's theorem there exists an element $g$ of that group such that $b^{-1}\mathbf{F}(b)=g^{-1}\mathbf{F}(g)$, hence $bg^{-1}$ is $\mathbf{F}$-fixed. Hence we have:
\[p=bnb'=bg^{-1}n(n^{-1}gnb'),\]
and applying the previous case to $bg^{-1}$ and $n^{-1}gnb'$ gives us the result. $\Box$

\section{Double classes}\label{dc}

From now on we will assume that $F$ $F_{nr}$, $\underline{G}$, $K$ and $\underline{\mth{G}}$ are defined as in the proposition \ref{gfparah}. Let $\mathcal{O}_{nr}$ be the ring of integers of $F_{nr}$, and let $\varpi$ be an uniformizer of $F_{nr}$; for every $\alpha\in\underline{\Phi}$, we can choose a surjective morphism $\phi_\alpha$ from $\varpi^{f_0(\alpha)}\mathcal{O}_{nr}$ to $\underline{\mth{U}}_{nr}$. For every $i\geq f_0(\alpha)$, we then have $\phi_\alpha(\varpi^i\mathcal{O}_{nr})=\underline{\mth{U}}_{\alpha,i}$; in particular, the kernel of $\phi_\alpha$ is $\varpi^{h-f_0(-\alpha)}$.

Similarly, if $T_{nr}$ and $\underline{T}$ is defined as in the proof of the proposition \ref{gfparah}, we have a canonical isomorphism between their groups of cocharacters $X_*(T_{nr})$ and $X_*(\underline{\mth{T}})$ (resp. between their groups of characers $X^*(T_{nr})$ and $X^*(\underline{\mth{T}})$), from which we deduce that:
\begin{itemize}
\item every one-parameter subgroup $\xi$ of $\underline{\mth{T}}$ can be canonically extended to a $h$-dimensional subgroup of $\underline{\mth{H}}$ isomorphic to $\mathcal{O}_{nr}^*/(1+\varpi^h\mathcal{O}_{nr})$, which we will also call $\xi$;
\item every character $\chi$ of $\underline{\mth{T}}$ can be canonically lifted to a morphism from $\underline{\mth{H}}$ to $\mathcal{O}_{nr}^*/(1+\varpi^h\mathcal{O}_{nr})$, which we will also call $\chi$.
\end{itemize}

The usual properties of characters and cocharacters (like the existence of a pairing between them for example) work well with these extensions and liftings since in that case, they can be directly deduced from the similar properties of characters and cocharacters of $T_{nr}$.

Let ${\mth{G}}$ be the group of $k$-points of $\underline{\mth{G}}$. In this section we will establish some facts about double classes of $\underline{\mth{G}}$ and ${\mth{G}}$ modulo their generalized pseudo-parabolic subgroups. Note that contrary to the reductive case, a pseudo-Borel subgroup and the normalizer of a maximal torus do not form a $BN$-pair, so the theory of Tits systems cannot be applied here. In fact, in most cases, the set of such double classes in $\underline{\mth{G}}$ is infinite. So we'll have to find other ways to deal with them.

Let $\underline{\Phi}$ be the root system of $\underline{\mth{G}}$ relatively to $\underline{\mth{T}}$, let $\underline{\mth{B}}$ be a pseudo-Borel subgroup of $\underline{\mth{G}}$ and $\underline{\Phi}^+$ (resp. $\underline{\Phi}^-$) the corresponding set of positive (resp. negative) roots of $\underline{\Phi}$, let $\underline{\mth{B}}_0$ be the unique Borel subgroup of $\underline{\mth{G}}$ containing $\underline{\mth{B}}$ and let $\underline{\mth{U}}_0$ be the unipotent radical of $\underline{\mth{B}}_0$. Let $W$ be the Weyl group of $\underline{\mth{G}}$ w.r.t. $\underline{\mth{T}}$. We have the Bruhat decomposition:
\[\underline{\mth{G}}=\bigsqcup_{w\in W}\underline{\mth{B}}_0w\underline{\mth{B}}_0=\bigsqcup_{w\in W}\underline{\mth{U}}_0w\underline{\mth{B}}_0\]
Let $\underline{\mth{U}}$ (resp. $\underline{\mth{U}}^-$) be the subgroup of $\underline{\mth{G}}$ generated by the $\underline{\mth{U}}_\alpha$ where $\alpha$ is positive (resp. negative) w.r.t $\underline{\mth{B}}$; we also have the following Iwahori decompositions:
\[\underline{\mth{B}}_0=\underline{\mth{H}}\underline{\mth{U}}(\underline{\mth{U}}^-\cap\underline{\mth{B}}_0).\]
\[\underline{\mth{U}}_0=R_u(\underline{\mth{H}})\underline{\mth{U}}(\underline{\mth{U}}^-\cap\underline{\mth{B}}_0).\]
Moreover, for every $w\in W$, $w^{-1}(\underline{\mth{U}}^-\cap\underline{\mth{B}}_0)w$ and $w^{-1}\underline{\mth{H}}w=\underline{\mth{H}}$ are contained in $\underline{\mth{B}}_0$; we finally obtain:
\[\underline{\mth{G}}=\bigsqcup_{w\in W}\underline{\mth{B}}w(\underline{\mth{U}}^-\cap\underline{\mth{B}}_0)\underline{\mth{B}}.\]
\[=\bigsqcup_{w\in W}\underline{\mth{U}}w(\underline{\mth{U}}^-\cap\underline{\mth{B}}_0)\underline{\mth{U}},\]

Hence every double class modulo $\underline{\mth{B}}$ (resp. modulo $\underline{\mth{U}}$ on the left and $\underline{\mth{B}}$ on the right) contains elements of the form $nu$, with $n$ belongs to a given system of representatives of the elements of $W$ and $u$ is an element of $\underline{\mth{U}}^-\cap\underline{\mth{B}}_0$. The problem of determining precisely those double classes in the most general setting seems to be difficult; we'll examine some simple cases.

First we will look at the case where $\underline{\Phi}$ is of rank $1$. Let $\alpha$ be the unique element of $\underline{\Phi}^-$; the generalized pseudo-parabolic subgroups of $\underline{\mth{G}}$ are the $\underline{\mth{P}}_i$ defined by the concave functions $f_i$ on $\underline{\Phi}=\{\alpha,-\alpha\}$ such that $f_i(\alpha)=i$, $i\geq f_0(\alpha)$; we have $\underline{\mth{G}}=\underline{\mth{P}}_{f_0(\alpha)}$, $\underline{\mth{B}}_0=\underline{\mth{P}}_{1-f_0(-\alpha)}$, $\underline{\mth{B}}=\underline{\mth{P}}_{h-f_0(\alpha)}$. Consider the Bruhat decomposition of $\underline{\mth{G}}$: if $\Psi$ is of rank $1$, the decomposition is:
\[\underline{\mth{G}}=\underline{\mth{B}}_0\sqcup\underline{\mth{B}}_0w\underline{\mth{B}}_0,\]
where $w$ is the nontrivial element of the Weyl group of $\underline{\Psi}$, and we'll examine which $\underline{\mth{B}}$-double classes are contained into each one of those two big double classes. If $\Psi$ is of rank $0$, then $\underline{\mth{G}}=\underline{\mth{B}}_0$ and its double classes are determined the same way.

Assume then that $\Psi$ is of rank $1$. We have the Iwahori decomposition:
\[\underline{\mth{P}}_i=\underline{\mth{B}}\underline{\mth{U}}_{-\alpha,i}\]
for every $i>-f_0(-\alpha)$. We deduce from this that we have:
\[\underline{\mth{B}}_0w=\underline{\mth{B}}\underline{\mth{U}}_{-\alpha,1-f_0(-\alpha)}w=\underline{\mth{B}}w\underline{\mth{U}}_{\alpha,1-f_0(\alpha)},\]
hence:
\[\underline{\mth{B}}_0w\underline{\mth{B}}_0=\underline{\mth{B}}w\underline{\mth{B}}_0.\]
Applying the Iwahori decomposition to $\underline{\mth{B}}_0$ once again we finally obtain:
\[\underline{\mth{B}}_0w\underline{\mth{B}}_0=\underline{\mth{B}}w\underline{\mth{B}}.\]
Hence $\underline{\mth{B}}_0w\underline{\mth{B}}_0$ consists of only one $\underline{\mth{B}}$-double class.
Since $\underline{\mth{H}}$ is normalized by $w$, we even have:
\[\underline{\mth{B}}_0w\underline{\mth{B}}_0=\underline{\mth{U}}w\underline{\mth{B}},\]
hence $\underline{\mth{B}}_0w\underline{\mth{B}}_0$ also consists of only one double class mod $\underline{\mth{U}}$ on the left and $\underline{\mth{B}}$ on the right.

Now look at the double classes contained in $\underline{\mth{B}}_0$. From the previous discussion we see that every such double class contains elements of $\underline{\mth{U}}_{-\alpha,1-f_0(-\alpha)}$. Let $u$ be such an element, and let $v(u)$ be its valuation. It can be easily checked that two elements $u,u'$ of $\underline{\mth{U}}_{-\alpha}$ belong to the same double class mod $\underline{\mth{B}}$
 if and only if there exists an element $t$ of $\underline{\mth{H}}$ such that $t^{-1}ut=u'$, which in particular implies that $u$ and $u'$ have the same valuation. When $\underline{\mth{G}}$ is adjoint modulo its center, this implication is an equivalence for rhe $\underline{\mth{B}}$-double classes; this is not true in general.

We'll now examine the double classes of $\underline{\mth{B}}_0$ mod $\underline{\mth{U}}$ on the left and $\underline{\mth{B}}$ on the right. For a given $i>-f_0(-\alpha)$, we can deduce from the commutator relations that for every $u\in\underline{\mth{U}}_{\alpha,i}$, $[u,\underline{\mth{U}}_{-\alpha}]$ is contained in $\underline{\mth{H}}\underline{\mth{U}}_{\alpha,2i+f_0(-\alpha)}$ and that its canonical projection on $\underline{\mth{U}}_{\alpha,2i+f_0(-\alpha)}$ is surjective; hence $u,u'\in\underline{\mth{U}}_{\alpha,i}$ belong to the same double class if and only if $u'\in u\underline{\mth{U}}_{\alpha,2i+f_0(-\alpha)}$.

Now we go back to the general case. Let $\underline{\Delta}$ be the set of negative simple roots in $\underline{\Phi}$; our goal will be to determine, for every $I\subset\underline{\Delta}$, the double classes of the generalized pseudo-parabolic subgroup $\underline{\mth{B}}_I$ of $\underline{\mth{G}}$ mod $\underline{\mth{B}}$; we will in fact prove a slightly more general result. First we have the following lemma:

\begin{lemme}
Let $\underline{\mth{P}}=\underline{\mth{M}}\underline{\mth{U}}_{\underline{\mth{P}}}$ be a pseudo-parabolic subgroup of $\underline{\mth{G}}$ and let $\underline{\mth{L}}_{\underline{\mth{M}},i}$ $i=1,2$, be subgroups of $\underline{\mth{M}}$; for each $i$, set $\underline{\mth{L}}_i=\underline{\mth{L}}_{\underline{\mth{M}},i}\underline{\mth{U}}_{\underline{\mth{P}}}$. Then the set of double classes of $\underline{\mth{P}}$ modulo $\underline{\mth{L}}_1$ on the left and $\underline{\mth{L}}_2$ on the right is in canonical $1-1$ correspondance with the set of double classes of $\underline{\mth{M}}$ modulo $\underline{\mth{L}}_{\underline{\mth{M}},1}$ on the left and $\underline{\mth{L}}_{\underline{\mth{M}},2}$ on the right.
\end{lemme}

Since $\underline{\mth{U}}_{\underline{\mth{P}}}$ is normal in $\underline{\mth{P}}$, the proof is obvious. $\Box$

We can immediately deduce from this lemma and the discussion about rank 1 the double classes of the $\underline{\mth{B}}_\alpha$ for every $\alpha\in\underline{\Delta}$. We now have to introduce a few notions.

Let $f$ be a concave function on $\underline{\Phi}$. For every subset $S$ of $\underline{\Phi}^-$, let $\underline{\mth{U}}_{f,S}$ be the subgroup of $\underline{\mth{G}}$ generated by the $\underline{\mth{U}}_{\alpha,f(\alpha)}$, $\alpha\in S$. We'll say that $S$ is $f$-closed if $\underline{\mth{U}}_{f,S}$ is precisely the product of the $\underline{\mth{U}}_{\alpha,f(\alpha)}$, $\alpha\in S$.

Obviously $S$ is $f$-closed if and only if for every $\alpha,\beta\in S$ such that $\alpha+\beta$ is a root not belonging to $S$, $\underline{\mth{U}}_{\alpha+\beta,f(\alpha+\beta)}$ is trivial. In particular every closed subset of $\underline{\Phi}^-$ is trivially $f$-closed.

We'll also say that $S$ is complete if for every $\alpha\in S$ and every $\beta\in\underline{\Phi}^-$ such that $\beta>\alpha$, $\beta\in S$. If $S$ is complete, it must in particular contain every simple root showing up in the decomposition of $\alpha$ into a sum of elements of $\underline{\Delta}$, which implies that $S$ is contained in the parabolic subsystem of $\underline{\Phi}$ generated by $\underline{\Phi}^+$ and the negative simple roots it contains.

Let $\underline{\mth{P}}$ be a generalized parabolic subgroup of $\underline{\mth{G}}$ containing $\underline{\mth{B}}$, let $f$ be the corresponding concave function on $\underline{\Phi}$ and let $S$ be the subset of elements $\alpha$ of $\underline{\Phi}^-$ such that $U_{\alpha,f(\alpha)}$ is nontrivial; $S$ is clearly $f$-closed and complete.

Assume now there exists two complete subsets $S_1$ and $S_2$ of $S$ such that $S_1\cup S_2=S$. Then for each $i$, every element in $S_i$ is a linear combination of the simple roots contained in $S_i$, which implies that the sum of two elements of $S_i$ is either in $S_i$ or not in $S$ at all; hence the $S_i$ are $f$-closed as well. Let $\underline{\mth{P}}_i$ be the pseudo-parabolic subgroup generated by $\underline{\mth{B}}$ and the root subgroups associated to the simple roots in $S_i$; $\underline{\mth{P}}$ is then contained in the product $\underline{\mth{P}}_1\underline{\mth{P}}_2$ (which is not necessarily a group). Moreover, the sum of an element of $S_1$ and an element of $S_2$ cannot be in $S$, hence $\underline{\mth{U}}_{f,S_1}$ and $\underline{\mth{U}}_{f,S_2}$ commute.

Let $\underline{\Psi}_f$ be the set of elements $\alpha\in\underline{\Psi}$ such that $f(\alpha)+f(-\alpha)=0$, and let $\underline{\Psi}_1$ (resp. $\underline{\Psi}_2$): be the intersection of $\underline{\Psi}_f$ with the root system of $\underline{\mth{P}}_1$ (resp. $\underline{\mth{P}}_2$). Since for $i=1,2$, $\underline{\Psi}_i\cap\underline{\Phi}^-$ is contained in $S_i$, they are disjoint; hence $\underline{\Psi_1}$ and $\underline{\Psi}_2$ are disjoint too; moreover, their union is $\underline{\Psi}_f$. We deduce from this that the Weyl group of $\underline{\mth{P}}_f$ w.r.t. $\underline{\mth{T}}$ is the direct product of the Weyl groups of the $\underline{\mth{P}}_i$.

We finally obtain that there exists a system of representatives of the double classes of $\underline{\mth{P}}_f$ modulo $\underline{\mth{B}}$ (resp. modulo $\underline{\mth{B}}$ on the left and $\underline{\mth{U}}$ on the right) which is made of elements of the form $n_1u_1n_2u_2$, where for $i=1,2$, $n_i$ belongs to a fixed set of representatives of the elements of the Weyl group of $\Psi_i$ and $u_i\in\underline{\mth{U}}_{f,S_i}$.

We have the following result::
\begin{prop}\label{dcprod}
For $i=1,2$, let $R_i$ be a system of representatives of the double classes of $\underline{\mth{P}}_f\cap \underline{\mth{P}}_i$ mod $\underline{\mth{U}}$ on the left and $\underline{\mth{B}}$ on the right; assume every element of $R_i$ is of the form $n_iu_i$, where $n_i$ and $u_i$ are as above. Then $R_1R_2$ is a system of representatives of the double classes of $\underline{\mth{P}}_f$ mod $\underline{\mth{U}}$ on the left and $\underline{\mth{B}}$ on the right.
\end{prop}

Let $g$ be an element of $\underline{\mth{P}}_f$. Since $g$ belongs to the product set $\underline{\mth{P}}_1\underline{\mth{P}}_2$ we can write it as $g=bn_1u_1b'n_2u_2b''$, where $b,b',b''\in\underline{\mth{B}}$ (and we can even assume $b,b'\in\underline{\mth{U}}$) and $n_i,u_i$ are as above.

Write $b'=b'_1b'_2$, where for each $i$, $b'_i$ belongs to the product $\underline{\mth{U}}_i$ of the root subgroups associated to the positive roots which are not linear combinations of the simple roots contained in $S_i$; this is always possible since those two sets of simple roots are disjoint. The conjugate of $b'_i$ by any representative of any element of the Weyl group of $\Psi_i$ is then still an element of $\underline{\mth{U}}_i$, and the commutators $[b'_1{}^{-1},u_1]$ and $[n_2^{-1}b'_2{}^{-1}n_2,u_2]$ also belong to $\underline{\mth{U}}_i$; we thus obtain:
\[g=bn_1b'_1[b'_1{}^{-1},u_1]u_1n_2(n_2^{-1}b'_2n_2)u_2b''\]
\[=b(n_1b'_1[b'_1{}^{-1},u_1]n_1^{-1})n_1u_1n_2u_2[n_2^{-1}b'_2{}^{-1}n_2,u_2]^{-1}b''\in\underline{\mth{U}}n_1u_1n_2u_2\underline{\mth{B}}.\]
Now assume we have $g\in\underline{\mth{U}}n'_1u'_1n'_2u'_2\underline{\mth{B}}$, where $n'_iu'_i\in R_i$ for $i=1,2$. There exists then $b_0\in\underline{\mth{U}}$ and $b'_0\in\underline{\mth{B}}$ such that $g'=u'_1{}^{-1}n'_1{}^{-1}b_0n_1u_1=n'_2u'_2b'_0u_2^{-1}n_2^{-1}$, and since $g'$ must then belong to $\underline{\mth{P}}_f\cap\underline{\mth{P}}_1\cap\underline{\mth{P}}_2=\underline{\mth{B}}$, this is only possible if $n_iu_i=n'_iu'_i$ for $i=1,2$. This proves the result. $\Box$

\begin{cor}
Let $S=\bigcup_{i=1}^nS_i$ be is a partition of $S$ into $n$ complete subsets, and for every $i$, let $\underline{\mth{P}}_i$ be the pseudo-parabolic subgroup of $\underline{\mth{G}}$ associated to $S_i$ and $R_i$ defined as above. Then $\prod_{i=1}^nR_i$ is a system of representatives of the double classes of $\underline{\mth{P}}_f$ mod $\underline{\mth{U}}$ on the left and $\underline{\mth{B}}$ on the right.
\end{cor}

Thos corollary comes from the proposition and an obvious induction, $\Box$

In particular, when $S$ is contained in $\underline{\Delta}$, it is complete if and only if for every $\alpha,\beta\in S$ such that $\alpha+\beta$ is a root, $\underline{\mth{U}}_{\alpha+\beta,f(\alpha+\beta)}$ is trivial, which is in particular true if $f(\alpha)=f_{\underline{\mth{B}}}(\alpha)-1$ for every $\alpha\in S$; moreover, all singletons contained in $\underline{\Delta}$ are trivially complete subsets of $\underline{\Phi}^-$. We thus obtain:

\begin{prop}
Let $f$ be a concave function of $\underline{\Phi}$ such that $f(\alpha)=f_{\underline{\mth{B}}}(\alpha)$ for every $\alpha\not\in\underline{\Delta}$, and let $\underline{\mth{P}}_f$ be the corresponding generalized pseudo-parabolic subgroup of $\underline{\mth{G}}$. Then the double classes of $\underline{\mth{P}}_f$ mod $\underline{\mth{U}}$ on the left and $\underline{\mth{B}}$ on the right admit as a set of representatives the set of products $\prod_{\alpha\in\underline{\Delta}}g_\alpha$, where for every $\alpha\in\underline{\Delta}$:
\begin{itemize}
\item if $f(\alpha)+f(-\alpha)>0$, $g_\alpha$ runs through the union of sets of representatives of the groups $\underline{\mth{U}}_{\alpha,i}$, $i\geq f(\alpha)$, respectively mod  $\underline{\mth{U}}_{\alpha,2i+f_0(-\alpha)}$;
\item if $f(\alpha)+f(-\alpha)=0$, $g_\alpha$ runs through the union of sets of representatives of the groups $\underline{\mth{U}}_{\alpha,i}$, $i>f(\alpha)$, respectively mod  $\underline{\mth{U}}_{\alpha,2i+f_0(-\alpha)}$, and of the singleton $\{n\}$, where $n$ is a representative of the nontrivial element of the Weyl group of the pseudo-parabolic group $\underline{\mth{P}}$ of $\underline{\mth{G}}$ generated by $\underline{\mth{B}}$ and $\underline{\mth{U}}_\alpha$.
\end{itemize}
\end{prop}

Moreover, as in the case of rank $1$, two elements $g=\prod_{\alpha\in\underline{\Delta}}g_\alpha$ and $g'=\prod_{\alpha\in\underline{\Delta}}g'_\alpha$ of this set of representatives are in the same double class mod $\underline{\mth{B}}$ on both sides if and only if there exists an element $t\in\underline{\mth{H}}$ such that $t^{-1}gt=g'$, which implies that for every $\alpha$, $g_\alpha$ and $g'_\alpha$ have the same valuation, and this implication is an equivalence when $\underline{\mth{G}}$ is adjoint modulo its center. In particular we obtain:

\begin{prop}\label{dcadj}
Assume $\underline{\mth{G}}$ is adjoint modulo its center and $\underline{\mth{U}}_{\alpha,f(\alpha)}$ is of dimension $1$ for every $\alpha\in\underline{\Delta}$. Then the double classes of $\underline{\mth{P}}_f$ mod $\underline{\mth{B}}$ are in $1-1$ correspondence with the subsets of $\underline{\Delta}$.
\end{prop}

Now we will consider the double classes of $\underline{\mth{G}}$ which are not contained in $\underline{\mth{B}}_{\underline{\Delta}}$. We'll only prove one useful result about them.

From now on we will assume that $\underline{\mth{G}}$ is special.

Let $R_u(\underline{\mth{B}})$ be the unipotent radical of $\underline{\mth{B}}$. We have the following result:

\begin{prop}\label{clbi}
Assume the following conditions on $\underline{\Phi}$ and $p$ are true:
\begin{itemize}
\item $p\neq 2$;
\item if $\underline{\Phi}$ has a least one component of type $A_n$, then $p$ doesn't divide the adjoint index of $\underline{\mth{G}}$;
\item if $\underline{\Phi}$ has a least one component of type $E_6$, $E_7$ or $F_4$, then $p\neq 3$;
\item if $\underline{\Phi}$ has a least one component of type $E_8$, then $p>5$.
\end{itemize}
Let $g$ be an element of $\underline{\mth{G}}$ which doesn't belong to $\underline{\mth{B}}_{\underline{\Delta}}$ ;there exists $\alpha\in\underline{\Delta}$, depending only on the double class of $g$ modulo $\underline{\mth{B}}_{\underline{\Delta}}$, such that $\underline{\mth{B}}g\underline{\mth{B}}=\underline{\mth{B}}g\underline{\mth{B}}_\alpha$ (resp. $R_u(\underline{\mth{B}})g\underline{\mth{B}}=R_u(\underline{\mth{B}})g\underline{\mth{B}}_\alpha$).
\end{prop}

For the (quite long) proof of this proposition, see the section $5$.

\subsection{Double classes of rational points}

Now we'll examine the double classes of ${\mth{G}}$ modulo the group of $k$-points of a pseudo-Borel subgroup. In this section, we will assume $k$ is finite; let $\mathbf{F}$ be the Frobenius map on $\underline{\mth{G}}$ such that ${\mth{G}}$ is the group of $\mathbf{F}$-fixed points of $\underline{\mth{G}}$. We will prove similar results as in the case of $\underline{\mth{G}}$, using the following proposition:

\begin{prop}\label{clrel}
Let $\underline{\mth{B}},\underline{\mth{B}}'$ be two $\mathbf{F}$-stable pseudo-Borel subgroups of $\underline{\mth{G}}$ and let ${\mth{B}},{\mth{B}}'$ be their respective groups of $k$-points. Let $R_u(\underline{\mth{B}})$ be the unipotent radical of $\underline{\mth{B}}$, and let $R_u({\mth{B}})$ be its group of $k$-points. There is a canonical $1-1$ correspondance between the $\mathbf{F}$-stable elements in $R_u(\underline{\mth{B}})\backslash\underline{\mth{G}}/\underline{\mth{B}}'$ and $R_u({\mth{B}})\backslash{\mth{G}}/{\mth{B}}'$ given by: the image of a $\mathbf{F}$-stable double class is its set of $k$-points.
\end{prop}

It is obvious that every element of $R_u(\underline{\mth{B}})\backslash\underline{\mth{G}}/\underline{\mth{B}}$ containing at least one element of ${\mth{G}}$ is $\mathbf{F}$-stable. Conversely, let $R_u(\underline{\mth{B}})g\underline{\mth{B}}'$ be a $\mathbf{F}$-stable double class; we will prove that its set of $k$-points is nonempty. Since it is $\mathbf{F}$-stable, it then also contains $\mathbf{F}(g)$, hence there exist $b\in R_u(\underline{\mth{B}})$ and $u\in\underline{\mth{B}}'$ such that $bgu=\mathbf{F}(g)$.

Since $R_u(\underline{\mth{B}})$ is connected, according to Lang's theorem, there exists $h\in R_u(\underline{\mth{B}})$ such that $b=\mathbf{F}(h)^{-1}h$, and similarly there exists $h'\in\underline{\mth{B}}'$ such that $u=h'\mathbf{F}(h')^{-1}$. We thus obtain $hgh'=\mathbf{F}(h)\mathbf{F}(g)\mathbf{F}(h')$, hence $hgh'$ is an element of ${\mth{G}}$ belonging to $R_u(\underline{\mth{B}})u\underline{\mth{B}}'$.

It remains to prove that the set of $k$-points of $R_u(\underline{\mth{B}})g\underline{\mth{B}}'$ contains only one element of $R_u({\mth{B}})\backslash{\mth{G}}/{\mth{B}}'$. Assume $g\in{\mth{G}}$ and let $g'$ be another element of ${\mth{G}}$ belonging to $R_u(\underline{\mth{B}})g\underline{\mth{B}}'$; write $g=bg'u$. $b\in R_u(\underline{\mth{B}})$, $u\in\underline{\mth{B}}'$. We then also have $g=\mathbf{F}(b)g'\mathbf{F}(u)$, hence:
\[g^{-1}\mathbf{F}(b)b'^{-1}g=\mathbf{F}(u)^{-1}u\in\underline{\mth{B}}'\cap g^{-1}R_u(\underline{\mth{B}})g.\]
Assume $\underline{\mth{B}}'\cap g^{-1}R_u(\underline{\mth{B}})g$ is connected. Then by Lang's theorem, there exists $h\in\underline{\mth{B}}'\cap g^{-1}R_u(\underline{\mth{B}})g$ such that $\mathbf{F}(u)^{-1}u=\mathbf{F}(h)h^{-1}$, hence $uh=\mathbf{F}(uh)\in{\mth{B}}'$. We also have $g^{-1}\mathbf{F}(b)b^{-1}g=\mathbf{F}(h)h^{-1}$, hence $b^{-1}ghg^{-1}=\mathbf{F}(b^{-1}ghg^{-1})\in R_u({\mth{B}})$. We finally obtain:
\[g'=b^{-1}gu^{-1}=(u^{-1}ghg^{-1})g(uh)^{-1}\in R_u({\mth{B}})g{\mth{B}}'.\]

We now have to prove that $\underline{\mth{B}}'\cap g^{-1}R_u(\underline{\mth{B}})g$ is connected. By eventually replacing $\underline{\mth{B}}'$ by a conjugate, we can make the following assumptions:
\begin{itemize}
\item there is a maximal torus $\underline{\mth{T}}$ of $\underline{\mth{G}}$ contained in both $\underline{\mth{B}}$ and $\underline{\mth{B}}'$;
\item $g$ is contained in the product of the $\underline{\mth{U}}_\alpha$, where the $\alpha$ are the roots of $\underline{\mth{G}}$ w.r.t. $\underline{\mth{T}}$ which are negative w.r.t. both $\underline{\mth{B}}$ and $\underline{\mth{B}}'$.
\end{itemize}

For every $i\in\{1,\dots,h\}$, let $\underline{\mth{G}}^i$ be the normal subgroup of $\underline{\mth{G}}$ which is the image in $\underline{\mth{G}}$ of the $i$-th congruence subgroup $K^i$ of $K$. It is enough to prove that for every $i\in\{0,\dots,h-1\}$, the quotient $(\underline{\mth{B}}'\cap g^{-1}R_u(\underline{\mth{B}})g\cap\underline{\mth{G}}^i)/(\underline{\mth{B}}'\cap g^{-1}R_u(\underline{\mth{B}})g\cap\underline{\mth{G}}^{i+1})$ is connected.
For every $u=u_{\mth{H}}\prod_{\alpha}u_\alpha(x_\alpha)\in\underline{\mth{B}}'\cap\underline{\mth{G}}^i$, with $u_{\mth{H}}\in{\mth{H}}$, the class of $u$ belongs to that group if and only if $gu'g^{-1}\in R_u(\underline{\mth{B}})$ for some $u'$ in that class. But then $u_{\mth{H}}$ belongs to $R_u({\mth{H}})$, and by the commutator relations, this condition is linear, which is enough to prove that our group is connected and finishes the proof of the proposition. $\Box$

\begin{cor}
Assume $\underline{\mth{G}}$ is adjoint modulo its center. There is then a canonical $1-1$ correspondance between the $\mathbf{F}$-stable elements in $\underline{\mth{B}}\backslash\underline{\mth{G}}/\underline{\mth{B}}'$ and ${\mth{B}}\backslash{\mth{G}}/{\mth{B}}'$ given by: the image of a $\mathbf{F}$-stable double class is its set of $k$-points.
\end{cor}

Let $\Delta$ be the set of negative simple roots of $\Phi$; for every $\alpha\in\Delta$ (resp for every $I\subset\Delta$, define ${\mth{B}}_\alpha$ 'resp. ${\mth{B}}_I$) in a similar way as in the absolute case. For every $\alpha$ (resp. $I$), let $J$ be the set of elements of $\underline{\Delta}$ whose image in $\Phi$ is $\alpha$ (resp. belongs to $I$); the group ${\mth{B}}_\alpha$ (resp. ${\mth{B}}_J$) is the group of $\mathbf{F}$-fixed points of $\underline{\mth{B}}_J$.

\begin{prop}\label{clbirel}
Assume $\underline{\Phi}$ and $p$ satisfy the same conditions as in the proposition \ref{clbi}. Let $g$ be an element of ${\mth{G}}$ which doesn't belong to ${\mth{B}}_{\Delta}$ ;there exists $\alpha\in\Delta$, depending only on the double class of $g$ modulo ${\mth{B}}_{\Delta}$, such that ${\mth{B}}g{\mth{B}}={\mth{B}}g{\mth{B}}_\alpha$ (resp. $R_u({\mth{B}})g{\mth{B}}=R_u({\mth{B}})g{\mth{B}}_\alpha$).
\end{prop}

According to the proposition \ref{clrel}, $R_u({\mth{B}})g{\mth{B}}$ is the group of $\mathbf{F}$-stable poinrs of $R_u(\underline{\mth{B}})g\underline{\mth{B}}$, and according to the proposition \ref{clbi}, there exists $\beta\in\underline{\Delta}$, depending only of the double class of $g$ modulo $\underline{\mth{B}}_{\underline{\Delta}}$, which contains the double class of $g$ modulo ${\mth{B}}_{\Delta}$, such that $R_u(\underline{\mth{B}})g\underline{\mth{B}}=R_u(\underline{\mth{B}})g\underline{\mth{B}}_\beta$. Since $R_u(\underline{\mth{B}})g\underline{\mth{B}}$ is $\mathbf{F}$-stable, it remains true if we replace $\beta$ by $\mathbf{F}^i(\beta)$ for every integer $i$; if $I$ is the set of such roots, we then have $R_u(\underline{\mth{B}})g\underline{\mth{B}}=R_u(\underline{\mth{B}})g\underline{\mth{B}}_I$. On the other hand, the elements of $I$ are precisely the elements of $\underline{\Delta}$ whose image in $\Phi$ is $\alpha$, hence $\underline{\mth{B}}_I$ contains ${\mth{B}}_\alpha$, which proves that $R_u({\mth{B}})g{\mth{B}}=R_u({\mth{B}})g{\mth{B}}_\alpha$. We deduce from this that ${\mth{B}}g{\mth{B}}$ is stable by ${\mth{B}}_\alpha$ on the right, hence ${\mth{B}}g{\mth{B}}={\mth{B}}g{\mth{B}}_\alpha$. $\Box$

\begin{prop}\label{dcadjrel}
Assume $\underline{\mth{G}}$ is adjoint modulo its center. Then the double classes of ${\mth{B}}_\Delta$ mod ${\mth{B}}$ are in $1-1$ correspondence with the subsets of $\Delta$.
\end{prop}

According to the proposition \ref{dcadj}, the double classes of $\underline{\mth{B}}_\Delta$ mod $\underline{\mth{B}}$ are in $1-1$ correspondence with the subsets of $\underline{\Delta}$; obviously, the $\mathbf{F}$-stable ones correspond to the $\mathbf{F}$-stable subsets of $\underline{\Delta}$, which are themselves in $1-1$ correspondence with the subsets of $\Delta$. We usr the corollaty of the proposition \ref{clrel} to conclude. $\Box$

\section{Steinberg representations}

\subsection{Generalities}

In this section, we'll prove the main results of this paper. Let's start with some general facts about Steinberg-like complex representations of finite groups in the most general setting. If $\pi,\pi'$ are two representations of a finite group $G$, we'll write $(\pi,\pi')$ for the dimension of the space of intertwining operators from $\pi$ to $\pi'$.

Let ${\mth{G}}$ be a finite group and let ${\mth{L}}$ be a subgroup of ${\mth{G}}$. Set $1_{\mth{L}}=1_{\mth{L}}^{\mth{G}}=Ind_{\mth{L}}^{\mth{G}}1$; its space $V$ is the space of left ${\mth{L}}$-invariant ${\mth{C}}$-valued functions on ${\mth{G}}$, and $1_{\mth{L}}$ acts on $V$ by $1_{\mth{L}}(g)f=f(.g)$. We have the following result:

\begin{prop}\label{entdc}
Let ${\mth{L}},{\mth{L}}'$ be two subgroups of ${\mth{G}}$; we have:
\[(1_{\mth{L}},1_{{\mth{L}}'})=\#({\mth{L}}'\backslash {\mth{G}}/{\mth{L}}).\]
\end{prop}

This is just the theorem 44.5 of \cite{cr} applied to the particular case of the trivial character. $\Box$

We also have:

\begin{lemme}\label{intind}
For every subgroups ${\mth{L}},{\mth{L}}'$ of ${\mth{G}}$, $1_{\mth{L}}\cap 1_{{\mth{L}}'}=1_{<{\mth{L}},{\mth{L}}'>}$, where $<{\mth{L}},{\mth{L}}'>$ is the subgroup of ${\mth{G}}$ generated by ${\mth{L}}$ and ${\mth{L}}'$.
\end{lemme}
The subspace $1_{\mth{L}}\cap 1_{{\mth{L}}'}$ is the space of elements of $V$ which are left-invariant both by ${\mth{L}}$ and by ${\mth{L}}'$, hence by $<{\mth{L}},{\mth{L}}'>$. $\Box$

In particular, for every subgroup ${\mth{L}}'$ of ${\mth{G}}$ containing ${\mth{L}}$, $1_{{\mth{L}}'}$ is a subrepresentation of $1_{\mth{L}}$. Le $st_{\mth{L}}$ be the following quotient:
\[st_{\mth{L}}=1_{\mth{L}}/(\sum_{{\mth{L}}'\supsetneq {\mth{L}}}1_{{\mth{L}}'}).\]
Since we are dealing here with complex representations of finite groups, all these representations are unitary. The representation $st_{\mth{L}}$ can then also be viewed as a subrepresentation of $1_{\mth{L}}$; more precisely, we have the following proposition:

\begin{prop}
The representation $st_{\mth{L}}$ is isomorphic to the subrepresentation of $1_{\mth{L}}$ formed by the elements $v$ of the space $V$ of $1_{\mth{L}}$ such that for every ${\mth{L}}'\supsetneq {\mth{L}}$ and for every $g\in {\mth{G}}$, we have $\sum_{h\in {\mth{L}}}v(hg)=0$.
\end{prop}

Consider the following hermitian product on $V$:
\[(v,v')=\sum_{g\in {\mth{G}}}v(g)\overline{v'(g)};\]
t is obviously $G$-stable and positive definite. Let $V'$ be the subspace of the elements $v$ satisfying the condition of the proposition.  It is easy to see that $V'$ is the orthogonal in $V$ for the above hermitian product of the sum of the spaces of the representations $1_{{\mth{L}}'}$, where ${\mth{L}}'$ runs through the set of subgroups of ${\mth{G}}$ strictly containing ${\mth{L}}$; Since the hermitian product is positive definite, we have:
\[V=V'\oplus(\sum_{{\mth{L}}'\supsetneq {\mth{L}}}1_{{\mth{L}}'}).\]
The proposition follows immediately. $\Box$

We also have:

\begin{prop}
Let $\mathcal{G}_{\mth{L}}$ be the set of subgroups of ${\mth{G}}$ containing ${\mth{L}}$, and for every $I\subset\mathcal{G}_{\mth{L}}$, let ${\mth{L}}_I$ be the subgroup generated by the elements of $I$. In the Grothendieck group of ${\mth{G}}$. we have:
\[st_{\mth{L}}=\sum_{I\subset\mathcal{G}_{\mth{L}}}(-1)^{\#(I)}1_{{\mth{L}}_I}.\]
\end{prop}

By lemma \ref{intind} and an obvious induction, for every $I\subset\mathcal{G}_{\mth{L}}$, we have $\bigcap_{ð\in I}1_{\mth{L}}=1_{{\mth{L}}_I}$. The lemma folows easily. $\Box$

We will use this proposition to prove the following one:

\begin{prop}
We have $1_{\mth{L}}=\oplus_{{\mth{L}}'\supset {\mth{L}}}st_{{\mth{L}}'}$.
\end{prop}

We will prove this result by induction on $[{\mth{G}}:{\mth{L}}]$, the case ${\mth{L}}={\mth{G}}$ being trivial. We have:
\[1_{\mth{L}}=st_{\mth{L}}\oplus(\sum_{{\mth{L}}'\supsetneq {\mth{L}}}1_{{\mth{L}}'}).\]
Using the induction hypothesis, we obtain:
\[1_{\mth{L}}=st_{\mth{L}}\oplus(\sum_{{\mth{L}}'\supsetneq {\mth{L}}}\bigoplus_{{\mth{L}}''\supset {\mth{L}}'}st_{{\mth{L}}''})\]
\[=\sum_{{\mth{L}}'\supset {\mth{L}}}st_{{\mth{L}}'}.\]
We still have to prove that the sum is direct, which is equivalent to say that the above equality holds in the Grothendieck group of ${\mth{G}}$. In this group, according to the previous proposition, we have, with the same notations as in that proposition:
\[1_{\mth{L}}=st_{\mth{L}}+\sum_{I\subset\mathcal{G}_{\mth{L}},I\neq\emptyset}(-1)^{\#(I)+1}1_{{\mth{L}}_I}.\]
Using the induction hypothesis once again, we obtain:
\[1_{\mth{L}}=st_{\mth{L}}+\sum_{I\subset\mathcal{G},I\neq\emptyset}(-1)^{\#(I)+1}(\sum_{{\mth{L}}'\supset {\mth{L}}_I}st_{{\mth{L}}'})\]
\[=1_{\mth{L}}+\sum_{{\mth{L}}'\subsetneq {\mth{G}}}(\sum_{J\neq\emptyset, {\mth{L}}_I\subset {\mth{L}}'}(-1)^{\#(J)+1})st_{{\mth{L}}'}.\]
Consider now the set $I_{{\mth{L}}'}$ of elements of $\mathcal{G}_{\mth{L}}$ contained in a given ${\mth{L}}'\supset {\mth{L}}$. Then obviously ${\mth{L}}_{I_{{\mth{L}}'}}={\mth{L}}'$ and for every $I$, ${\mth{L}}_I\subset {\mth{L}}'$ if and only if $I\subset I_{{\mth{L}}'}$. We obtain:
\[\sum_{J\neq\emptyset, {\mth{L}}_I\subset {\mth{L}}'}(-1)^{\#(J)+1}=1\]
for every ${\mth{L}}'$, from which we deduce the desired equality. $\Box$

Let $\Delta=\Delta_{\mth{L}}$ be the set of subgroups ${\mth{L}}'$ of ${\mth{G}}$ which strictly contain ${\mth{L}}$ and which are minimal for that property. We obviously have:
\[st_{\mth{L}}=1_{\mth{L}}/(\sum_{{\mth{L}}'\in\Delta}1_{{\mth{L}}'})\]
from which we deduce the following formula in the Grothendieck group:
\[st_{\mth{L}}=\sum_{I\subset\Delta}(-1)^{\#(I)}1_{{\mth{L}}_I}.\]
We deduce from this the following proposition:

\begin{prop}\label{stent}
Let ${\mth{L}},{\mth{L}}'$ be two subgroups of ${\mth{G}}$. We have:
\[(st_{\mth{L}},st_{{\mth{L}}'})=\sum_{I\in\Delta_{\mth{L}},J\in\Delta_{{\mth{L}}'}}(-1)^{\#(I)+\#(J)}\#({\mth{L}}_I\backslash {\mth{G}}/{\mth{L}}'_J);\]
\[(1_{\mth{L}},st_{{\mth{L}}'})=\sum_{I\in\Delta_{\mth{L}}'}(-1)^{\#(I)}\#({\mth{L}}\backslash {\mth{G}}/{\mth{L}}'_I).\]
\end{prop}

This is an immediate consequence of the above formula and the proposition \ref{entdc}. $\Box$

\subsection{The main results}

We will now prove the main result of the paper. Now $\underline{\mth{G}}$ is defined as in the proposition \ref{gfparah}, with $k$ finite, and ${\mth{G}}$ is the group of its $k$-points; let ${\mth{B}}$ be a pseudo-Borel subgroup of ${\mth{G}}$, and let $\Delta$ be the set of negative simple roots associated to ${\mth{B}}$. First we observe the following fact:

\begin{lemme}
The unipotent radical $R_u({\mth{B}})$ of ${\mth{B}}$ is a normal subgroup of ${\mth{B}}_\Delta$.
\end{lemme}

This is an easy consequence of the commutator relations. $\Box$

Any representation of ${\mth{B}}_\Delta$ trivial on $R_u({\mth{B}})$ can thus be viewed as a representation of the quotient group ${\mth{B}}_\Delta/R_u({\mth{B}})$, which is isomorphic to the semi-direct product of ${\mth{T}}$ with the group ${\mth{U}}_\Delta=\prod_{\alpha\in\Delta}{\mth{U}}_{\alpha,h-1}$.

Let $\chi$ be a character of ${\mth{U}}_\Delta$; we'll say $\chi$ is regular if for every $\alpha\in\Delta$, the restriction of $\chi$ to ${\mth{U}}_{\alpha,h-1}$ is nontrivial.

We can now prove the following theorem:

\begin{theo}
Assume the conditions on $p$ and $\underline{\Phi}$ are the same as in the proposition \ref{clbi}. Then $st_{\mth{B}}$ is multiplicity-free, and its irreducible components are in $1-1$ correspondance with the orbits of the action of ${\mth{T}}$ on the regular characters of ${\mth{U}}_\Delta$. In particular, if $\underline{\mth{G}}$ is adjoint modulo its center, then $st_{\mth{B}}$ is irreducible.
\end{theo}

According to the proposition \ref{stent}, we have:
\[(1_{\mth{B}},st_{\mth{B}})=\sum_{I\in\Delta_{\mth{B}}}(-1)^{\#(I)}\#({\mth{B}}\backslash{\mth{G}}/{\mth{B}}_I).\]
Let $\overline{g}$ be an element of ${\mth{B}}\backslash{\mth{G}}/{\mth{B}}_\Delta$ distinct from ${\mth{B}}_\Delta$. According to the proposition \ref{clbirel} applied to $g^{-1}$, there exists $\alpha\in\Delta$ such that for every $g\in\overline{g}$, ${\mth{B}}g{\mth{B}}={\mth{B}} g{\mth{B}}_\alpha$. Let $I$ be a subset of $\Delta-\{i\}$; we then have:
\[{\mth{B}}g{\mth{B}}_I=\bigcup_{h\in{\mth{B}}_I}{\mth{B}}gh{\mth{B}}\]
\[=\bigcup_{h\in{\mth{B}}_I}{\mth{B}} gh{\mth{B}}_\alpha={\mth{B}}g{\mth{B}}_I{\mth{B}}_\alpha\]
which is right-invariant by both ${\mth{B}}_I$ and ${\mth{B}}_\alpha$, hence by ${\mth{B}}_{I\cup\{\alpha\}}$, and we finally obtain:
\[{\mth{B}}g{\mth{B}}_I={\mth{B}}g{\mth{B}}_{I\cup\{\alpha\}}.\]
Since this is true for every $g\in\overline{g}$, we obtain:
\[\#({\mth{B}}\backslash\overline{g}/{\mth{B}}_I)=\#({\mth{B}}\backslash\overline{g}/{\mth{B}}_{I\cup\{\alpha\}})\]
for every $I$, hence:
\[\sum_{I\in\Delta_{\mth{B}}}(-1)^{\#(I)}\#({\mth{B}}\backslash\overline{g}/{\mth{B}}_I)=0.\]
Since this is true for every $\overline{g}\not\in{\mth{B}}_\Delta$, we obtain:
\[(1_{\mth{B}},st_{\mth{B}})=\sum_{I\in\Delta_{\mth{B}}}(-1)^{\#(I)}\#({\mth{B}}\backslash{\mth{B}}_\Delta/{\mth{B}}_I)\]
\[=(1_{\mth{B}}^{{\mth{B}}_\Delta},st_{\mth{B}}^{{\mth{B}}_\Delta}).\]
Tnis equality implies that for every irreducible component $\pi$ of $st_{\mth{B}}^{{\mth{B}}_\Delta}$, $Ind_{{\mth{B}}_\Delta}^{\mth{G}}\pi$ is irreducible, and that the multiplicity of $Ind_{{\mth{B}}_\Delta}^{\mth{G}}\pi$ in $st_{\mth{B}}$ is the same as the multiplicity of $\pi$ in $st_{\mth{B}}^{{\mth{B}}_\Delta}$. We then only have to prove that the assertions of the theorem are true for that last representation.

Since $R_u({\mth{B}})$ is normal in ${\mth{B}}_\Delta$, for every $I\subset\Delta$, $1_{{\mth{B}}_I}^{{\mth{B}}_\Delta}$ is trivial on this group, hence $st_{\mth{B}}^{{\mth{B}}_\Delta}$ is also trivial on it. These representatnins can thus be viewed as representations of the quotient ${\mth{B}}_\Delta/R_u({\mth{B}})$.

Consider first their restriction to ${\mth{U}}_\Delta$. The restriction of $1_{\mth{B}}^{{\mth{B}}_\Delta}$ is simply the regular representation of ${\mth{U}}_\Delta$; since this group is abelian, it is the direct sum of all characters of ${\mth{U}}_\Delta$, each one of them occuring with multiplicity one. For every $I\subset\Delta$, the restriction of $1_{{\mth{B}}_I}^{{\mth{B}}_\Delta}$ is the subrepresentation of that regular representation containing exactly the characters of ${\mth{U}}_\Delta$ which are trivial on $\prod_{\alpha\in I}{\mth{U}}_{\alpha,h-1}$; the restriction of $st_{\mth{B}}^{{\mth{B}}_\Delta}$ is then the direct sum of the characters of ${\mth{U}}_\Delta$ which doesn't satisfy any such condition, that is the regular characters. The representation $st_{\mth{B}}^{{\mth{B}}_\Delta}$ is then the direct sum of the ${\mth{T}}$-orbits of such characters, and any two such subrepresentations are nonisomorphic, hence the first assertion.

Now assume $\underline{\mth{G}}$ is adjoint modulo its center. According to proposition \ref{dcadjrel}, the cardinal of ${\mth{B}}\backslash{\mth{B}}_\Delta/{\mth{B}}$ is $2^{\#(\Delta)}$, and an obvious consequence of the same proposition is that for every $I\subset\Delta$, the cardinal of ${\mth{B}}\backslash{\mth{B}}_\Delta/{\mth{B}}_I$ is $2^{\#(\Delta)-\#(I)}$. We thus obtain from proposition \ref{stent}:
\[(1_{\mth{B}}^{{\mth{B}}_\Delta},st_{\mth{B}}^{\mth{\Delta}})=\sum_{I\subset\Delta}2^{\#(\Delta)-\#(I)}=1,\]
hence $st_{\mth{B}}^{{\mth{B}}_\Delta}$ is irreducible, as required. $\Box$

We deduce from this the following results:

\begin{prop}
Assume $\underline{\mth{G}}$ is adjoint modulo its center. Let ${\mth{Z}}_{\mth{G}}$ be the center of ${\mth{G}}$. The quotient ${\mth{T}}/{\mth{Z}}_{\mth{G}}$ acts transitively and faithfully on the set of regular characters of ${\mth{U}}_\Delta$.
\end{prop}

The transitivity of the action is an immediate consequence of the previous proposition; we will now prove its faithfulness. Assume first $\underline{\mth{G}}$ is split. For each $\alpha\in\Delta$, let $\phi_\alpha$ be a group isomorphism between ${\mth{F}}_q$ and ${\mth{U}}_{\alpha,h-1}$, and let $\eta$ be a primitive $p$-th root of unity in ${\mth{C}}$; a character of ${\mth{U}}_\Delta$ is of the following form:
\[\chi:\prod_{\alpha\in\Delta}\phi(\alpha)(x_\alpha)\mapsto\prod_{\alpha\in\Delta}\eta^{tr(\lambda_\alpha x_\alpha)},\]
with the $\lambda_\alpha$ being elements of ${\mth{F}}_q$; $\chi$ is regular iff all the $\lambda_\alpha$ are nonzero. The group ${\mth{T}}$ acts on these characters by:
\[Ad(t)\chi:\prod_{\alpha\in\Delta}\phi(\alpha)(x_\alpha)\mapsto\prod_{\alpha\in\Delta}\eta^{tr(\alpha(t)\lambda_\alpha x_\alpha)}.\]
Since the action of ${\mth{T}}$ is transitive, we see by cardinality that the stabilisator of any regular character is reduced to ${\mth{Z}}_{\mth{G}}$, hence the action of the quotient is faithful.

Now consider the general case. Let $k'$ be a finite extension of $k$ on which $\underline{\mth{G}}$ splits. Set ${\mth{G}}'=\underline{\mth{G}}(k')$ and ${\mth{T}}'=\underline{\mth{T}}(k')$; let $\Delta'$ be the set of negative simple roots of ${\mth{G}}'$ relatively to ${\mth{T}}'$ and set ${\mth{U}}'_{\Delta'}=\prod_{\alpha'\in\Delta'}{\mth{U}}_{\alpha',h-1}$. The regular characters of ${\mth{U}}_\Delta$ are in $1-1$ correspondance with the $\mathbf{F}$-stable regular characters of ${\mth{U}}'_{\Delta'}$, and we see from the split case that ${\mth{T}}/{\mth{Z}}_{\mth{G}}$ acts faithfully on those characters. $\Box$

\begin{cor}\label{ftf}
Assume now $\underline{\mth{G}}$ is not necessarily adjoint modulo its center. Then ${\mth{T}}/{\mth{Z}}_{\mth{G}}$ acts faithfully (but not necessarily transitively) on the set of regular characters of ${\mth{U}}_\Delta$.
\end{cor}

Let $\underline{\mth{G}}_{ad}$ be an adjoint group defined over $k$ and such that there exists a $k$-isogeny between $\underline{\mth{G}}$ and $\underline{\mth{G}}_{ad}$, let $\underline{\mth{T}}_{ad}$ be a maximal torus of $\underline{\mth{T}}$ containing the image of $\underline{\mth{T}}$ by that isogeny, and let ${\mth{G}}_{ad}$ (resp; ${\mth{T}}_{ad}$) be the group of $k$-points of $\underline{\mth{G}}_{ad}$ (resp. $\underline{\mth{T}}_{ad}$). Since the kernel of the isogeny is the center ${\mth{Z}}_{\underline{\mth{G}}}$ of $\underline{\mth{G}}$, it induces an injection from ${\mth{T}}/{\mth{Z}}_{\mth{G}}$ into ${\mth{T}}_{ad}/{\mth{Z}}_{{\mth{G}}_{ad}}={\mth{T}}_{ad}$; moreover, ${\mth{U}}_\delta$ is isomorphic to the corresponding subgroup of ${\mth{G}}_{ad}$. The faitfulness of the action of ${\mth{T}}/{\mth{Z}}_{\mth{G}}$ on the set of regular characters of this group follows then from the faithfulness of the action of ${\mth{T}}_{ad}$. $\Box$

Let now ${\mth{P}}$ be any generalized pseudo-parabolic subgroup contained in ${\mth{B}}_0$, let $f$ be the corresponding concave function, let $\Delta=\Delta_{\mth{P}}$ be the corresponding subset of $\Phi^-$, and let ${\mth{U}}_{\Delta}$ be the quotient $({\mth{P}}_\Delta\cap{\mth{U}}^-)/({\mth{P}}\cap{\mth{U}}^-)$. Assume this quotient is an abelian group; by the same reasoning as in the previous proposition, the representation $st_{\mth{P}}^{{\mth{P}}_\Delta}$ is the direct sum of the ${\mth{T}}$-orbits of regular characters of ${\mth{U}}_\Delta$, hence $st_{\mth{P}}$ has at least as many irreducible components as there are such ${\mth{T}}$-orbits. Now we'll see that the number of such components can be quite large.

We'll say ${\mth{P}}$ is generic if $\Delta=\Phi^-$ and ${\mth{U}}_\Delta$ is abelian. This can be translated in terms of concave functions the following way:
\begin{itemize}
\item for every $\alpha\in\Phi^-$, $f(\alpha)\geq 2$, and for every $\beta<\alpha$, $f(\beta)>f(\alpha)$;
\item for every $\alpha,\beta\in\Phi^-$ such that $\alpha+\beta\in\Phi^-$, $f(\alpha+\beta)<f(\alpha)+f(\beta)-1$.
\end{itemize}

We'll say such a concave function is generic. The fact that, assuming the first property is true, the second property is equivalent to ${\mth{U}}_\Delta$ being an abelian group is an immediate consequence of the commutator relations.

Generic generalized pseudo-parabolic subgroups exist when $h$ is large enough, more precisely, we have:

\begin{prop}
Assume $\Phi$ is irreducible and not of type $A_1$. The group ${\mth{G}}$ admits generic generalized pseudo-parabolic subgroups if and only if $h\geq h_0+1$, where $h_0$ is the Coxeter number of $\Phi$.
\end{prop}

Let $f$ be the concave function on $\Phi^-$ defined by $f(\alpha)=l(\alpha)+2$, where $l(\alpha)$ is the length of $\alpha$ as defined in the proof of the proposition \ref{clbi}. We can easily check that $f$ is generic, and for every $\alpha\in\Phi^-$, we have $f(\alpha)\leq h_0+1$, hence if $h\geq h_0+1$, ${\mth{P}}_f$ is generic.

Now asssume $f$ is a generic concave function on $\Phi^-$ such that $f(\alpha)\leq h_0$ for every $\alpha\in\Phi^-$; we'll prove that such a function cannot exist. From the first property and an easy induction, we must have $f(\alpha)\leq 2$ for every simple root $\alpha$. Let $\alpha,\beta$ be two simple roots such that $\alpha+\beta\in\Phi$; we then have:
\[f(\alpha)<f(\alpha+\beta)<f(\alpha)+f(\beta)-1,\]
which is impossible since $f(\beta)\leq 2$. Hence if $h\leq h_0$, ${\mth{P}}_f$ donesn't admit any generic generalized pseudo-parabolic subgroups. $\Box$

When $\Phi$ is of type $A_1$, all generalized pseudo-parabolic subgroups which are neither ${\mth{G}}$ itself nor a Borel subgroup are generic. If $\Phi$ is reducible, we can check the genericity of generalized pseudo-parabolic subgroups componentwise.

\begin{prop}\label{generic}
Assume ${\mth{P}}$ is generic. Then the number of irreducible components of $st_{\mth{P}}$ is at least $(q-1)^{\#(\underline{\Phi}^-)-rg(\underline{\Phi})}$.
\end{prop}

Let $\Delta_0$ be the set of negative simple roots of $\Phi$. By the corollary \ref{ftf}, ${\mth{T}}$ acts faithfully on the regular characters of ${\mth{U}}_{\Delta_0}=\prod_{\alpha\in\Delta_0}({\mth{U}}_{\alpha,f(\alpha)-1}/{\mth{U}}_{\alpha,f(\alpha)})$; hence in a given ${\mth{T}}$-orbit of regular characters of ${\mth{U}}_\Delta$, all characters must have different restrictions to ${\mth{U}}_{\Delta_0}$. We deduce from this that the number of such ${\mth{T}}$-orbits is at least equal to the number of regular characters of the quotient ${\mth{U}}_\Delta/{\mth{U}}_{\Delta_0}$, which is at least $(q-1)^{\#(\underline{\Phi}^-)-rg(\underline{\Phi})}$. $\Box$

\section{Proof of the proposition \ref{clbi}}

We will prove the result for the double classes mod $R_u(\underline{\mth{B}})$ on the left and $\underline{\mth{B}}$ on the right; the result for double classes mod $\underline{\mth{B}}$ will then follow immediately.

First we will assume that $\underline{\Phi}$ is irreducible. Let $g$ be any element of $\underline{\mth{G}}$ not belonging to $\underline{\mth{B}}_{\underline{\Delta}}$. By eventually replacing $g$ by another element of its double class, according to the discussion at the beginning of section \ref{dc}, we can assume $g$ is of the form $nu$, where $u$ is an element of $\underline{\mth{U}}_1^-$ and $n$ is a representative of an element $w$ of $W$.

Assume first $w=1$; we can then also assume $n=1$, hence $g=u$.
For every $\alpha<0$, we can define its length as the number of negative simple roots it is the sum of. Write $u=\prod_{\alpha<0}u_\alpha$, the product being taken in any arbitrarily chosen order. Consider the couple of integers $(v,l)$ defined the following way:
\begin{itemize}
\item $v$ is the minimal valuation of the $u_\alpha$, $\alpha\in\underline{\Phi}^-$;
\item $l$ is the maximal length of the $\alpha$ such that the valuation of $u_\alpha$ is $v$.
\end{itemize}
As an easy consequence of the commutator relations, $v$ and $l$ don't depend on the choice of the order on the roots. Moreover, we can easily see that if $u'\in\underline{\mth{U}}^-$ belongs to the same double class modulo $\underline{\mth{B}}_{\underline{\Delta}}$ as $u$, then the couple $(v,l)$ associated to $u'$ is the same.

First assume $l>1$. Let $\Gamma_l$ (resp. $\Gamma_{l-1}$) be the subset of $\alpha\in\underline{\Phi}^-$ such that $l(\alpha)=l$ (resp. $l-1$). As an easy consequence of the commutator relations, if $\beta$ is an element of $\Gamma_{l-1}$ and $u'$ an element of $\underline{\mth{U}}_{-\beta,h-v-1}$, $[u,u']$ belongs to $\underline{\mth{B}}_{\underline{\Delta}}$. Writing $[u,u']=u_0b$, with $b\in R_u(\underline{\mth{B}})$ and $u_0\in{\mth{U}}_{h-1}^-=\prod_{\alpha<0}\underline{\mth{U}}_{\alpha,h-1}$, we obtain that $u_0^{-1}u=b[u,u']^{-1}u=bu'uu'^{-1}\in R_u(\underline{\mth{B}})u\underline{\mth{B}}$. We will prove the following claim: when $u'$ runs through $\prod_{\beta\in\Gamma_{l-1}}\underline{\mth{U}}_{-\beta,h-v-1}$, $u_0$ runs through a subgroup $\underline{\mth{U}}_u$ of $\underline{\mth{U}}_{h-1}^-$ which contains a subgroup of the form $\prod_{\delta\in I}\underline{\mth{U}}_{\delta,h-1}$, with $I$ being a nonempty subset of $\underline{\Delta}$; moreover, the subset $I$ depends only on the triplet $(l,v,\Gamma_l)$, which depends only on the class of $u$ modulo $\prod_{\alpha,l(\alpha)<l}\underline{\mth{U}}_{\alpha,v}\prod_{\alpha,l(\alpha)\geq l}\underline{\mth{U}}_{\alpha,v+1}$ (for any $x$ in this group and any $y\in\underline{\mth{U}}_{-\beta,h-v-1}$ $[x,y]\in R_u(\underline{\mth{B}})$ by the commutator relarions), which itself obviously depends only on the double class modulo $\underline{\mth{B}}_{\underline{\Delta}}$ containing $u$, this will be enough to prove the assertion of the proposition in this case.

First we will prove some lemmas.

\begin{lemme}\label{glinf}
For every $l$, the cardinal of $\Gamma_l$ is smaller than or equal to the cardinal of $\Gamma_{l-1}$, and strictly smaller if $l=2$.
\end{lemme}
We will use the notations of \cite[4]{bou}.
\begin{itemize}
\item Assume $\underline{\Phi}$ is of type $A_n$. Then for every $i,j$, $j<i$, the root $\varepsilon_i-\varepsilon_j$ is of length $i-j$; hence for every $l$, there are exactly $n-l-1$ negative roots of length $l$.
\item Assume $\underline{\Phi}$ is of type $B_n$. Then for every $i,j$, $j<i$, the root $\varepsilon_i-\varepsilon_j$ is of length $i-j$ and the root $-\varepsilon_i-\varepsilon_j$ is of length $2n-i-j+2$, and for every $i$, the root $-\varepsilon_i$ is of length $n-i+1$. Combining those results, we obtain that for every $l$, there are exactly $[\frac{2n-l+1}2]$ negative roots of length $l$.
\item Assume $\underline{\Phi}$ is of type $C_n$. Then for every $i,j$, $j<i$, the root $\varepsilon_i-\varepsilon_j$ is of length $i-j$ and the root $-\varepsilon_i-\varepsilon_j$ is of length $2n-i-j+1$, and for every $i$, the root $-2\varepsilon_i$ is of length $2(n-i)+1$. Combining those results, we obtain that for every $l$, there are exactly $[\frac{2n-l+1}2]$ negative roots of length $l$, exactly as in the $B_n$ case.
\item Assume $\underline{\Phi}$ is of type $D_n$. Then for every $i,j$, $j<i$, the root $\varepsilon_i-\varepsilon_j$ is of length $i-j$ and the root $-\varepsilon_i-\varepsilon_j$ is of length $2n-i-j$. Combining those results, we obtain that for every $l$, there are exactly $[\frac{2n-l+1}2]$ negative roots of length $l$ if $l\leq n-1$ and $[\frac{2n-l-1}2]$ if $l\geq n$.
\end{itemize}
We deduce from this that in all classical cases, the assertion of the lemma is true. The exceptional cases can easily be checked directly by counting the number of roots of each length.

For the case $l=2$, we can also remark that there is a $1-1$ correspondence between the roots of length $2$ (resp $1$) and the edges (resp. vertices) of the Dynkin diagram of $\underline{\Phi}$. Since that diagram has no cycles, it has strictly less edges than vertices. $\Box$

In particular, when $\underline{\Phi}$ is of type $G_2$, there is exactly one negative root of length $l$ for every $l\in\{2,\dots,5\}$. We deduce then immediately from the commutator relations that if $l=2$, we have $\underline{\mth{U}}_u=\underline{\mth{U}}_{\underline{\Delta}}$, and if $l>2$, if $\alpha$ (resp. $\beta$) is the unique root of length $l$ (resp $l-1$), we deduce from \cite[1, prop. 19]{bou} that $\alpha-\beta$ must belong to $\underline{\Delta}$, and we have $\underline{\mth{U}}_u=\underline{\mth{U}}_{{\alpha-\beta}}$, which proves the claim in this case. In the sequel we will assume that $\underline{\Phi}$ is not of type $G_2$.

Consider the graph $\mathcal{G}_l$ defined the following way: its vertices are the elements of $\Gamma_l$ and there is an edge between two distinct vertices $\alpha$ and $\alpha'$ if and only if there exists $\beta\in\Gamma_{l-1}$ such that both $\alpha-\beta$ and $\alpha'-\beta$ are elements of $\underline{\Delta}$.

\begin{lemme}
If such a $\beta$ exists, then it is unique.
\end{lemme}
Assume there exist $\beta,\beta'$ such that $\alpha-\beta,\alpha-\beta',\alpha'-\beta,\alpha'-\beta'$ all belong to $\underline{\Delta}$. Since we have $(\alpha-\beta)+(\alpha'-\beta')=(\alpha-\beta')+(\alpha'-\beta)$, by linear independence of the elements of $\underline{\Delta}$ we must have either $\alpha-\beta=\alpha-\beta'$ or $\alpha-\beta=\alpha'-\beta$. Since $\alpha\neq\alpha'$, we obtain $\beta=\beta'$ and the result follows. $\Box$

In other words, to every edge of $\mathcal{G}_l$ we can attach a unique element $\beta\in\Gamma_{l-1}$. In general this correspondance is not $1-1$: for a given element $\beta$ of $\Gamma_{l-1}$, if there is only one (resp. at least three) elements $\alpha$ of $\Gamma-l$ such that $\alpha-\beta$ is a root, then $\beta$ is not attached to any edge (resp. is attached to several different edges) of $\mathcal{G}_l$.

\begin{lemme}
The graph $\mathcal{G}_l$ is connected.
\end{lemme}

We will prove the lemma by induction on $l$. If $l=2$, there is a $1-1$ correspondance between $\Gamma_l$ and the set of edges of the Dynkin diagram of $\underline{\Phi}$, and there is an edge between two elements of $\Gamma_l$ if and only if the corresponding edges of the Dynkin diagram share a common vertex; since the Dynkin diagram is connected, $\mathcal{G}_l$ is also connected.

Assume now $l\geq 3$ and let $\alpha,\alpha'$ be distinct elements of $\Gamma_l$. By (\cite[1. prop.19]{bou} there exists at least one element $\beta$ (resp. $\beta'$) of $\Gamma_{l-1}$ such that $\alpha-\beta$ (resp. $\alpha'-\beta'$) is an element of $\underline{\Delta}$, and since by induction hypothesis the graph $\mathcal{G}_{l-1}$ is connected, there exists $\beta_0=\beta,\beta_1,\dots,\beta_s=\beta'$ such that for every $i\in\{1,\dots,r\}$, $\beta_i\neq\beta_{i-1}$ and there exists $\gamma_i\in\Gamma_{l-2}$ such that $\beta_{i-1}-\gamma_i=\delta_i$ and $\beta_i-\gamma_i=\delta'_i$ are elements of $\underline{\Delta}$. But then $\alpha_i=\gamma_i+\delta_i+\delta'_i$ is an element of $\underline{\Phi}$, and more precisely an element of $\Gamma_l$, and, setting $\alpha_0=\alpha$ and $\alpha_{r+1}=\alpha'$, for every $i\in\{1,\dots,r+1\}$; $\alpha_{i-1}-\beta_{i-1}$ and $\alpha_i-\beta_{i-1}$ are both elements of $\underline{\Delta}$; hence there is a path between every pair of vertices of $\mathcal{G}$, which proves the lemma. $\Box$

Now we will prove the claim. Assume first there are no cycles in $\mathcal{G}_l$. Then its number of edges is exactly the cardinal of $\Gamma_l$ minus one, which is strictly smaller than the cardinal of $\Gamma_{l-1}$, hence there exists at least one $\beta\in\Gamma_{l-1}$ which isn't attached to any edge, which is equivalent to say that there is exactly one $\alpha\in\Gamma_l$ such that $\alpha-\beta$ is a root. If $v(u_\alpha)=v$, then the commutator relations imply that $\underline{\mth{U}}_u$ contains $\underline{\mth{U}}_{{\alpha-\beta}}$. Assume now $v(u_\alpha)>v$; the connectedness of $\mathcal{G}_l$ implies then that there exist $\alpha_0=\alpha,\dots,\alpha_r$ such that:
\begin{itemize}
\item for every $i\in\{1,\dots,r\}$, there is an edge between $\alpha_{i-1}$ and $\alpha_i$;
\item for every $i<r$, $v(\alpha_i)>v$, and $v(\alpha_r)=v$.
\end{itemize}
Let $\beta_r$ be the element of $\Gamma_{l-1}$ corresponding to the edge between $\alpha_{r-1}$ and $\alpha_r$. Since $\mathcal{G}_l$ has no cycles, there can be no other $\alpha$ than these two such that $\alpha-\beta_r$ is a root. We conclude the same way, replacing $\beta$ and $\alpha$ by $\beta_r$ and $\alpha_r$.

We will now check that when $\underline{\Phi}$ is of type $A_n$, $B_n$ or $C_n$, the graphs $\mathcal{G}_l$ don't have any cycles, which means that the claim is now proved in these cases. Assume first $\underline{\Phi}$ is of type $A_n$; we then have $\Gamma_l=\{-\varepsilon_i+\varepsilon_{i+l}|1\leq i\leq n-l\}$ and the graph $\mathcal{G}_l$ is the following one:

\begin{picture}(300,50)(0,0)
\put(25,20){\circle{6}}
\put(5,7){$\varepsilon_{l+1}-\varepsilon_1$}
\put(28,20){\line(1,0){64}}
\put(45,26){\tiny $\varepsilon_{l+1}-\varepsilon_2$}
\put(95,20){\circle{6}}
\put(75,7){$\varepsilon_{l+2}-\varepsilon_2$}
\put(98,20){\dashbox{2}(64,0){}}
\put(165,20){\circle{6}}
\put(135,7){$\varepsilon_{n-1}-\varepsilon_{n-l-1}$}
\put(168,20){\line(1,0){64}}
\put(175,26){\tiny $\varepsilon_{n-1}-\varepsilon_{n-l}$}
\put(235,20){\circle{6}}
\put(215,7){$\varepsilon_n-\varepsilon_{n-l}$}
\end{picture}

(Edges are labeled with the names of the corresponding elements of $\Gamma_{l-1}$.)

There are no cycles in this graph.

Assume now $\underline{\Phi}$ is of type $B_n$. Then $\Gamma_l=\{\varepsilon_{l+i}-\varepsilon_i|1\leq i\leq n-l\}\cup\{-\varepsilon_{n-l+1}\}\cup\-\varepsilon_{n-l+1+i}-\varepsilon_{n+1-i}|1\leq i\leq[\frac l2]\}$ (with the leftmost two subsets being eventually empty). The graph $\mathcal{G}_l$ is the following graph:

\begin{picture}(400,50)(0,0)
\put(25,20){\circle{6}}
\put(5,7){$\varepsilon_{l+1}-\varepsilon_1$}
\put(28,20){\line(1,0){64}}
\put(45,26){\tiny $\varepsilon_{l+1}-\varepsilon_2$}
\put(95,20){\circle{6}}
\put(75,7){$\varepsilon_{l+2}-\varepsilon_2$}
\put(98,20){\dashbox{2}(64,0){}}
\put(165,20){\circle{6}}
\put(135,7){$\varepsilon_n-\varepsilon_{n-l}$}
\put(168,20){\line(1,0){64}}
\put(185,26){\tiny $-\varepsilon_{n-l+2}$}
\put(235,20){\circle{6}}
\put(225,7){$-\varepsilon_{n-l+1}$}
\put(238,20){\line(1,0){64}}
\put(255,26){\tiny $-\varepsilon_{n-l+3}-\varepsilon_n$}
\put(305,20){\circle{6}}
\put(285,7){$-\varepsilon_{n-l+2}-\varepsilon_n$}
\put(308,20){\dashbox{2}(64,0){}}
\put(375,20){\circle{6}}
\put(345,7){$-\varepsilon_{n-[\frac{l+1}2]}-\varepsilon_{n-[\frac l2]+1}$}
\end{picture}

if $l<n$, and if $l\geq n$, the rightmost part of the above graph. This graph has no cycles.

Assume now $\underline{\Phi}$ is of type $C_n$. Then $\Gamma_l=\{\varepsilon_{l+i}-\varepsilon_i|1\leq i\leq n-l\}\cup\{\varepsilon_{n-l+i}-\varepsilon_{n+1-i}|1\leq i\leq\frac l2\}$ if $l$ is even, and $\Gamma_l=\{\varepsilon_{l+i}-\varepsilon_i|1\leq i\leq n-l\}\cup\-\varepsilon_{n-l+i}-\varepsilon_{n+1-i}|1\leq i\leq\frac {l-1}2\}\cup\{-2\varepsilon_{n-\frac{l-1}2}\}$ if $l$ is odd (in both cases, with the leftmost subset being eventually empty). If $l<n$, $\mathcal{G}_l$ is the following graph:

\begin{picture}(400,50)(0,0)
\put(25,20){\circle{6}}
\put(5,7){$\varepsilon_{l+1}-\varepsilon_1$}
\put(28,20){\line(1,0){64}}
\put(45,26){\tiny $\varepsilon_{l+1}-\varepsilon_2$}
\put(95,20){\circle{6}}
\put(75,7){$\varepsilon_{l+2}-\varepsilon_2$}
\put(98,20){\dashbox{2}(64,0){}}
\put(165,20){\circle{6}}
\put(135,7){$\varepsilon_n-\varepsilon_{n-l}$}
\put(168,20){\line(1,0){64}}
\put(185,26){\tiny $\varepsilon_n-\varepsilon_{n-l+1}$}
\put(235,20){\circle{6}}
\put(210,7){$-\varepsilon_{n-l+1}-\varepsilon_n$}
\put(238,20){\dashbox{2}(64,0){}}
\put(305,20){\circle{6}}
\put(280,7){$-\varepsilon_{n-\frac l2}-\varepsilon_{n-\frac l2+1}$}
\end{picture}

if $l$ is even, and:

\begin{picture}(400,50)(0,0)
\put(25,20){\circle{6}}
\put(5,7){$\varepsilon_{l+1}-\varepsilon_1$}
\put(28,20){\line(1,0){64}}
\put(45,26){\tiny $\varepsilon_{l+1}-\varepsilon_2$}
\put(95,20){\circle{6}}
\put(75,7){$\varepsilon_{l+2}-\varepsilon_2$}
\put(98,20){\dashbox{2}(64,0){}}
\put(165,20){\circle{6}}
\put(135,7){$\varepsilon_n-\varepsilon_{n-l}$}
\put(168,20){\line(1,0){64}}
\put(185,26){\tiny $\varepsilon_n-\varepsilon_{n-l+1}$}
\put(235,20){\circle{6}}
\put(210,7){$-\varepsilon_{n-l+1}-\varepsilon_n$}
\put(238,20){\dashbox{2}(64,0){}}
\put(305,20){\circle{6}}
\put(270,7){$-\varepsilon_{n-\frac{l+1}2}-\varepsilon_{n-\frac{l-3}2}$}
\put(308,20){\line(1,0){64}}
\put(310,29){\tiny $-\varepsilon_{n-\frac{l+1}2}-\varepsilon_{n-\frac{l-1}2}$}
\put(375,20){\circle{6}}
\put(365,7){$-2\varepsilon_{n-\frac{l-1}2}$}
\end{picture}

if $l$ is odd; if $l\geq n$, $\mathcal{G}_l$ is the rightmost part of a similar graph. This graph obviously has no cycles.

Now we'll assume the graph $\mathcal{G}_l$ contains at least one cycle. We will first consider separately the case when $\underline{\Phi}$ is of type $F_4$; in this case, the only $l$ such that $\mathcal{G}_l$ contains at least one cycle are $l=3$ and $l=4$, and in both cases the graph has exactly $3$ vertices and $3$ edges forming a $3$-cycle.

Assume first $l=3$. Then $\Gamma_3=\{\delta_2+\delta_3+\delta_4,\delta_1+\delta_2+\delta_3,\delta_2+2\delta_3\}$ and all three edges are associated to the same $\beta=\delta_2+\delta_3\in\Gamma_2$. Moreover, there are two elements of $\Gamma_2$ which aren't associated to any edge: $\delta_1+\delta_2$ and $\delta_3+\delta_4$.

(From now on all such elements of $\Gamma_{l-1}$ will also be displayed on the graphs, with their names in parentheses.)

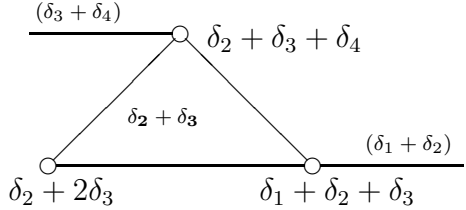
\begin{picture}(400,100)(0,0)
\put(45,20){\circle{6}}
\put(30,7){$\delta_2+2\delta_3$}
\put(48,20){\line(1,0){94}}
\put(75,36){\tiny $\mathbf{\delta_2+\delta_3}$}
\put(145,20){\circle{6}}
\put(125,7){$\delta_1+\delta_2+\delta_3$}
\put(47,22){\line(1,1){46}}
\put(143,22){\line(-1,1){46}}
\put(95,70){\circle{6}}
\put(105,65){$\delta_2+\delta_3+\delta_4$}
\put(148,20){\line(1,0){54}}
\put(165,26){\tiny $(\delta_1+\delta_2)$}
\put(38,70){\line(1,0){54}}
\put(40,76){\tiny $(\delta_3+\delta_4)$}
\end{picture}

We deduce from this that if $v(u_{\delta_1+\delta_2+\delta_3})=v$ (resp. $v(u_{\delta_2+\delta_3+\delta_4})=v$), then $\underline{\mth{U}}_u$ contains $U_{{\delta_3}}$ (resp. $U_{{\delta_2}}$), and if both those valuations are greater than $v$, then $v(u_{\delta_2+2\delta_3})=v$ and $\underline{\mth{U}}_u$ contains $U_{{\delta_3}}$.

Assume now $l=4$. Then $\Gamma_4=\{\delta_1+\delta_2+\delta_3+\delta_4,\delta_1+\delta_2+2\delta_3,\delta_2+2\delta_3+\delta_4\}$ and there is a $1-1$ correspondance between the edges and the elements of $\Gamma_3$.

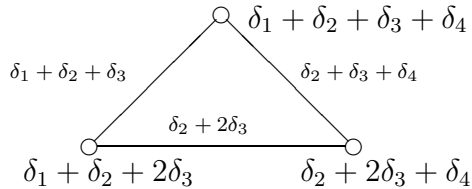
\begin{picture}(400,100)(0,0)
\put(45,20){\circle{6}}
\put(20,7){$\delta_1+\delta_2+2\delta_3$}
\put(48,20){\line(1,0){94}}
\put(75,26){\tiny $\delta_2+2\delta_3$}
\put(145,20){\circle{6}}
\put(125,7){$\delta_2+2\delta_3+\delta_4$}
\put(47,22){\line(1,1){46}}
\put(15,46){\tiny $\delta_1+\delta_2+\delta_3$}
\put(143,22){\line(-1,1){46}}
\put(125,46){\tiny $\delta_2+\delta_3+\delta_4$}
\put(95,70){\circle{6}}
\put(105,65){$\delta_1+\delta_2+\delta_3+\delta_4$}
\end{picture}

Let $\alpha_1,\alpha_2,\alpha_3$ (resp. $\beta_1,\beta_2,\beta_3$) be the elements of $\Gamma_4$ (resp. $\Gamma_3$), in the order in which they are written in the description of $\Gamma_4$ (resp. $\Gamma_3$) above; we will also set $\gamma=\delta_2+\delta_3$, $\delta'_1=\delta_1$, $\delta'_2=\delta_3$, $\delta'_3=\delta_4$.

Let's consider more closely the commutator relations: according to \cite[III]{chev}, there exists a set of isomorphisms $x\mapsto u_\alpha(x)$ from $F_{nr}$ to the $U_\alpha$, $\alpha\in\underline{\Phi}$, satisfying the following condition: for every $\alpha,\beta$ such that $\alpha+\beta$ is a root and every $x,y\in F_{nr}$, we have $[u_\alpha(x),u_\beta(y)]=u_{\alpha+\beta}(\varepsilon_{\alpha,\beta}p_{\alpha,\beta}xy)u'$, where $\varepsilon_{\alpha,\beta}=\pm 1$, $p_{\alpha,\beta}$ is $1$ or $2$ and $u'$ is an element of either $U_{2\alpha+\beta}$ or $U_{\alpha+2\beta}$ (and is $1$ if both are trivial). Moreover, the $\varepsilon_{\alpha,\beta}$ satisfy the following conditions:
\begin{itemize}
\item $\varepsilon_{\beta,\alpha}=-\varepsilon_{\alpha,\beta}$;
\item $\varepsilon_{-\alpha,-\beta}=-\varepsilon_{\alpha,\beta}$;
\item $\varepsilon_{\alpha,\beta}=\varepsilon_{\beta,-\alpha-\beta}=\varepsilon_{-\alpha-\beta,\beta}$;
\item $\varepsilon_{\alpha,\beta}\varepsilon_{\alpha+\beta,\gamma}=\varepsilon_{\beta,\gamma}\varepsilon_{\alpha,\beta+\gamma}$ for every $\alpha,\beta,\gamma$ such that the equality makes sense.
\end{itemize}
Write $u=\prod_\alpha u_\alpha(x_\alpha)$, and for every $i=1,2,3$, write $x_i=x_{\alpha_i}$. By conjugating $u$ by the element $\prod_{i=1}^3u_{-\beta_i}(y_i)$, for some $y_i$ whose valuation is $h-v_u-1$, we see that the element $u\prod_{i=1}^3u_{\delta_i}(\varepsilon_{\alpha_i,-\beta_i}p_{\alpha_i,-\beta_i}x_iy_i+\varepsilon_{\alpha_{i-1},-\beta_{i+1}}p_{\alpha_{i-1},-\beta_{i+1}}x_{i-1}y_{i+1})$ belongs to our double class; to prove the claim, we must prove that the application:
\[(y_1,y_2,y_3)\longmapsto(\varepsilon_{\alpha_i,-\beta_i}p_{\alpha_i,-\beta_i}x_iy_i+\varepsilon_{\alpha_{i-1},-\beta_{i+1}}p_{\alpha_{i-1},-\beta_{i+1}}x_{i-1}y_{i+1})_{i=1,2,3}\]
induces an isomorphism of $\overline{k}$-vector spaces between the space $({\mathfrak{p}}_{nr}^{v_u}/{\mathfrak{p}}_{nr}^{v_u+1})^3\simeq \overline{k}^3$ and the space $({\mathfrak{p}}_{nr}^{h-1}/{\mathfrak{p}}_{nr}^h)^3\simeq \overline{k}^3$, hence that the matrix:
\[\left(\begin{array}{ccc}\varepsilon_{\alpha_1,-\beta_1}p_{\alpha_1,-\beta_1}x_1&\varepsilon_{\alpha_3,-\beta_2}p_{\alpha_3,-\beta_2}x_3\\&\varepsilon_{\alpha_2,-\beta_2}p_{\alpha_2,-\beta_2}x_2&\varepsilon_{\alpha_1,-\beta_3}p_{\alpha_1,-\beta_3}x_3\\\varepsilon_{\alpha_2,-\beta_1}p_{\alpha_2,-\beta_1}x_2&&\varepsilon_{\alpha_3,-\beta_3}p_{\alpha_3,-\beta_3}x_3\end{array}\right)\]
is invertible. Its determinant is:
\[(\prod_{i=1}^3\varepsilon_{\alpha_i,-\beta_i}p_{\alpha_i,-\beta_i}+\prod_{i=1}^3\varepsilon_{\alpha_i,-\beta_{i-1}}p_{\alpha_i,\beta_{i-1}})x_1x_2x_3.\]
All $p_{\alpha,\beta}$ are equal to $1$ except $p_{\alpha_1,-\beta_1}$ which is $2$. Now we will check the equality $\prod_{i=1}^3\varepsilon_{\alpha_i,-\beta_i}=\prod_{i=1}^3\varepsilon_{\alpha_i,-\beta_{i-1}}$: using the properties of the $\varepsilon$ stated above, we obtain for every $i$:
\[\varepsilon_{\alpha_i,-\beta_i}=\varepsilon_{-\beta_i,-\delta_i}\]
\[=\varepsilon_{-\delta_{i+1},-\gamma}\varepsilon_{-\gamma,-\delta_i}\varepsilon_{-\delta_{i-1},-\beta_{i+1}}\]
\[=\varepsilon_{-\gamma,-\delta_{i+1}}\varepsilon_{-\gamma,-\delta_i}\varepsilon_{\alpha_i,-\beta_{i-1}}.\]
By taking the product over $i=1,2,3$, we obtain the desired equality. The determinant of the matrix is then $\pm 3$, hence invertible since we have assumed $p\neq 3$ in the $F_4$ case.

Now the only cases which are not completely solved yet are $D_n$ and $E_n$; we first remark that both are simply-laced. This implies the following properties (see \cite[1]{bou}):
\begin{itemize}
\item there exists a $W$-invariant scalar product $(.,.)$ on $X^*(\underline{\mth{T}})$ such that $(\alpha,\alpha)=2$ for every $\alpha\in\underline{\Phi}$;
\item for every $\alpha\neq\pm\beta$, $(\alpha,\beta)\in\{-1,0,1\}$ and $(\alpha,\beta)=1$ (resp. $-1$) if and only if $\alpha-\beta$ (resp. $\alpha+\beta$) is a root;
\item for every roots $\alpha,\beta$ such that $\alpha+\beta$ is also a root, we have $p_{\alpha,\beta}=1$, and neither $\alpha+2\beta$ nor $2\alpha+\beta$ is a root.
\end{itemize}

We will now examine the cases where $\underline{\Phi}$ is of type $D_n$ or $E_n$ and there are cycles in $\mathcal{G}_l$: the discussion about $F_4$ already gives us some hints about how to proceed. Let's first prove some general facts about those cycles.

Let $\alpha_1,\dots,\alpha_n$, $n\geq 3$, be distinct elements of $\Gamma_l$ forming a cycle in $\mathcal{G}_l$; from now on, we will denote by $\beta_1,\dots,\beta_n$ the elements of $\Gamma_{l-1}$ such that for every $i$, $\beta_i$ corresponds to the edge between $\alpha_{i-1}$ and $\alpha_i$ (the indices being taken modulo $n$). We'll write, for every $i$, $\delta_i=\alpha_i-\beta_i$ and $\delta'_i=\alpha_i-\beta_{i+1}$.

We'll first prove the following lemma:

\begin{lemme}
Let $\alpha$ be any element of $\underline{\Phi}$. There are at most three elements $\delta\in\underline{\Delta}$ such that $\alpha+\delta\in\underline{\Phi}$ (resp. $\alpha-\delta\in\underline{\Phi}$).
\end{lemme}

It is enough to prove it with $\alpha+\delta$, the proof with $\alpha-\delta$ being symmetrical. Let $\delta_1,\dots,\delta_r$ be distinct elements of $\underline{\Delta}$ satisfying that condition. For every $i$, we have $(\alpha,\delta_i)<0$; on the other hand, for $i\neq j$, $(\delta_i,\delta_j)\leq 0$, hence $(\alpha+\delta_i,\delta_j)<0$, from which we deduce that $\alpha+\delta_i+\delta_j\in\underline{\Phi}$. By an obvious induction we obtain that $\alpha+\delta_1+\dots+\delta_r$ is a root. We thus have, since $\alpha$ and $\alpha+\delta_1+\dots+\delta_r$ are obviously not opposite to each other:
\[(\alpha,\alpha+\delta_1+\dots+\delta_r)\geq -1.\]
But $(\alpha,\alpha)=2$ and $(\alpha,\delta_i)<0$ for every $i$; we must then have $r\leq 3$, which proves the lemma. $\Box$

Assume now we have $\beta_i=\beta_j$ for some $i\neq j$. For every $\delta\in\{\delta_i,\delta_j,\delta'_{i-1},\delta'_{j-1}\}$, we tnen have $\beta_i+\delta\in\underline{\Phi}$. Since, according to the previous lemma, the nomber of elements of $\underline{\Delta}$ satisfying that condition is at most three, it implies that either $\delta_i=\delta'_{j-1}$ or $\delta_j=\delta'_{i-1}$, hence $j=i\pm 1$. Assume for example $j=i+1$. Then $\alpha_1,\dots,\hat{\alpha_j},\dots,\alpha_n$ also form a cycle. If at least two of the $\beta_i$ are distinct, by an obvious induction, we see that the cycle contains a subcycle such that all the $\beta_i$ attached to that subcycle are distinct. We'll say a cycle is reduced if it satisfies that last property, nonreduced if it doesn't.

Assume the cycle is reduced and $l\geq 3$; we will then prove that for every $i$, $\gamma_i=\beta_{i-1}-\delta'_{i-1}$ is a root. We have for every $i$:
\[(\alpha_i,\delta'_i)=(\beta_i,\delta'_i)+(\delta_i,\delta'_i).\]
Since the cycle is reduced, $\beta_i\neq\beta_{i+1}$, hence $\delta_i\neq\delta'_i$, from which we deduce $(\delta_i,\delta'_i)\leq 0$, hence $(\beta_i,\delta'_i)\geq(\alpha_i,\delta'_i)=1$, which proves the assertion.

For every $i$, we then have:
\[\beta_i=\gamma_i+\delta'_i=\gamma_{i-1}+\delta_{i-1}.\]
We thus see that the graph $\mathcal{G}_{l-1}$ contains a cycle of length $n$ as well; but that cycle is not necessarily reduced. By iterating the process (reducing the cycle, then decreasing the length of the roots), after a finite number of steps we reach one of the following two situations:
\begin{itemize}
\item $l=2$;
\item all the $\beta_i$ are equal to each other.
\end{itemize}
We will prove that we in fact always reach the second situation. If it is not the case, then the cycle we get for $l=2$ contains a nontrivial reduced subcycle. Consider now the case of a reduced cycle $(\alpha_1,\dots,\alpha_r)$ for $l=2$. With the $\beta_i,\delta_i,\delta'_i$ defined as above, we have for every $i$, $\beta_i=\delta'_i=\delta_{i-1}$, hence $\alpha_i=\delta_{i-1}+\delta_i$. That amounts to say that there is a cycle in the Dynkin diagram of $\underline{\Phi}$, which is known to be impossible, hence a contradiction.

Let $l'$ be the length of the $\beta_i$ for which the second situation occurs; we will call the integer $l-l'$ the level of the cycle.

Now we can examine the different cycle cases; we can always assume that there is at least one $\alpha$ in the cycle such that $v(u_\alpha)=v$, since if it is not true, we are obviously reduced to a simpler case. Assume first that $\mathcal{G}_l$ contains only 3-cycles. First we'll examine the case when one of them is reduced: we will first prove that it is necessarily of level $2$. Assume the contrary; since all the cycles occuring in the iteration process before the last step will be reduced $3$-cycles, it is enough to find a contradiction when our cycle is of level $3$. Let the $\beta_i$ and the $\gamma_i$ be defined as above, and let $\zeta$ be the element of $\Gamma_{l-3}$ such that $\gamma_i-\zeta\in\underline{\Delta}$ for $i=1,2,3$. We must then have for every $i$ (the indices being taken modulo $3$):
\[\beta_i=\gamma_i+\delta'_i=\gamma_{i-1}+\delta_{i-1}=\zeta+\delta'_i+\delta_{i-1},\]
hence:
\[\gamma_i=\zeta+\delta_{i-1}=\zeta+\delta'_{i+1},\]
which proves $\delta_{i-1}=\delta'_{i+1}$. We thus obtain:
\[\alpha_i=\beta_i+\delta_i=\zeta+\delta_{i+1}+\delta_{i-1}+\delta_i.\]
But then all three $\alpha_i$ are equal, which is impossible.

Assume now the cycle is of level $2$. Let $\alpha_1,\alpha_2,\alpha_3$ be the vertices of the cycle and let $\beta_1,\beta_2,\beta_3$ be the elements of $\Gamma_{l-1}$ corresponding to its edges. There can be no other $\alpha\in\Gamma_l$ such that $\alpha-\beta_i\in\underline{\Delta}$ for any $i$, because if it was the case $\alpha$ and the $\alpha_i$ would form a (nonreduced) 4-cycle in the graph; we can then solve this case the same way as the case $F_4$ and $l=4$, the only difference being that the determinant of the matrix is now $\pm 2$ instead of $\pm 3$.

Assume now there are only nonreduced $3$-cycles. First we remark that, if we remove exactly one edge of every $3$-cycle in the graph, we obtain a graph without cycles and the set of $\beta$ associated to the remaining edges is still the same; hence the nomber of $\beta\in\Gamma_{l-1}$ corresponding to edges of the graph is strictly smaller than the cardinal of $\Gamma_l$, which implies by the lemma \ref{glinf} that there exists at least one $\beta\in\Gamma_{l-1}$ which is not attached to any edge, and we can proceed the same way as in the case without cycles, except in a particular case which is not covered by that proof: the case (C) when the set $\mathcal{S}$ of $\alpha\in\Gamma_l$ such that $v(u_\alpha)=v$ is such that:
\begin{itemize}
\item  every edge between an element of $\mathcal{S}$ and an elment of $\Gamma_l-\mathcal{S}$ is part of a nonreduced $3$-cycle involving exactly two elements of $\mathcal{S}$;
\item there is no $\beta\in\Gamma_{l-1}$ such that there's only one $\alpha\in\Gamma_l$ such that $\alpha-\beta\in\underline{\Delta}$ and that $\alpha$ belongs to $\mathcal{S}$.
\end{itemize}

We will prove that the case $(C)$ doesn't show up on graphs with only nonreduced $3$-cycles. By removing one edge from every $3$-cycle, we see that $\mathcal{S}$ must contain at least two pseudo-leaves, that is two roots $\alpha$ such that there is exactly one $\beta\in\Gamma_{l-1}$ satisfying $\alpha-\beta\in\underline{\Delta}$. (Note that such a pseudo-leaf is not necessarily a leaf of $\mathcal{G}_l$ since it may be part of a nonreduced $3$-cycle.)

First we'll prove the following lemma:

\begin{lemme}
Let $\alpha$ be an element of $\underline{\Phi}^-$. The following conditions are equivalent:
\begin{itemize}
\item $\alpha$ is a pseudo-leaf of $\Gamma_l$, where $l=l(\alpha)$;
\item there is exactly one element $\delta$ of $\underline{\Delta}$ such that $\alpha-\delta
$ is a root;
\item there exists $\delta\in\underline{\Delta}$ such that $\alpha$ is a maximal element of $\underline{\Phi}'$ relatively to the root subsystem generated by $\underline{\Delta}-\{\delta\}$.
\end{itemize}
\end{lemme}

The equivalence of the first two conditions is obvious. Assume the second condition is true and the third one is not. Then there exist elements $\delta_1,\dots,\delta_r$ of $\underline{\Delta}-\{\delta\}$ such that $\alpha-\sum_{i=1}^r\delta_i$ is a root. By \cite[1, prop. 19]{bou}, there exists an $i$ such that $\alpha-\delta_i$ is a root, hence a contradiction. The converse implication is obvious. $\Box$

We will now examine the sets $\Gamma_l$ containing at least two such elements. We still use the notations of \cite[4]{bou}.

\begin{itemize}
\item Assume $\underline{\Phi}$ is of type $D_n$. The only elements of $\underline{\Phi}$ which satisfy the conditions of the lemma and are not simple roots are the roots of the form $\varepsilon_i+\varepsilon_{i+1}$, which are of respective length $2(n-i)$. No two of them are of the same length.

\item Assume now $\underline{\Phi}$ is of type $E_6$. We'll write $(abcdef)$ for the root $\alpha=a\delta_1+b\delta_2+c\delta_3+d\delta_4+e\delta_5+f\delta_6$. The elements of $\underline{\Phi}^-$ satisfying the conditions of the lemma and which are not simple roots are ($i$ being the unique index such that $\alpha-\delta_i$ is a root):
\begin{itemize}
\item for $i=2$, $(122321)$;
\item for $i=3$, $(112210)$;
\item for $i=4$, $(011210)$ and $(112321)$;
\item for $i=5$, $(011221)$.
\end{itemize}
The only two which are of the same length are $(112210)$ and $(011221)$, for which $l=7$; the only other root of length $7$ in $\underline{\Phi}^-$ is $(1.1.1.2.1.1)$. The graph is the following one:

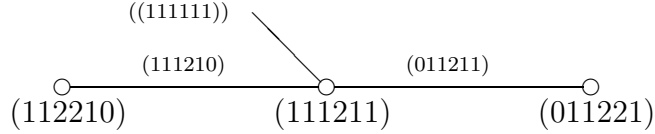
\begin{picture}(300,70)(0,0)
\put(45,20){\circle{6}}
\put(25,7){$(112210)$}
\put(48,20){\line(1,0){94}}
\put(75,26){\tiny $(111210)$}
\put(145,20){\circle{6}}
\put(125,7){$(111211)$}
\put(148,20){\line(1,0){94}}
\put(175,26){\tiny $(011211)$}
\put(245,20){\circle{6}}
\put(225,7){$(011221)$}
\put(143,22){\line(-1,1){26}}
\put(70,46){\tiny $((111111))$}
\end{picture}

 There is no $3$-cycle in $\mathcal{G}_l$ in this case.

\item Assume now $\underline{\Phi}$ is of type $E_7$. We'll write $(abcdefg)$ for the root $\alpha=a\delta_1+b\delta_2+c\delta_3+d\delta_4+e\delta_5+f\delta_6+g\delta_7$. The elements of $\underline{\Phi}^-$ satisfying the conditions of the lemma and which are not simple roots are:
\begin{itemize}
\item for $i=1$, $(2234321)$;
\item for $i=2$, $(1223210)$;
\item for $i=3$, $(1122100)$ and $(1234321)$;
\item for $i=4$, $(0112100)$, $(1123210$ and $(1224321)$;
\item for $i=5$, $0112210)$ and $(1123321)$;
\item for $i=6$, $(0112221)$.
\end{itemize}
The only two of the same length are $(1122100)$ and $(0112210)$; the other roots of length $7$ in $\underline{\Phi}^-$ are $(1112110)$, $(0112111)$ and $(1111111)$. Since this graph contains a $4$-cycle, we'll check it later. 

\item Assume now $\underline{\Phi}$ is of type $E_8$. We'll write $(abcdefgh)$ for the root $\alpha=a\delta_1+b\delta_2+c\delta_3+d\delta_4+e\delta_5+f\delta_6+g\delta_7+h\delta_8$. The elements of $\underline{\Phi}^-$ satisfying the conditions of the lemma and which are not simple roots are:
\begin{itemize}
\item for $i=1$, $(22343210)$;
\item for $i=2$, $(12232100)$ and $(13354321)$;
\item for $i=3$, $(11221000)$, $(12343210)$ and $23454321)$;
\item for $i=4$, $(01121000)$, $(11232100)$, $(12243210)$, $(12354321)$ and

$(23464321)$;
\item for $i=5$, $(01122100)$, $(11233210)$, $(12344321)$ and $(23465321)$;
\item for $i=6$, $(01122210)$, $(11233321)$ and $(23465421)$;
\item for $i=7$, $(01122221)$ and $(23465431)$;
\item for $i=8$, $(23465432)$.
\end{itemize}
There are three pairs of roots of the same length, respectively for lengths $7$, $11$ and $16$. The graph $\mathcal{G}_7$ contains a $6$-cycle and the graph $\mathcal{G}_{11}$ a $4$-cycle; we'll deal with those two graphs later. The graph $\mathcal{G}_{16}$ contains, in addition to $(12343210)$ and $(11233321)$, $(12243211)$ and $(12233221)$:

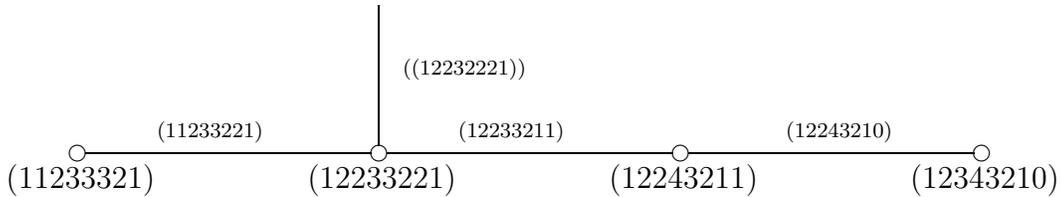
\begin{picture}(370,100)(0,0)
\put(2,20){\circle{6}}
\put(-25,7){$(11233321)$}
\put(5,20){\line(1,0){108}}
\put(32,26){\tiny $(11233221)$}
\put(116,20){\circle{6}}
\put(89,7){$(12233221)$}
\put(119,20){\line(1,0){108}}
\put(146,26){\tiny $(12233211)$}
\put(230,20){\circle{6}}
\put(203,7){$(12243211)$}
\put(233,20){\line(1,0){108}}
\put(270,26){\tiny $(12243210)$}
\put(344,20){\circle{6}}
\put(317,7){$(12343210)$}
\put(116,23){\line(0,1){53}}
\put(125,50){\tiny $((12232221))$}
\end{picture}

There's no cycle in this graph.
\end{itemize}

Assume now $\mathcal{G}_l$ contains a $4$-cycle and no strictly larger cycles. Let $\alpha_1,\dots,\alpha_4$ (resp. $\beta_1,\dots,\beta_4$) be the corresponding elements of $\Gamma_l$ (resp. $\Gamma_{l-1}$). We already know that the $\beta_i$ can't be all equal to each other. We will also show that they can't be all different either.

By an obvious induction on the level of the cycle, we see that it is enough to prove the result when the cycle of the $\beta_i$ is nonreduced. Let $\gamma_1,\dots,\gamma_4$ be the elements of $\Gamma_{l-2}$ corresponding to its edges (for every $i$, $gamma_i$ corresponds to the edge between $\beta_i$ and $\beta_{i-1}$); exactly two consecutive of them are equal to each other, say for example $\gamma_1=\gamma_2$.
But then $(\beta_4,\beta_1,\beta_2)$ is a nonreduced $3$-cycle in $\mathcal{G}_{l-1}$, which implies that there exists $\zeta\in\Gamma_{l-3}$ such that $\gamma_1=\gamma_2=\zeta+\delta_3$, $\gamma_3=\zeta+\delta_1$ and $\gamma_4=\zeta+\delta_2$; we finally obtein that both $\alpha_2$ and $\alpha_3$ must be equal to $\zeta+\delta_1+\delta_2+\delta_3$, hence a contradiction.

If $\mathcal{G}_l$ contains a reduced $3$-cycle which is not part of a nonreduced $4$-cycle, we proceed as in the case in which there are only $3$-cycles, one of them being reduced. If all reduced $3$-cycles are part of nonreduced $4$-cycles, then the nomber of $\beta\in\Gamma_{l-1}$ corresponding to edges of the graph is strictly smaller than the cardinal of $\Gamma_l$, and as previously, we proceed as in the cycle-free case except for the case $(C)$.

In particular we now have completely solved the case $D_n$: assuming we are in that case, when $l\leq n-2$, $\Gamma_l=\{\varepsilon_{l+i}-\varepsilon_i|1\leq i\leq n-l\}\cup\{-\varepsilon_{n-l+i}-\varepsilon_{n-i}|0\leq i\leq[\frac{l-1}2]\}$, which the leftmost two sets being eventually empty, and the graph $\mathcal{G}_l$ is the following one:

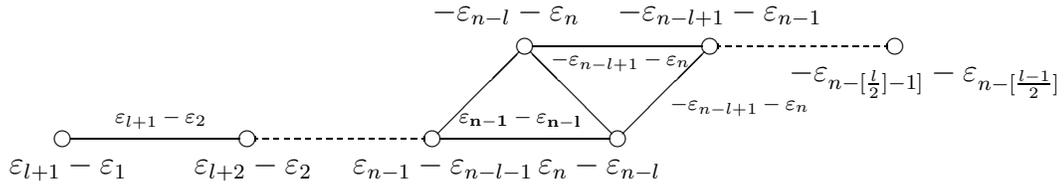
\begin{picture}(400,100)(0,0)
\put(25,20){\circle{6}}
\put(5,7){$\varepsilon_{l+1}-\varepsilon_1$}
\put(28,20){\line(1,0){64}}
\put(45,26){\tiny $\varepsilon_{l+1}-\varepsilon_2$}
\put(95,20){\circle{6}}
\put(75,7){$\varepsilon_{l+2}-\varepsilon_2$}
\put(98,20){\dashbox{2}(64,0){}}
\put(165,20){\circle{6}}
\put(135,7){$\varepsilon_{n-1}-\varepsilon_{n-l-1}$}
\put(168,20){\line(1,0){64}}
\put(175,26){\tiny $\mathbf{\varepsilon_{n-1}-\varepsilon_{n-l}}$}
\put(235,20){\circle{6}}
\put(205,7){$\varepsilon_n-\varepsilon_{n-l}$}
\put(167,22){\line(1,1){31}}
\put(233,22){\line(-1,1){31}}
\put(200,55){\circle{6}}
\put(165,65){$-\varepsilon_{n-l}-\varepsilon_n$}
\put(203,55){\line(1,0){64}}
\put(210,48){\tiny $-\varepsilon_{n-l+1}-\varepsilon_n$}
\put(270,55){\circle{6}}
\put(235,65){$-\varepsilon_{n-l+1}-\varepsilon_{n-1}$}
\put(237,22){\line(1,1){31}}
\put(255,31){\tiny $-\varepsilon_{n-l+1}-\varepsilon_n$}
\put(273,55){\dashbox{2}(64,0){}}
\put(340,55){\circle{6}}
\put(300,42){$-\varepsilon_{n-[\frac l2]-1]}-\varepsilon_{n-[\frac {l-1}2]}$}
\end{picture}

For $3\leq l\leq n-2$, this graph contains a $4$-cycle $\varepsilon_{n-1}-(\varepsilon_{n-l-1},\varepsilon_n-\varepsilon_{n-l},-\varepsilon_{n-l+1}-\varepsilon_{n-1},-\varepsilon_{n-l}-\varepsilon_n)$ whose reduced part is $(\varepsilon_n-\varepsilon_{n-l},-\varepsilon_{n-l+1}-\varepsilon_{n-1},-\varepsilon_{n-l}-\varepsilon_n)$, and no other vertices are part of any cycle. (For $l=2$ or $l=n-1$ there's only the reduced $3$-cycle, and for $l\geq n$ there are no cycles at all). Moreover, for $l$ even there are no pseudo-leaves in $\mathcal{G}_l$, and for $l$ odd the only pseudo-leaf in $\mathcal{G}_l$ is $-\varepsilon_{n-\frac{l+1}2}+\varepsilon_{n-\frac{l-1}2}$; this implies that if we are in the case (C), $\mathcal{S}$ must contain the reduced $3$-cycle but not $\varepsilon_{n-1}-\varepsilon_{n--l1}$, and we can then conclude as in the reduced $3$-cycle case.

Now the only case we haven't fully solved yet is $E_n$; since $n$ is then at most $8$, it follows then from the lemma \ref{glinf} that the cardinal of $\Gamma_l$ is at most $7$, which means in particular that we can only have cycles of length up to $7$.

But we haven't finished to deal with the $4$-cycles yet. Assume now we are in the case $(C)$ and the subgraph $\mathcal{S}$ contains no $4$-cycles. If it contains a reduced $3$-cycle we can proceed as in the corresponding case; we thus can assume it doesn't either. According to the previous discussion, there are only two occurrences of this case: either $\underline{\Phi}=E_7$ and $l=7$ or $\underline{\Phi}=E_8$ and $l=11$. In the first case, the graph $\mathcal{G}_l$ is the following one:

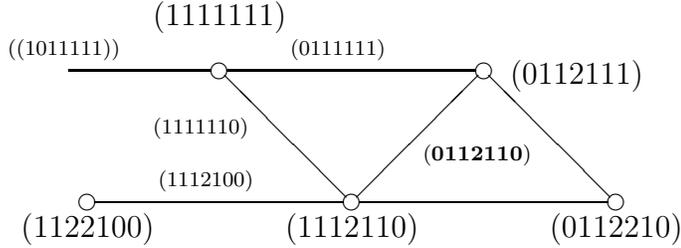
\begin{picture}(400,100)(0,0)
\put(45,20){\circle{6}}
\put(20,7){$(1122100)$}
\put(48,20){\line(1,0){94}}
\put(72,26){\tiny $(1112100)$}
\put(145,20){\circle{6}}
\put(120,7){$(1112110)$}
\put(148,20){\line(1,0){94}}
\put(172,36){\tiny $\mathbf{(0112110)}$}
\put(245,20){\circle{6}}
\put(220,7){$(0112210)$}
\put(147,22){\line(1,1){46}}
\put(243,22){\line(-1,1){46}}
\put(195,70){\circle{6}}
\put(205,65){$(0112111)$}
\put(98,70){\line(1,0){94}}
\put(122,76){\tiny $(0111111)$}
\put(95,70){\circle{6}}
\put(70,87){$(1111111)$}
\put(143,22){\line(-1,1){46}}
\put(70,46){\tiny $(1111110)$}
\put(38,70){\line(1,0){54}}
\put(15,76){\tiny $((1011111))$}
\end{picture}

The set $\mathcal{S}$ cannot contain $(1111111)$, hence cannot contain $(1112110)$ and $(0112111)$ either.  But then it contains at most one element of the reduced $3$-cycle, which contradicts the definition of the case $(C)$.

In the second case, the graph $\mathcal{G}_l$ is the following one:

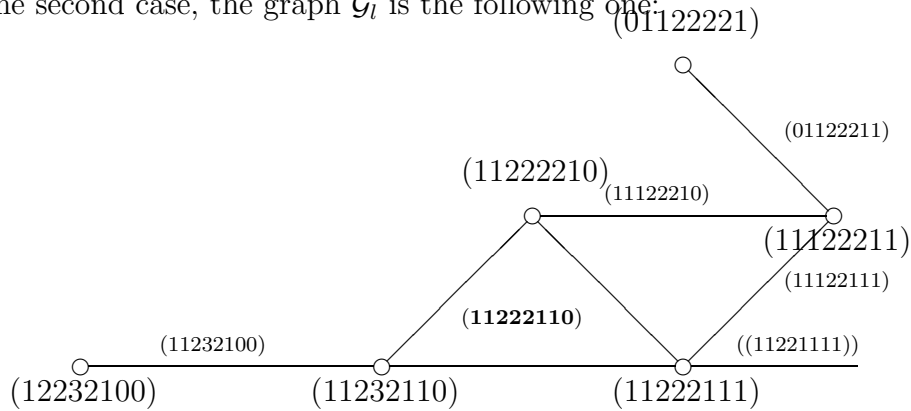
\begin{picture}(400,150)(0,0)
\put(52,20){\circle{6}}
\put(25,7){$(12232100)$}
\put(55,20){\line(1,0){108}}
\put(82,26){\tiny $(11232100)$}
\put(166,20){\circle{6}}
\put(139,7){$(11232110)$}
\put(169,20){\line(1,0){108}}
\put(196,36){\tiny $\mathbf{(11222110)}$}
\put(280,20){\circle{6}}
\put(253,7){$(11222111)$}
\put(168,22){\line(1,1){53}}
\put(278,22){\line(-1,1){53}}
\put(223,77){\circle{6}}
\put(196,90){$(11222210)$}
\put(283,20){\line(1,0){63}}
\put(300,26){\tiny $((11221111))$}
\put(226,77){\line(1,0){108}}
\put(250,83){\tiny $(11122210)$}
\put(337,77){\circle{6}}
\put(310,64){$(11122211)$}
\put(282,22){\line(1,1){53}}
\put(318,50){\tiny $(11122111)$}
\put(335,79){\line(-1,1){53}}
\put(318,107){\tiny $(01122211)$}
\put(280,134){\circle{6}}
\put(253,147){$(01122221)$}
\end{picture}

Obviously, a subset $\mathcal{S}$ of $\Gamma_l$ satisfying the required conditions and containing both $(12232100)$ and $(01122221)$ cannot exist.

Assume now $\mathcal{S}$ contains a $4$-cycle; $\mathcal{G}$ must then also contain a nonreduced $3$-cycle not entirely included in $\mathcal{S}$, and since $\mathcal{G}$ contains no cycles of length greater than $4$, we are in one of the following two cases:

\begin{itemize}
\item that nonreduced $3$-cycle contains two opposites vertices of the $4$-cycle;
\item that nonreduced $3$-cycle shares at most one vertex with the $4$-cycle.
\end{itemize}

Assume we are in the first case; let $(\alpha_1,\alpha_2,\alpha_3,\alpha_4)$ be the $4$-cycle and $(\alpha_1,\alpha_3,\alpha_5)$ be the nonreduced $3$-cycle. Then $\mathcal{G}_l$ contains the $4$-cycles $(\alpha_1,\alpha_2;\alpha_3,\alpha_5)$ and $(\alpha_1,\alpha_3,\alpha_4,\alpha_5)$ as well, hence $(\alpha_1,\alpha_2,\alpha_3)$ and $(\alpha_1,\alpha_3,\alpha_4)$ must both be reduced $3$-cycles, which is impossible since $(\alpha_1,\alpha_2,\alpha_3,\alpha_4)$ cannot be reduced..

Assume now we are in the second case. We have the following result:

\begin{lemme}
Assume $\mathcal{G}_l$ contains two $3$-cycles $(\alpha_1,\alpha_2;\alpha_3)$ and $(\alpha_1,\alpha_4,\alpha_5)$ sharing the vertex $\alpha_1$, all the $\alpha_i$ being different.. Then there is a $5$-cycle in $\mathcal{G}$ whose vertices are the $\alpha_i$.
\end{lemme}

It is enough to prove that there is an edge between $\alpha_i$, $i=2$ or $3$, and $\alpha_j$, $j=4$ or $5$. We will first prove the following lemma:

\begin{lemme}
Let $\beta$ be an element of $\Gamma_{l-1}$ corresponding to an edge between two elements $\alpha$ and $\alpha'$ of $\Gamma_l$, and write $\delta=\alpha-\beta$, $\delta'=\alpha'-\beta$. Then $\alpha+\delta'=\alpha'+\delta\in\Phi$.
\end{lemme}

We have $(\alpha,\delta')=(\beta,\delta')+(\delta,\delta')\leq -1$, which is enough to prove the assertion. $\Box$

For $i=2,3$ (resp. $i=4,5$), write $\delta_i=\alpha_i-\beta$ (resp. $\alpha_i-\beta'$). Then for every $i$, according to the previous lemma, $\alpha_1+\delta_i\in\Phi$. This is only possible if $\delta_i=\delta_j$ for some $(i,j)$, $i<j$; we already know that $(i,j)$ cannot be $(2,3)$ or $(4,5)$, and by symmetry we may assume $(i,j)$ is any other pair, say for example $(2,4)$. We then have:
\[\alpha_4=\beta"+\delta_2=\alpha_1-\delta'_1+\delta_2=\beta+\delta_1-\delta'_1+\delta_2=\alpha_2+\delta_1-\delta'_1.\]
We then only have to check that $\alpha_2-\delta'_1$ is a root. We have:
\[(\alpha_2,\delta'_1)=(\alpha_4,\delta'_1)+(\delta'_1,\delta'_1)-(\delta_1,\delta'_1)\geq -1+2-0=1,\]
which proves the desired assertion. $\Box$

According to that lemma, since we already know that the $4$-cycle can be broken into two $3$-cycles, it must be disjoint from the nonreduced $3$-cycle, which implies that the cardinal of $\Gamma_l$ is precisely $7$ and $\mathcal{S}$ contains all elements of $\Gamma_l$ but one, that one belonging to the isolated $3$-cycle. In particular there must be an element of $\mathcal{S}$ attached to only one $\beta\in\Gamma_{l-1}$, which, according to the previous discussion, leaves only two possible cases: $\underline{\Phi}$ is of type $E_8$ and either $l=5$ or $l=7$. In both cases, $\mathcal{G}$ contains a cycle of length greater than $4$, which leads to a contradiction.

Before considering the cases of larger cycles, we will prove the following result:

\begin{prop}\label{triang}
Let $(\alpha_1,\dots,\alpha_r)$ be a cycle in $\mathcal{G}_l$. There exist $r-3$ edges in $\mathcal{G}_l$ not belonging to this cycle and dividing it into $r-2$ $3$-cycles of the form $(\alpha_i,\alpha_{i+1},\alpha_j)$. Moreover, two $3$-cycles sharing a common edge can be neither both reduced nor both nonreduced.
\end{prop}

The second assertion comes fron the discussion on $4$-cycles: two $3$-cycles sharing a common edge form a $4$-cycle, and such a cycle can be neither reduced (hence the $3$-cycles are not both reduced) nor of level $1$ (hence the $3$-cycles are not both nonreduced).

We will now prove the first assertion by double induction, first on $l$, then on $r$. Assume first $l=2$: since we have already proved that $\Gamma_2$ cannot contain reduced cycles, there can be only $3$-cycles and there's nothing to prove. Assume now $l>2$ and $r>3$ (if $r=3$, there is nothing to prove) and let $\beta_1,\dots,\beta_r$ be the elements of $\Gamma_{l-1}$ attached to the edges of the cycle. If the cycle is nonreduced, say for example if $\beta_1=\beta_2$, then there is an edge between $\alpha_1$ and $\alpha_3$, hence $(\alpha_1,\dots,\alpha_3)$ is a $3$-cycle and we conclude by applying the induction hypothesis to the $r-1$-cycle $(\alpha_1,\alpha_3,\alpha_4,\dots,\alpha_r)$.

Now we will assume the cycle is reduced. By the induction hypothesis, the $r$-cycle $(\beta_1,\dots,\beta_r)$ in $\mathcal{G}_{l-1}$ can be divided into $r-2$ $3$-cycles, and since these cycles contain all $r$ edges of the $r$-cycle, at least one contains two of them, hence is made of three consecutive $\beta_i$, say $\beta_1,\beta_2$ and $\beta_3$. Assume that $3$-cycle is reduced; then $l\geq 4$ and there exist $\zeta\in\Gamma_{l-3}$ and $\delta_1,\delta_2,\delta_3\in\underline{\Delta}$ such that $\beta_i=\zeta+\delta_i+\delta_{i-1}$ for each $i$ (the indices being taken modulo $3$). But then we must have $\alpha_1=\alpha_2=\zeta+\delta_1+\delta_2+\delta_3$, hence a contradiction. Hence that $3$-cycle is nonreduced, which implies that there exist $\gamma\in\Gamma_{l-2}$ and $\delta_1,\delta_2,\delta_3\in\underline{\Delta}$ such that $\beta_i=\gamma+\delta_i$ for each $i$. We deduce from this that we have $\alpha_1=\zeta+\delta_1+\delta_2$ and $\alpha_2=\zeta+\delta_2+\delta_3$.

Consider now the character $\alpha'=\zeta+\delta_1+\delta_3$; we can easily see that it is a root, and there are edges $(\alpha_1,\alpha')$ and $(\alpha_2,\alpha')$ in $\Gamma_l$ associated respectively to $\beta_1$ and $\beta_3$. If $\alpha'=\alpha_i$ for some $i$, then we have a $3$-cycle $(\alpha_1,\alpha_2,\alpha_i)$ and we conclude by applying the induction hypothesis to the cycles $(\alpha_1,\alpha_i,\dots,\alpha_r)$ and $(\alpha_2,\dots,\alpha_i)$.

Assume now $\alpha'$ does not belong to the set of the $\alpha_i$; we will prove that it leads to a contradiction. We have two edges $(\alpha_r,\alpha')$ and $(\alpha_1,\alpha')$ associated to $\beta_1$, and two edges $(\alpha_2,\alpha')$ and $(\alpha_3,\alpha')$ associated to $\beta_3$. If there is an edge between $\alpha_3$ and $\alpha_r$, associated to some element $\beta'$ of $\Gamma_{l-1}$, then $(\alpha_3,\alpha',\alpha_r)$ is a reduced $3$-cycle, hence $(\beta_1,\beta_3,\beta')$ is a nonreduced $3$-cycle in $\mathcal{G}_{l-1}$; we must then have $\beta'=\beta_2$, but then $\alpha_r=\alpha_2$, hence a contradiction.

Assume now there is no such edge. From the induction hypothesis applied to $(\alpha',\alpha_3,\dots,\alpha_r)$ we deduce that there exists indices $i_0=3<i_1<\dots<i_s=r$, $s\geq 2$, such that $(\alpha',\alpha_{i_{j-1}},\alpha_{i_j})$ is a $3$-cycle for every $j\in\{1,\dots,s\}$. From the discussion about $4$-cycles we deduce that neither two nonreduced $3$-cycles nor two reduced $3$-cycles can share a common edge; hence the reduced and nonreduced cycles must alternate, and the first and last must be reduced, which implies that $s$ is odd; since the cardinal of $\Gamma_l$ is not greater than $7$, we must have $s=3$, $r=6$ and $i_j=j+3$ for every $j$.

Consider then the cycle $(\beta_1,\dots,\beta_6)$ in $\mathcal{G}_{l-1}$. We deduce from the above conditions that $(\beta_1,\beta_2,\beta_3)$, $(\beta_3,\beta_4,\beta_5)$ and $(\beta_5,\beta_6,\beta_1)$  are nonreduced $3$-cycles; the cycle $(\beta_1,\beta_3,\beta_5)$ must then be a reduced $3$-cycle. We will prove that this case is impossible. Let $\gamma_1,\gamma_3,\gamma_5$ be the elements of $\Gamma_{l-2}$ associated respectively to the edges of those three nonreduced $3$-cycles, and let $\zeta$ be the element of $\Gamma_{l-3}$ associated to the edges of the nonreduced $3$-cycle $(\gamma_1,\gamma_3,\gamma_5)$ in $\mathcal{G}_{l-2}$. Write for $i=1,2,3$:
\[\delta_i=\gamma_{2i-1}-\zeta,\delta'_i=\beta_{2i}-\gamma_{2i-1}.\]
We have $(\zeta,\delta_i)=-1$ and $(\zeta,\delta'_i)=0$ for every $i$. Moreover, $(\gamma_{2i-1},\delta'_i)=-1$ hence $(\delta_i,\delta'_i)=-1$, and if $i\neq j$, $(\delta_i,\delta_j)=0$ and $(\alpha_i,\delta_j)=-1$ which imposes $(\delta_i,\delta'_j)=0$. Now consider the character $\chi=3\zeta+\sum_{i=1}^3(2\delta_i+\delta'_i)$; it is a sum of negative roots, hence nonzero, but we have:
\[(\chi,\chi)=9(\zeta,\zeta)+\sum_{i=1}^3(12(\zeta,\delta_i)+4(\delta_i,\delta_i)+4(\delta_i,\delta'_i)+(\delta'_i,\delta'_i))\]
\[=18+3(-12+8-4+2)=0.\]
We obtain a contradiction, which concludes the proof of the first assertion. $\Box$

Assume now $\mathcal{G}$ contains a $5$-cycle, and no larger cycle. Assume first the $5$-cycle is reduced; and define the $\alpha_i,\beta_i,\gamma_i,\delta_i,\delta'_i$ as usual; we may assume that the $3$-cycles given by the previous proposition are $(\alpha_5,\alpha_1,\alpha_2)$, $(\alpha_2,\alpha_3,\alpha_5)$ and $(\alpha_3,\alpha_4,\alpha_5)$. Then the first and third must be reduced, which implies the second is not and thus the edges $(\alpha_2,\alpha_5)$ and $(\alpha_3,\alpha_5)$ are both associated to $\beta_3$. Moreover, $(\beta_1,\beta_2,\beta_3)$ and $(\beta_3,\beta_4,\beta_5)$ are both nonreduced $3$-cycles in $\mathcal{G}_{l-1}$, hence $\gamma_1=\gamma_2$ and $\gamma_3=\gamma_4$.

We then have:
\[\beta_2=\gamma_1+\delta_1=\gamma_1+\delta'_2\]
hence $\delta_1=\delta'_2$, and similarly $\delta_3=\delta'_4$. Moreover:
\[\beta_3=\gamma_1+\delta_2=\gamma_3+\delta'_3;\]
\[\beta_5=\gamma_3+\delta_4=\gamma_5+\delta'_5;\]
\[\beta_1=\gamma_1+\delta'_1=\gamma_5+\delta_5.\]
We deduce from this that $(\beta_3,\beta_5,\beta_1)$ is a reduced $3$-cycle in $\mathcal{G}_{l-1}$; hence $\delta_2=\delta'_5$, $\delta_4=\delta'_1$ and $\delta_5=\delta'_3$, and there exists $\zeta\in\Gamma_{l-3}$ such that $\gamma_1=\zeta+\delta_5$, $\gamma_3=\zeta+\delta_2$ and $\gamma_5=\zeta+\delta_4$. We finally obtain:
\[\alpha_1=\zeta+\delta_1+\delta_4+\delta_5;\]
\[\alpha_2=\zeta+\delta_1+\delta_2+\delta_5;\]
\[\alpha_3=\zeta+\delta_2+\delta_3+\delta_5;\]
\[\alpha_4=\zeta+\delta_2+\delta_3+\delta_4;\]
\[\alpha_5=\zeta+\delta_2+\delta_4+\delta_5.\]
By a similar reasoning as in the case of the reduced $3$-cycle, we obtain that the following matrix must be invertible for every $x_1,\dots,x_5\in\overline{k}^*$: 
\[\left(\begin{array}{ccccc}\varepsilon_{\alpha_1,-\beta_1}x_1&&\varepsilon_{\alpha_2,-\beta_3}x_2\\\varepsilon_{\alpha_5,-\beta_1}x_5&\varepsilon_{\alpha_2,-\beta_2}x_2\\&&\varepsilon_{\alpha_3,-\beta_3}x_3&&\varepsilon_{\alpha_4,-\beta_5}x_4\\&\varepsilon_{\alpha_1,-\beta_2}x_1&\varepsilon_{\alpha_5,-\beta_3}x_5&\varepsilon_{\alpha_4,-\beta_4}x_4\\&&&\varepsilon_{\alpha_3,-\beta_4}x_3&\varepsilon_{\alpha_5,-\beta_5}x_5\end{array}\right).\]
The determinant of this matrix is $\pm 3 x_1x_2x_3x_4x_5$, which is invertible since we have assumed $p\neq 3$ for the case $E_n$.

Assume now the graph contains a nonreduced $5$-cycle and no larger cycle. With the same notations as above, we may assume that $(\alpha_5,\alpha_1,\alpha_2)$ is nonreduced, which implies $(\alpha_2,\alpha_3,\alpha_5)$ is reduced and $(\alpha_3,\alpha_4,\alpha_5)$ is not. Once again we only have to check the case $(C)$; we deduce from its definition that the intersection of the subgraph $\mathcal{S}$ with the $5$-cycle is one of the following parts:
\begin{itemize}
\item the reduced $3$-cycle,
\item $\{\alpha_4,\alpha_5,\alpha_1\}$;
\item one of the $4$-cycles, say $(\alpha_1,\alpha_2,\alpha_3,\alpha_5)$.
\end{itemize}
In the first case, we conclude as in the case of reduced $3$-cycles. The second case would imply that $\mathcal{S}$ contains two pseudo-leaves, and we know from the discussion of graphs satisfying that condition that none of them contains a $5$-cycle and no larger cycle. We thus only have to examine the latter case: it implies that exactly one of the elements of $\mathcal{S}$ is a pseudo-leaf, and thus, since $\Gamma_l$ must contain at least $5$ elements, that we are in one of the following cases: $\underline{\Phi}$ is of type $E_7$ and $l=5$, or $\underline{\Phi}$ is of type $E_8$ and $l\in\{5,9,10,13\}$.

The first three contain a cycle of length greater than $6$; we'll deal with them later. In the case $E_8$ and $l=10$, the graph is the following one:

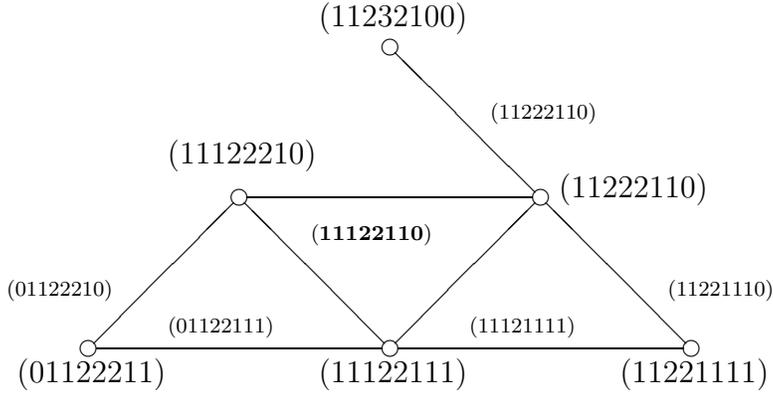
\begin{picture}(400,160)(0,0)
\put(52,20){\circle{6}}
\put(25,7){$(01122211)$}
\put(55,20){\line(1,0){108}}
\put(82,26){\tiny $(01122111)$}
\put(166,20){\circle{6}}
\put(139,7){$(11122111)$}
\put(169,20){\line(1,0){108}}
\put(196,26){\tiny $(11121111)$}
\put(280,20){\circle{6}}
\put(253,7){$(11221111)$}
\put(54,22){\line(1,1){53}}
\put(21,40){\tiny $(01122210)$}
\put(164,22){\line(-1,1){53}}
\put(109,77){\circle{6}}
\put(82,90){$(11122210)$}
\put(168,22){\line(1,1){53}}
\put(278,22){\line(-1,1){53}}
\put(271,40){\tiny $(11221110)$}
\put(223,77){\circle{6}}
\put(230,77){$(11222110)$}
\put(112,77){\line(1,0){108}}
\put(136,61){\tiny $\mathbf{(11122110)}$}
\put(221,79){\line(-1,1){53}}
\put(204,107){\tiny $(11222110)$}
\put(166,134){\circle{6}}
\put(139,142){$(11232100)$}
\end{picture}

This graph contains a $5$-cycle, but it is reduced. In the case $E_8$ and $l=13$, the graph is the following one:

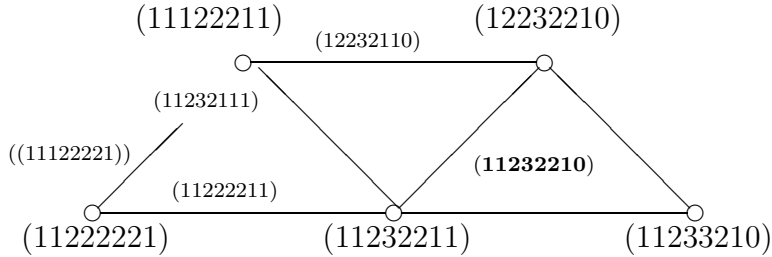
\begin{picture}(400,120)(0,0)
\put(52,20){\circle{6}}
\put(25,7){$(11222221)$}
\put(55,20){\line(1,0){108}}
\put(82,26){\tiny $(11222211)$}
\put(166,20){\circle{6}}
\put(139,7){$(11232211)$}
\put(169,20){\line(1,0){108}}
\put(196,36){\tiny $\mathbf{(11232210)}$}
\put(280,20){\circle{6}}
\put(253,7){$(11233210)$}
\put(168,22){\line(1,1){53}}
\put(278,22){\line(-1,1){53}}
\put(223,77){\circle{6}}
\put(196,90){$(12232210)$}
\put(54,22){\line(1,1){32}}
\put(20,40){\tiny $((11122221))$}
\put(112,77){\line(1,0){108}}
\put(136,83){\tiny $(12232110)$}
\put(168,22){\line(-1,1){53}}
\put(75,60){\tiny $(11232111)$}
\put(109,77){\circle{6}}
\put(68,90){$(11122211)$}
\end{picture}

This graph contains no $5$-cycle.

Assume now that the graph contains a $6$-cycle and no larger cycle. If this cycle is reduced, then amonng the four $3$-cycles given by proposition \ref{triang}, at least three must be reduced (if it was not the case, at most five of the $\beta_i$ would be different); the only possibility (up to a circular permutation of the indices) is that we have reduced $3$-cycles $(\alpha_1,\alpha_2,\alpha_3)$, $(\alpha_3,\alpha_4,\alpha_5)$ and $(\alpha_5,\alpha_6,\alpha_1)$, and thus a nonreduced $3$-cycle $(\alpha_1,\alpha_3,\alpha_5)$, whose edges are associated to an element $\beta'$ of $\Gamma_{l-1}$ which must necessarily be different from all $\beta_i$, $i\in\{1,\dots,6\}$. But then, in the graph $\mathcal{G}_{l-1}$, we have six $3$-cycles $(\beta_i,\beta_{i+1},\beta')$, which are reduced if $i$ is even and nonreduced if $i$ is odd, and this case has already explicitely been ruled out during the proof of the proposition \ref{triang}.

Hence the $6$-cycle must be nonreduced. Once again we only have to check the case $(C)$. Assume we are in this case; then either the intersection of $\mathcal{S}$ with the $6$-cycle contains a reduced cycle which is not part of a larger cycle in $\mathcal{S}$, in which case we conclude as in the corresponding case, or $\mathcal{S}$ contains a pseudo-leaf of the graph. In that last case, we have only a few cases to check, and we will check them directly.
\begin{itemize}
\item Assume $\underline{\Phi}$ is of type $E_7$ and $l=5$. The graph is the following one:

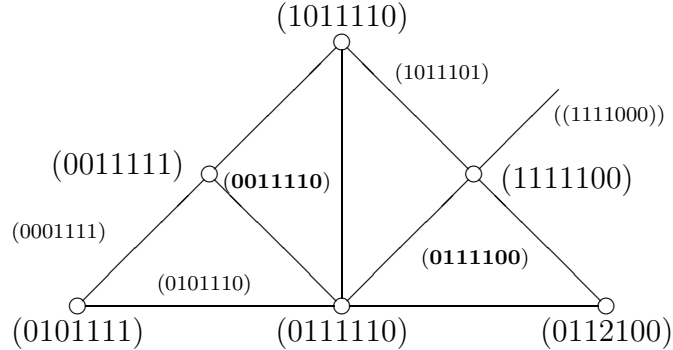
\begin{picture}(360,150)(0,0)
\put(45,20){\circle{6}}
\put(20,7){$(0101111)$}
\put(48,20){\line(1,0){94}}
\put(75,26){\tiny $(0101110)$}
\put(145,20){\circle{6}}
\put(120,7){$(0111110)$}
\put(148,20){\line(1,0){94}}
\put(175,36){\tiny $\mathbf{(0111100)}$}
\put(245,20){\circle{6}}
\put(220,7){$(0112100)$}
\put(147,22){\line(1,1){46}}
\put(243,22){\line(-1,1){46}}
\put(195,70){\circle{6}}
\put(205,65){$(1111100)$}
\put(47,22){\line(1,1){46}}
\put(20,46){\tiny $(0001111)$}
\put(143,22){\line(-1,1){46}}
\put(95,70){\circle{6}}
\put(35,70){$(0011111)$}
\put(193,72){\line(-1,1){46}}
\put(165,106){\tiny $(1011101)$}
\put(97,72){\line(1,1){46}}
\put(145,23){\line(0,1){94}}
\put(100,65){\tiny $\mathbf{(0011110)}$}
\put(145,120){\circle{6}}
\put(120,126){$(1011110)$}
\put(197,72){\line(1,1){30}}
\put(225,90){\tiny $((1111000))$}
\end{picture}

The pseudo-leaf is $(0112100)$, which is part of the nonreduced $3$-cycle $((0112100),)(1111100),(0111110))$.,Since the element $(1111000)$ of $\Gamma_4$ is attached to $(1111100)$ and no other element of $\Gamma_5$, $\mathcal{S}$ cannot contain $(1111100)$, hence must be $\{(0112100),(0111110),(0101111),(0011111)\}$, which contains the reduced $3$-cycle $((0111110),(0011111),(0101111))$ and no larger cycle, and we conclude as in the case of a reduced $3$-cycle.
\item Assume $\underline{\Phi}$ is of type $E_8$ and $l=7$. The graph is the following one:

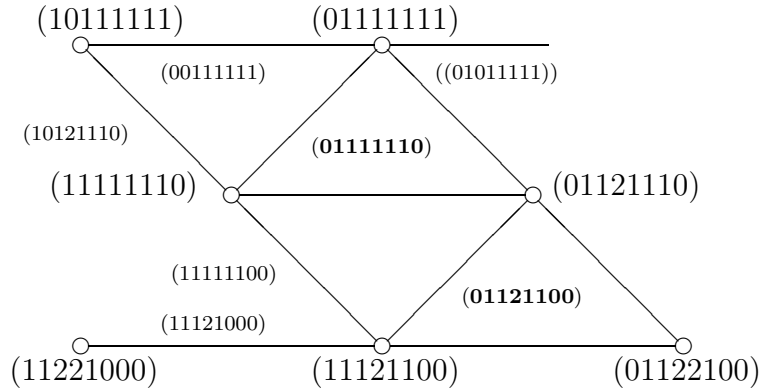
\begin{picture}(350,160)(0,0)
\put(52,20){\circle{6}}
\put(25,7){$(11221000)$}
\put(55,20){\line(1,0){108}}
\put(82,26){\tiny $(11121000)$}
\put(166,20){\circle{6}}
\put(139,7){$(11121100)$}
\put(169,20){\line(1,0){108}}
\put(196,36){\tiny $\mathbf{(01121100})$}
\put(280,20){\circle{6}}
\put(253,7){$(01122100)$}
\put(164,22){\line(-1,1){53}}
\put(86,45){\tiny $(11111100)$}
\put(109,77){\circle{6}}
\put(40,77){$(11111110)$}
\put(168,22){\line(1,1){53}}
\put(278,22){\line(-1,1){53}}
\put(223,77){\circle{6}}
\put(230,77){$(01121110)$}
\put(112,77){\line(1,0){108}}
\put(139,93){\tiny $\mathbf{(01111110)}$}
\put(111,79){\line(1,1){53}}
\put(221,79){\line(-1,1){53}}
\put(166,134){\circle{6}}
\put(139,140){$(01111111)$}
\put(169,134){\line(1,0){60}}
\put(186,121){\tiny $((01011111))$}
\put(55,134){\line(1,0){108}}
\put(82,121){\tiny $(00111111)$}
\put(107,79){\line(-1,1){53}}
\put(30,98){\tiny $(10121110)$}
\put(52,134){\circle{6}}
\put(35,140){$(10111111)$}
\end{picture}

There are two pseudo-leaves: $(11221000)$ and $(01122100)$, and there are edges between both of them and $(11121100)$. If $\mathcal{S}$ contains $(11121100)$, then by following the edges which are not part of nonreduced $3$-cycles, we see that $\mathcal{S}$ must contain as well $(11111110)$, $(10111111)$ and $(01111111)$; but the element $(01011111)$ of $\Gamma_6$ is attached to $(01111111)$ and no other element of $\Gamma_7$, hence a contradiction. If $\mathcal{S}$ doesn't contain $(11121100)$, then it must contain $(01122100)$, and by following the edges we see once again that it must contain $(01111111)$; we thus reach the same contradiction.
\item Assume $\underline{\Phi}$ is of type $E_8$ and $l=9$. The graph is the following one:

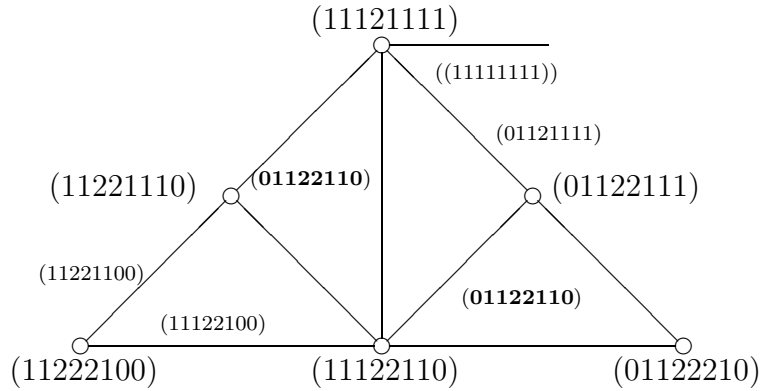
\begin{picture}(350,160)(0,0)
\put(52,20){\circle{6}}
\put(25,7){$(11222100)$}
\put(55,20){\line(1,0){108}}
\put(82,26){\tiny $(11122100)$}
\put(166,20){\circle{6}}
\put(139,7){$(11122110)$}
\put(169,20){\line(1,0){108}}
\put(196,36){\tiny $\mathbf{(01122110})$}
\put(280,20){\circle{6}}
\put(253,7){$(01122210)$}
\put(54,22){\line(1,1){53}}
\put(36,45){\tiny $(11221100)$}
\put(164,22){\line(-1,1){53}}
\put(109,77){\circle{6}}
\put(40,77){$(11221110)$}
\put(168,22){\line(1,1){53}}
\put(278,22){\line(-1,1){53}}
\put(223,77){\circle{6}}
\put(230,77){$(01122111)$}
\put(111,79){\line(1,1){53}}
\put(166,23){\line(0,1){108}}
\put(116,81){\tiny $\mathbf{(01122110)}$}
\put(221,79){\line(-1,1){53}}
\put(209,98){\tiny $(01121111)$}
\put(166,134){\circle{6}}
\put(139,140){$(11121111)$}
\put(169,134){\line(1,0){60}}
\put(186,121){\tiny $((11111111))$}
\end{picture}

The pseudo-leaf is $(01122210)$. Since the element $(11111111)$ of $\Gamma_8$ is attached to $(11121111)$ and no other vertex, $\mathcal{S}$ cannot contain $(11121111)$, hence must be $\{(01122210),(11,122110),(11221110),(11222100)\}$, which contains the reduced $3$-cycle $((11122110),(11221110),(11222100))$ and no larger cycle, and we conclude as in the case of a reduced $3$-cycle.
\end{itemize}

Now we finally consider the cases when the graph contains a $7$-cycle. Since there are only two occurrences of such cycles, one being reduced and the other one nonreduced, we will examine them directly.

First assume $\underline{\Phi}$ is of type $E_8$ and $l=6$. The graph is the following one:

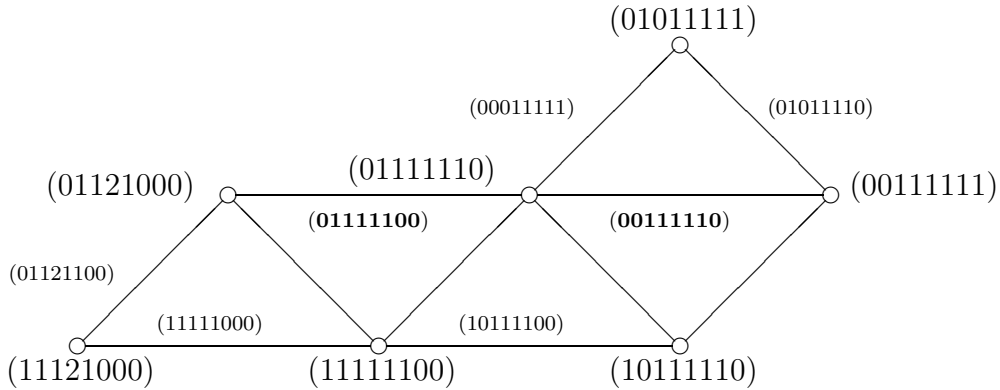
\begin{picture}(400,160)(0,0)
\put(52,20){\circle{6}}
\put(25,7){$(11121000)$}
\put(55,20){\line(1,0){108}}
\put(82,26){\tiny $(11111000)$}
\put(166,20){\circle{6}}
\put(139,7){$(11111100)$}
\put(169,20){\line(1,0){108}}
\put(196,26){\tiny $(10111100)$}
\put(280,20){\circle{6}}
\put(253,7){$(10111110)$}
\put(54,22){\line(1,1){53}}
\put(26,45){\tiny $(01121100)$}
\put(164,22){\line(-1,1){53}}
\put(109,77){\circle{6}}
\put(40,77){$(01121000)$}
\put(168,22){\line(1,1){53}}
\put(278,22){\line(-1,1){53}}
\put(223,77){\circle{6}}
\put(154,83){$(01111110)$}
\put(112,77){\line(1,0){108}}
\put(139,65){\tiny $\mathbf{(01111100)}$}
\put(282,22){\line(1,1){53}}
\put(226,77){\line(1,0){108}}
\put(253,65){\tiny $\mathbf{(00111110)}$}
\put(337,77){\circle{6}}
\put(344,77){$(00111111)$}
\put(225,79){\line(1,1){53}}
\put(200,108){\tiny $(00011111)$}
\put(335,79){\line(-1,1){53}}
\put(313,108){\tiny $(01011110)$}
\put(280,134){\circle{6}}
\put(253,140){$(01011111)$}
\end{picture}

Write $\alpha_1=(11121000)$, $\alpha_2=(01112100)$, $\alpha_3=(01111110)$, $\alpha_4=(01011111)$, $\alpha_5=(00111111)$, $\alpha_6=(10111110)$, $\alpha_7=(11111100)$; we see that $(\alpha_1,\dots,\alpha_7)$ is a reduced $7$-cycle, and the associated $\beta_i$ in $\Gamma_5$ are $\beta_1=(11111000)$, $\beta_2=(01121000)$, $\beta_3=(01111100)$, $\beta_4=(01011110)$, $\beta_5=(00011111)$, $\beta_6=(00111110)$, $\beta_7=(10111100))$. We obtain:

\[\alpha_1=\beta_1+\delta_4=\beta_2+\delta_1;\]
\[\alpha_2=\beta_2+\delta_6=\beta_3+\delta_4;\]
\[\alpha_3=\beta_3+\delta_7=\beta_4+\delta_3=\beta_6+\delta_2;\]
\[\alpha_4=\beta_4+\delta_8=\beta_5+\delta_2;\]
\[\alpha_5=\beta_5+\delta_3=\beta_6+\delta_8;\]
\[\alpha_6=\beta_6+\delta_1=\beta_7+\delta_7;\]
\[\alpha_7=\beta_7+\delta_2=\beta_1+\delta_6=\beta_3+\delta_1,\]
and we obtain that the following matrix:
\[\left(\begin{array}{ccccccc}&&&\varepsilon{\alpha_1,-\beta_1}x_1&\varepsilon_{\alpha_7,-\beta_1}x_7\\\varepsilon_{\alpha_1,-\beta_2}x_1&&&&\varepsilon_{\alpha_2,-\beta_2}x_2\\\varepsilon_{\alpha_7,-\beta_3}x_7&&&\varepsilon_{\alpha_2,-\beta_3}x_2&&\varepsilon_{\alpha_3,-\beta_3}x_3\\&&\varepsilon_{\alpha_3,-\beta_4}x_3&&&&\varepsilon_{\alpha_4,-\beta_4}x_4\\&\varepsilon_{\alpha_4,-\beta_5}x_4&\varepsilon_{\alpha_5,-\beta_5}x_5\\\varepsilon_{\alpha_6,-\beta_6}x_6&\varepsilon_{\alpha_3,-\beta_6}x_3&&&&&\varepsilon_{\alpha_5,-\beta_6}x_5\\&\varepsilon_{\alpha_7,-\beta_7}x_7&&&&\varepsilon_{\alpha_6,-\beta_7}x_6\end{array}\right)\]
must be invertible for every $x_1,\dots,x_7\in\overline{k}^*$. The determinant of that matrix is $\pm 5x_1x_2x_3x_4x_5x_6x_7$, which is nonzero since we have assumed $p\neq 5$ in the case of a group of type $E_8$, which proes the claim in that case.

Assume now $\underline{\Phi}$ is of type $E_8$ and $l=5$. The graph is the following one:

\begin{picture}(400,160)(0,0)
\put(116,20){\circle{6}}
\put(89,7){$(10111100)$}
\put(119,20){\line(1,0){108}}
\put(146,36){\tiny $\mathbf{(00111100})$}
\put(230,20){\circle{6}}
\put(203,7){$(00111110)$}
\put(114,22){\line(-1,1){53}}
\put(36,45){\tiny $(10111000)$}
\put(59,77){\circle{6}}
\put(-10,83){$(11111000)$}
\put(118,22){\line(1,1){53}}
\put(228,22){\line(-1,1){53}}
\put(173,77){\circle{6}}
\put(180,83){$(0111100)$}
\put(62,77){\line(1,0){108}}
\put(94,93){\tiny $\mathbf{(01111000)}$}
\put(61,79){\line(1,1){53}}
\put(171,79){\line(-1,1){53}}
\put(116,134){\circle{6}}
\put(89,140){$(01121000)$}
\put(56,77){\line(-1,0){60}}
\put(0,64){\tiny $((11110000))$}
\put(176,77){\line(1,0){108}}
\put(203,64){\tiny $(01011100)$}
\put(232,22){\line(1,1){53}}
\put(287,77){\circle{6}}
\put(294,77){$(01011110)$}
\put(233,20){\line(1,0){108}}
\put(260,36){\tiny $\mathbf{(00011110})$}
\put(342,22){\line(-1,1){53}}
\put(344,20){\circle{6}}
\put(319,7){$(00011111)$}
\put(346,22){\line(1,1){30}}
\put(341,56){\tiny $((00001111))$}
\end{picture}

The graph $\mathcal{G}_6$ contains the nonreduced $7$-cycle of the $\beta_i$ stated above. As usual when dealing with nonreduced cycles we only have to check the case $(C)$. Since the element $(11110000)$ (resp. $(00001111)$ of $\Gamma_4$ is related to $\beta_1$ (resp. $\beta_5$) and no other elements of $\Gamma_5$, $\mathcal{S}$ cannot contain these two elements, and since there is an edge between $\beta_1$ and $\beta_7$ which is not part of a nonreduced $3$-cycle, $\mathcal{S}$ cannot contain $\beta_7$ either. The only remaining possibility is that $\mathcal{S}$ contains all four remaining elements of $\Gamma_5$; this subgraph contains the reduced $3$-cycle $(\beta_3,\beta_4,\beta_6)$ and no larger cycle, and we conclude as in the case of reduced $3$-cycles.

Now we will consider the case $l=1$.
Since $u$ doesn't belong to $\underline{\mth{B}}_{\underline{\Delta}}$ we must have $v\leq h-2$. Let $\xi$ be any element of $X_*(\underline{\mth{T}})$, and write $u=\prod_{\delta\in\underline{\Delta}}u_\delta(\varpi^{v}x_\delta)$; we obtain, for every $y\in{\got{p}}^{h-1-v}$:
\[[\xi(1+y),u]=u\prod_{\delta\in\underline{\Delta}}u_\delta(<\delta,\xi>\varpi^{h-1}x_\delta y).\]
Obviously $\xi(1+y)$ is an element of $R_u(\underline{\mth{B}})$; moreover, if there exists some $\xi$ such that for some $\delta$ such that $v(u_\delta)=v$, $<\delta,\xi>$ is not a multiple of $p$, and $<\delta',\xi>=0$ for every $\delta'\neq\delta$ satisfying the same condition, then we obtain $[Im(\xi),u]\supset\underline{\mth{U}}_{\delta,h-1}$, which proves the assertion. If $p$ does not divide the adjoint index of $\underline{\mth{G}}$, such a $\xi$ always exists; that's precisely the assumption we have made in the $A_n$ case, and the hypotheses on $p$ always imply it in the other cases too.

Assume now $w$ is nontrivial; there exists then $\alpha\in\underline{\Delta}$ such that $w(\alpha)>0$, hence $\underline{\mth{U}}_{w(\alpha),h-1}\subset R_u(\underline{\mth{B}})$. Moreover, we have $n^{-1}\underline{\mth{U}}_{w(\alpha),h-1}n=\underline{\mth{U}}_{\alpha,h-1}$, which commutes with $u$ by the commutator relations; hence:
\[R_u(\underline{\mth{B}})nu\underline{\mth{B}}=R_u(\underline{\mth{B}})nu\underline{\mth{U}}_{\alpha,h-1}\underline{\mth{B}}=R_u(\underline{\mth{B}})nu\underline{\mth{B}}_\alpha.\]
Since $\alpha$ doesn't depend on $u$, it depends only on the double class of $g$ modulo $\underline{\mth{B}}_{\underline{\Delta}}$. Hence the proposition is proved in the irreducible case.

Assume now $\underline{\Phi}$ is not irreducible. Let $\underline{\Phi}_1,\dots,\underline{\Phi}_s$ be its irreducible components. For every $i$, the subset $\underline{\Phi}_i^-=\underline{\Phi}_i\cap\underline{\Phi}^-$ of $\underline{\Phi}^-$ is $f$-closed and complete, and we can thus apply the corollary to proposition \ref{dcprod} to the $f$-closed and complete set $\underline{\Phi}^-$ which is their disjoint union.

Now let $n$ be an element of the normalizzer of $\underline{\mth{T}}$ and $u$ be an element of $\underline{\mth{U}}^-$; we have $nu=\prod_{i=1}^sn_iu_i$, where for every $i$ $n_i$ is a representative of an element of the Weyl group of $\underline{\Psi}_i=\underline{\Psi}\cap\underline{\Phi}_i$ and $u_i$ is an element of the group $\underline{\mth{U}}_i$ generated by the $\underline{\mth{U}}_\alpha$, $\alpha\in\underline{\Phi}_i$. Assume we have $R_u(\underline{\mth{B}})n_iu_i\underline{\mth{B}}=R_u(\underline{\mth{B}})n_iu_i\underline{\mth{B}}_\alpha$ for some $i$ and some negative simple root $\alpha$ in $\underline{\Phi}_i$; we then also have $\bigsqcup_{u'}R_u(\underline{\mth{B}})nu'\underline{\mth{B}}=\bigcup_{u'}R_u(\underline{\mth{B}})nu'\underline{\mth{B}}_\alpha$, where $u'$ runs over a set of representatives of the $R_u(\underline{\mth{H}})$-conjugacy classes of $\underline{\mth{U}}^-$ of the form $u_1t_1u_2\dots t_{s-1}u_st_s$, with the $t_s$ being elements of $R_u(\underline{\mth{H}})$ such that $\prod_{i=1}^st_i=1$. Moreover, since $u$ doesn't belong to $\underline{\mth{B}}_{\underline{\Delta}}$, then neither do the $u'$, and we deduce from this that if $u'$ and $u''$ are not conjugated by any element of $R_u(\underline{\mth{H}})$, for every $u_\alpha\in\underline{\mth{U}}_{\alpha,h-1}$, $u'u_\alpha$ and $u''$ don't belong to the same element of $R_u(\underline{\mth{B}}\backslash\underline{\mth{G}}/\underline{\mth{B}}$ either; hence the second union is disjoint too, and in particular we have $R_u(\underline{\mth{B}})nu\underline{\mth{B}}=R_u(\underline{\mth{B}})nu\underline{\mth{B}}_\alpha$. We are thus reduced to the case when $\underline{\Phi}$ is irreducible. $\Box$

\end{document}